\documentclass{article}
\usepackage{amsmath}
\usepackage{amssymb}
\begin{document}
\def\e#1\e{\begin{equation}#1\end{equation}}
\def\ea#1\ea{\begin{align}#1\end{align}}
\def\eq#1{{\rm(\ref{#1})}}
\newtheorem{thm}{Theorem}[section]
\newtheorem{lem}[thm]{Lemma}
\newtheorem{prop}[thm]{Proposition}
\newtheorem{cor}[thm]{Corollary}
\newenvironment{dfn}{\medskip\refstepcounter{thm}
\noindent{\bf Definition \thesection.\arabic{thm}\ }}{\medskip}
\newenvironment{ex}{\medskip\refstepcounter{thm}
\noindent{\bf Example \thesection.\arabic{thm}\ }}{\medskip}
\newenvironment{proof}[1][,]{\medskip\ifcat,#1
\noindent{\it Proof.\ }\else\noindent{\it Proof of #1.\ }\fi}
{\hfill$\square$\medskip}
\def\dim{\mathop{\rm dim}}
\def\Re{\mathop{\rm Re}}
\def\Im{\mathop{\rm Im}}
\def\Ker{\mathop{\rm Ker}}
\def\vol{\mathop{\rm vol}}
\def\sind{{\ts\mathop{\text{\rm s-ind}}}}
\def\id{\mathop{\rm id}}
\def\U{\mathbin{\rm U}}
\def\SU{\mathop{\rm SU}}
\def\ge{\geqslant} 
\def\le{\leqslant} 
\def\N{\mathbin{\mathbb N}}
\def\R{\mathbin{\mathbb R}}
\def\Z{\mathbin{\mathbb Z}}
\def\C{\mathbin{\mathbb C}}
\def\g{\mathfrak{g}} 
\def\su{\mathfrak{su}} 
\def\al{\alpha}
\def\be{\beta}
\def\ga{\gamma}
\def\de{\delta}
\def\ep{\epsilon}
\def\ve{\varepsilon}
\def\la{\lambda}
\def\ka{\kappa}
\def\th{\theta}
\def\ze{\zeta}
\def\up{\upsilon}
\def\vp{\varphi}
\def\si{\sigma}
\def\om{\omega}
\def\De{\Delta}
\def\La{\Lambda}
\def\Om{\Omega}
\def\Ga{\Gamma}
\def\Si{\Sigma}
\def\Th{\Theta}
\def\Up{\Upsilon}
\def\d{{\rm d}}
\def\pd{\partial}
\def\db{{\bar\partial}}
\def\ts{\textstyle}
\def\sst{\scriptscriptstyle}
\def\w{\wedge}
\def\lt{\ltimes}
\def\sm{\setminus}
\def\op{\oplus}
\def\ot{\otimes}
\def\bigot{\bigotimes}
\def\iy{\infty}
\def\ra{\rightarrow}
\def\hookra{\hookrightarrow}
\def\longra{\longrightarrow}
\def\t{\times}
\def\na{\nabla}
\def\ha{{\textstyle\frac{1}{2}}}
\def\ti{\tilde}
\def\wt{\widetilde}
\def\ov{\overline}
\def\ovB{\,\overline{\!B}}
\def\sC{{\smash{\sst C}}}
\def\sCi{{\smash{\sst C_i}}}
\def\sE{{\smash{\sst\cal E}}}
\def\sL{{\smash{\sst L}}}
\def\sLi{{\smash{\sst L_i}}}
\def\sSi{{\smash{\sst\Si}}}
\def\sSii{{\smash{\sst\Si_i}}}
\def\sN{{\smash{\sst N}}}
\def\sX{{\smash{\sst X}}}
\def\sNt{{\smash{\sst N^t}}}
\def\sNst{{\smash{\sst N^{s,t}}}}
\def\sW{{\smash{\sst W}}}
\def\sWt{{\smash{\sst W^t}}}
\def\sWst{{\smash{\sst W^{s,t}}}}
\def\sXp{{\smash{\sst X'}}}
\def\sF{{\smash{\sst\cal F}}}
\def\sFp{{\smash{\sst{\cal F}'}}}
\def\B{{\cal B}}
\def\D{{\cal D}}
\def\E{{\cal E}}
\def\F{{\cal F}}
\def\G{{\cal G}}
\def\H{{\cal H}}
\def\I{{\cal I}}
\def\K{{\cal K}}
\def\M{{\cal M}}
\def\O{{\cal O}}
\def\V{{\cal V}}
\def\ms#1{\vert#1\vert^2}
\def\bms#1{\bigl\vert#1\bigr\vert^2}
\def\md#1{\vert #1 \vert}
\def\bmd#1{\big\vert #1 \big\vert}
\def\cnm#1#2{\Vert #1 \Vert_{C^{#2}}} 
\def\lnm#1#2{\Vert #1 \Vert_{L^{#2}}} 
\def\snm#1#2#3{\Vert #1 \Vert_{L^{#2}_{#3}}} 
\def\blnm#1#2{\bigl\Vert #1 \bigr\Vert_{L^{#2}}} 
\def\an#1{\langle#1\rangle}
\def\ban#1{\bigl\langle#1\bigr\rangle}
\title{Special Lagrangian submanifolds with isolated\\
conical singularities. IV. Desingularization, \\
obstructions and families}
\author{Dominic Joyce \\ Lincoln College, Oxford}
\date{}
\maketitle

\section{Introduction}
\label{co1}

{\it Special Lagrangian $m$-folds (SL\/ $m$-folds)} are a
distinguished class of real $m$-dimensional minimal submanifolds
which may be defined in $\C^m$, or in {\it Calabi--Yau $m$-folds},
or more generally in {\it almost Calabi--Yau $m$-folds} (compact
K\"ahler $m$-folds with trivial canonical bundle). We write an
almost Calabi--Yau $m$-fold as $M$ or $(M,J,\om,\Om)$, where
the manifold $M$ has complex structure $J$, K\"ahler form
$\om$ and holomorphic volume form~$\Om$.

This is the fourth in a series of five papers
\cite{Joyc3,Joyc4,Joyc5,Joyc6} studying SL $m$-folds with
{\it isolated conical singularities}. That is, we consider an
SL $m$-fold $X$ in an almost Calabi--Yau $m$-fold $M$ for $m>2$
with singularities at $x_1,\ldots,x_n$ in $M$, such that for some
special Lagrangian cones $C_i$ in $T_{\smash{x_i}}M\cong\C^m$
with $C_i\sm\{0\}$ nonsingular, $X$ approaches $C_i$ near $x_i$
in an asymptotic $C^1$ sense. Readers are advised to begin with
the final paper \cite{Joyc6}, which surveys the series, and
applies the results to prove some conjectures.

The first paper \cite{Joyc3} laid the foundations for the series,
and studied the {\it regularity} of SL $m$-folds with conical
singularities near their singular points. The second paper
\cite{Joyc4} discussed the {\it deformation theory} of compact
SL $m$-folds $X$ with conical singularities in an almost
Calabi--Yau $m$-fold~$M$.

The third paper \cite{Joyc5} and this one study
{\it desingularizations} of compact SL $m$-folds $X$ with
conical singularities. That is, we construct a family of
compact, {\it nonsingular} SL $m$-folds $\smash{\ti N^t}$
in $M$ for $t\in(0,\ep]$ with $\smash{\ti N^t}\ra X$
as $t\ra 0$, in the sense of currents. In \cite{Joyc5} we
did this for simple situations, working in a single almost
Calabi--Yau $m$-fold $(M,J,\om,\Om)$, and making topological
assumptions to avoid problems with obstructions to the
existence of~$\smash{\ti N^t}$.

This paper extends the results of \cite{Joyc5} to more
complicated situations, in which there are {\it topological
obstructions} to the existence of desingularizations
$\smash{\ti N^t}$ of $X$, and to desingularizations in
{\it families} of almost Calabi--Yau $m$-folds
$(M,J^s,\om^s,\Om^s)$ for $s\in\F$.

Having a good understanding of the singularities of special
Lagrangian submanifolds will be essential in clarifying the
Strominger--Yau--Zaslow conjecture on the Mirror Symmetry
of Calabi--Yau 3-folds \cite{SYZ}, and also in resolving
conjectures made by the author \cite{Joyc1} on defining
new invariants of Calabi--Yau 3-folds by counting special
Lagrangian homology 3-spheres with weights. The series aims
to develop such an understanding for simple singularities
of SL $m$-folds.

Here is the basic idea of \cite{Joyc5}. Let $X$ be a compact
SL $m$-fold with conical singularities $x_1,\ldots,x_n$ in an
almost Calabi--Yau $m$-fold $(M,J,\om,\Om)$. Choose an
isomorphism $\up_i:\C^m\ra T_{x_i}M$ for $i=1,\ldots,n$.
Then there is a unique {\it SL cone} $C_i$ in $\C^m$ with
$X$ asymptotic to $\up_i(C_i)$ at $x_i$. Let $L_i$ be an
{\it Asymptotically Conical SL\/ $m$-fold\/} ({\it AC SL\/
$m$-fold\/}) in $\C^m$, asymptotic to $C_i$ at infinity.

Now $tL_i=\{t\,{\bf x}:{\bf x}\in L_i\}$ is also an
AC SL $m$-fold asymptotic to $C_i$ for $t>0$. We construct a
1-parameter family of compact, nonsingular {\it Lagrangian}
$m$-folds $N^t$ in $(M,\om)$ for $t\in(0,\de)$ by gluing
$tL_i$ into $X$ at $x_i$. When $t$ is small, $N^t$ is close
to being special Lagrangian, but also close to being singular.
We prove using analysis that for small $t$ we can deform
$N^t$ to a {\it special\/} Lagrangian $m$-fold $\smash{\ti
N^t}$ in $M$, using a small Hamiltonian deformation.

In this paper we shall study the following issues, not tackled
in \cite{Joyc5}. The AC SL $m$-folds $L_i$ have topological
invariants $Y(L_i),Z(L_i)$ defined in \S\ref{co41},
which measure the relative de Rham cohomology classes
of $\om$ and $\Im\Om$ in $H^k(\C^m,L_i;\R)$. In
\cite{Joyc5} we assumed that $Y(L_i)=0$ for $i=1,\ldots,n$.
Section \ref{co6} shows how to extend the results of
\cite{Joyc5} to the case $Y(L_i)\ne 0$, so that they
are applicable to a much larger class of AC SL $m$-folds.

Doing this is a problem in symplectic geometry. If $Y(L_i)=0$
then we can choose the Lagrangian $m$-folds $N^t$ to coincide
with $X$ away from $x_i$, and work locally near $x_i$. But if
$Y(L_i)\ne 0$ we cannot do this. Instead, we must define
$N^t$ away from $x_i$ as the graph of a closed 1-form on
$X'=X\sm\{x_1,\ldots,x_n\}$ with nonzero cohomology class,
and there can be {\it topological obstructions} to the
existence of $N^t$ as a Lagrangian $m$-fold.

In \S\ref{co7} and \S\ref{co8} we extend the results to
{\it smooth families} of almost Calabi--Yau $m$-folds
$(M,J^s,\om^s,\Om^s)$, for $s\in\F\subset\R^d$ with $0\in\F$
and $(M,J,\om,\Om)=(M,J^0,\om^0,\Om^0)$. It turns out that
the cohomology classes $[\om^s]$ and $[\Im\Om^s]$ contribute
to the obstruction equations involving $Y(L_i)$ and $Z(L_i)$
for the existence of desingularizations~$\smash{\ti N^t}$.

Because of this, it can happen that a singular SL $m$-fold
$X$ in $(M,J,\om,\Om)$ admits no desingularizations
$\smash{\ti N^t}$ in $(M,J,\om,\Om)$, but does admit
desingularizations $\smash{\ti N^{s,t}}$ in $(M,J^s,\om^s,\Om^s)$
for small $s\ne 0$. Thus we can overcome obstructions to the
existence of desingularizations by varying the underlying
almost Calabi--Yau $m$-fold~$(M,J,\om,\Om)$.

We begin in \S\ref{co2} with an introduction to special
Lagrangian geometry. Sections \ref{co3} and \ref{co4}
define SL $m$-folds with conical singularities and
Asymptotically Conical SL $m$-folds, and review results
we need from \cite{Joyc3}. Section \ref{co5} recalls and
discusses the major definitions and theorems from the
previous paper~\cite{Joyc5}.

The new material of the paper is \S\ref{co6}--\S\ref{co8}.
Section \ref{co6} generalizes the results of \cite{Joyc5}
to the case when $Y(L_i)\ne 0$, and \S\ref{co7} to families
of almost Calabi--Yau $m$-folds $(M,J^s,\om^s,\Om^s)$ when
$Y(L_i)=0$. Finally \S\ref{co8} considers the most complicated
case, in families $(M,J^s,\om^s,\Om^s)$ when~$Y(L_i)\ne 0$.

For simplicity we generally take all submanifolds to be {\it
embedded}. However, all our results generalize immediately
to {\it immersed\/} submanifolds, with only cosmetic changes.

Other authors have also desingularized SL $m$-folds
using gluing techniques. Those known to me are Salur
\cite{Salu1,Salu2}, Butscher \cite{Buts} and Lee \cite{Lee},
which were discussed in \cite[\S 1]{Joyc5}. They all involve
connect sum constructions in Calabi--Yau $m$-folds or $\C^m$,
rather than the more general kinds of singularities we consider.
\medskip

\noindent{\it Acknowledgements.} I would like to thank
Stephen Marshall, Sema Salur and Adrian Butscher for useful
conversations. I was supported by an EPSRC Advanced Research
Fellowship whilst writing this paper.

\section{Special Lagrangian geometry}
\label{co2}

We introduce special Lagrangian submanifolds (SL $m$-folds)
in two different geometric contexts. First, in \S\ref{co21},
we define SL $m$-folds in $\C^m$. Then \S\ref{co22} discusses
SL $m$-folds in {\it almost Calabi--Yau $m$-folds}, compact
K\"ahler manifolds equipped with a holomorphic volume form,
which generalize Calabi--Yau manifolds. Some references for
this section are Harvey and Lawson \cite{HaLa} and the
author~\cite{Joyc2}.

\subsection{Special Lagrangian submanifolds in $\C^m$}
\label{co21}

We begin by defining {\it calibrations} and {\it calibrated 
submanifolds}, following~\cite{HaLa}.

\begin{dfn} Let $(M,g)$ be a Riemannian manifold. An {\it oriented
tangent $k$-plane} $V$ on $M$ is a vector subspace $V$ of
some tangent space $T_xM$ to $M$ with $\dim V=k$, equipped
with an orientation. If $V$ is an oriented tangent $k$-plane
on $M$ then $g\vert_V$ is a Euclidean metric on $V$, so 
combining $g\vert_V$ with the orientation on $V$ gives a 
natural {\it volume form} $\vol_V$ on $V$, which is a 
$k$-form on~$V$.

Now let $\vp$ be a closed $k$-form on $M$. We say that
$\vp$ is a {\it calibration} on $M$ if for every oriented
$k$-plane $V$ on $M$ we have $\vp\vert_V\le \vol_V$. Here
$\vp\vert_V=\al\cdot\vol_V$ for some $\al\in\R$, and 
$\vp\vert_V\le\vol_V$ if $\al\le 1$. Let $N$ be an 
oriented submanifold of $M$ with dimension $k$. Then 
each tangent space $T_xN$ for $x\in N$ is an oriented
tangent $k$-plane. We say that $N$ is a {\it calibrated 
submanifold\/} if $\vp\vert_{T_xN}=\vol_{T_xN}$ for all~$x\in N$.
\label{co2def1}
\end{dfn}

It is easy to show that calibrated submanifolds are automatically
{\it minimal submanifolds} \cite[Th.~II.4.2]{HaLa}. Here is the 
definition of special Lagrangian submanifolds in $\C^m$, taken
from~\cite[\S III]{HaLa}.

\begin{dfn} Let $\C^m$ have complex coordinates $(z_1,\dots,z_m)$, 
and define a metric $g'$, a real 2-form $\om'$ and a complex $m$-form 
$\Om'$ on $\C^m$ by
\e
\begin{split}
g'=\ms{\d z_1}+\cdots+\ms{\d z_m},\quad
\om'&=\ts\frac{i}{2}(\d z_1\w\d\bar z_1+\cdots+\d z_m\w\d\bar z_m),\\
\text{and}\quad\Om'&=\d z_1\w\cdots\w\d z_m.
\end{split}
\label{co2eq1}
\e
Then $\Re\Om'$ and $\Im\Om'$ are real $m$-forms on $\C^m$. Let $L$
be an oriented real submanifold of $\C^m$ of real dimension $m$. We
say that $L$ is a {\it special Lagrangian submanifold\/} of $\C^m$,
or {\it SL\/ $m$-fold}\/ for short, if $L$ is calibrated with respect
to $\Re\Om'$, in the sense of Definition~\ref{co2def1}.
\label{co2def2}
\end{dfn}

Harvey and Lawson \cite[Cor.~III.1.11]{HaLa} give the following
alternative characterization of special Lagrangian submanifolds:

\begin{prop} Let\/ $L$ be a real $m$-dimensional submanifold 
of\/ $\C^m$. Then $L$ admits an orientation making it into an
SL submanifold of\/ $\C^m$ if and only if\/ $\om'\vert_L\equiv 0$ 
and\/~$\Im\Om'\vert_L\equiv 0$.
\label{co2prop1}
\end{prop}

An $m$-dimensional submanifold $L$ in $\C^m$ is called {\it Lagrangian} 
if $\om'\vert_L\equiv 0$. Thus special Lagrangian submanifolds are 
Lagrangian submanifolds satisfying the extra condition that 
$\Im\Om'\vert_L\equiv 0$, which is how they get their name.

\subsection{Almost Calabi--Yau $m$-folds and SL $m$-folds} 
\label{co22}

We shall define special Lagrangian submanifolds not just in
Calabi--Yau manifolds, as usual, but in the much larger
class of {\it almost Calabi--Yau manifolds}.

\begin{dfn} Let $m\ge 2$. An {\it almost Calabi--Yau $m$-fold\/}
is a quadruple $(M,J,\om,\Om)$ such that $(M,J)$ is a compact
$m$-dimensional complex manifold, $\om$ is the K\"ahler form
of a K\"ahler metric $g$ on $M$, and $\Om$ is a non-vanishing
holomorphic $(m,0)$-form on~$M$.

We call $(M,J,\om,\Om)$ a {\it Calabi--Yau $m$-fold\/} if in
addition $\om$ and $\Om$ satisfy
\e
\om^m/m!=(-1)^{m(m-1)/2}(i/2)^m\Om\w\bar\Om.
\label{co2eq2}
\e
Then for each $x\in M$ there exists an isomorphism $T_xM\cong\C^m$
that identifies $g_x,\om_x$ and $\Om_x$ with the flat versions
$g',\om',\Om'$ on $\C^m$ in \eq{co2eq1}. Furthermore, $g$ is
Ricci-flat and its holonomy group is a subgroup of~$\SU(m)$.
\label{co2def3}
\end{dfn}

This is not the usual definition of a Calabi--Yau manifold, but
is essentially equivalent to it.

\begin{dfn} Let $(M,J,\om,\Om)$ be an almost Calabi--Yau $m$-fold,
and $N$ a real $m$-dimensional submanifold of $M$. We call $N$ a
{\it special Lagrangian submanifold}, or {\it SL $m$-fold\/} for
short, if $\om\vert_N\equiv\Im\Om\vert_N\equiv 0$. It easily
follows that $\Re\Om\vert_N$ is a nonvanishing $m$-form on $N$.
Thus $N$ is orientable, with a unique orientation in which
$\Re\Om\vert_N$ is positive.
\label{co2def4}
\end{dfn}

Again, this is not the usual definition of SL $m$-fold, but is
essentially equivalent to it. Suppose $(M,J,\om,\Om)$ is an
almost Calabi--Yau $m$-fold, with metric $g$. Let
$\psi:M\ra(0,\iy)$ be the unique smooth function such that
\e
\psi^{2m}\om^m/m!=(-1)^{m(m-1)/2}(i/2)^m\Om\w\bar\Om,
\label{co2eq3}
\e
and define $\ti g$ to be the conformally equivalent metric $\psi^2g$
on $M$. Then $\Re\Om$ is a {\it calibration} on the Riemannian manifold
$(M,\ti g)$, and SL $m$-folds $N$ in $(M,J,\om,\Om)$ are calibrated
with respect to it, so that they are minimal with respect to~$\ti g$.

If $M$ is a Calabi--Yau $m$-fold then $\psi\equiv 1$ by \eq{co2eq2},
so $\ti g=g$, and an $m$-submanifold $N$ in $M$ is special Lagrangian
if and only if it is calibrated w.r.t.\ $\Re\Om$ on $(M,g)$, as in
Definition \ref{co2def2}. This recovers the usual definition of
special Lagrangian $m$-folds in Calabi--Yau $m$-folds.

The {\it deformation theory} of special Lagrangian submanifolds
was studied by McLean \cite[\S 3]{McLe}, who proved the following
result in the Calabi--Yau case. The extension to the almost
Calabi--Yau case is described in~\cite[\S 9.5]{Joyc2}.

\begin{thm} Let\/ $N$ be a compact SL\/ $m$-fold in an almost
Calabi--Yau $m$-fold\/ $(M,J,\om,\Om)$. Then the moduli space
$\M_\sN$ of special Lagrangian deformations of\/ $N$ is a smooth
manifold of dimension $b^1(N)$, the first Betti number of\/~$N$.
\label{co2thm1}
\end{thm}

We shall often consider {\it families} of almost Calabi--Yau
$m$-folds.

\begin{dfn} Let $(M,J,\om,\Om)$ be an almost Calabi--Yau
$m$-fold. A {\it smooth family of deformations of\/}
$(M,J,\om,\Om)$ is a connected open set $\F\subset\R^d$
for $d\ge 0$ with $0\in\F$ called the {\it base space},
and a smooth family $\bigl\{(M,J^s,\om^s,\Om^s):s\in\F
\bigr\}$ of almost Calabi--Yau structures on $M$
with~$(J^0,\om^0,\Om^0)=(J,\om,\Om)$.
\label{co2def5}
\end{dfn}

If $N$ is an SL $m$-fold in $(M,J,\om,\Om)$, the moduli
spaces of deformations of $N$ in each $(M,J^s,\om^s,\Om^s)$
for $s\in\F$ fit together into a big moduli space~$\M_\sN^\sF$.

\begin{dfn} Let $\bigl\{(M,J^s,\om^s,\Om^s):s\in\F\bigr\}$
be a smooth family of deformations of an almost Calabi--Yau
$m$-fold $(M,J,\om,\Om)$, and $N$ be a compact SL $m$-fold
in $(M,J,\om,\Om)$. Define the {\it moduli space $\M_\sN^\sF$
of deformations of\/ $N$ in the family} $\F$ to be the set
of pairs $(s,\hat N)$ for which $s\in\F$ and $\hat N$ is a
compact SL $m$-fold in $(M,J^s,\om^s,\Om^s)$ which is
diffeomorphic to $N$ and isotopic to $N$ in $M$. Define a
{\it projection} $\pi^\sF:\M_\sN^\sF\ra\F$ by
$\pi^\sF(s,\hat N)=s$. Then $\M_\sN^\sF$ has a
natural topology, and $\pi^\sF$ is continuous.
\label{co2def6}
\end{dfn}

The following result is proved by Marshall \cite[Th.~3.2.9]{Mars},
using similar methods to Theorem~\ref{co2thm1}.

\begin{thm} Let\/ $\bigl\{(M,J^s,\om^s,\Om^s):s\in\F\bigr\}$ be
a smooth family of deformations of an almost Calabi--Yau $m$-fold\/
$(M,J,\om,\Om)$, with base space $\F\subset\R^d$. Suppose $N$ is a
compact SL\/ $m$-fold in $(M,J,\om,\Om)$ with\/ $[\om^s\vert_N]=0$
in $H^2(N,\R)$ and\/ $[\Im\Om^s\vert_N]=0$ in $H^m(N,\R)$ for all\/
$s\in\F$. Let\/ $\M_\sN^\sF$ be the moduli space of deformations
of\/ $N$ in $\F$, and\/ $\pi^\sF:\M_\sN^\sF\ra\F$ the natural
projection.

Then $\M_\sN^\sF$ is a smooth manifold of dimension\/ $d+b^1(N)$,
and\/ $\pi^\sF:\M_\sN^\sF\ra\F$ a smooth submersion. For small
$s\in\F$ the moduli space $\M_\sN^s=(\pi^\sF)^{-1}(s)$ of
deformations of\/ $N$ in $(M,J^s,\om^s,\Om^s)$ is a nonempty
smooth manifold of dimension $b^1(N)$, with\/~$\M_\sN^0=\M_\sN$.
\label{co2thm2}
\end{thm}

This describes the {\it obstructions} to the existence of SL $m$-folds
when we deform the underlying almost Calabi--Yau $m$-fold.

\section{SL $m$-folds with conical singularities}
\label{co3}

The preceding papers \cite{Joyc3,Joyc4,Joyc5} defined and studied
{\it compact SL\/ $m$-folds $X$ with conical singularities} in
an almost Calabi--Yau $m$-fold $(M,J,\om,\Om)$. We now recall
the definitions and results from \cite{Joyc3} that we will need
later. For brevity we shall keep explanations to a minimum, and
readers are referred to \cite{Joyc3} for further discussion and
motivation.

\subsection{Preliminaries on special Lagrangian cones}
\label{co31}

Following \cite[\S 2.1]{Joyc3} we give definitions
and results on {\it special Lagrangian cones}.

\begin{dfn} A (singular) SL $m$-fold $C$ in $\C^m$ is called a
{\it cone} if $C=tC$ for all $t>0$, where $tC=\{t\,{\bf x}:{\bf x}
\in C\}$. Let $C$ be a closed SL cone in $\C^m$ with an isolated
singularity at 0. Then $\Si=C\cap{\cal S}^{2m-1}$ is a compact,
nonsingular $(m\!-\!1)$-submanifold of ${\cal S}^{2m-1}$, not
necessarily connected. Let $g_\sSi$ be the restriction
of $g'$ to $\Si$, where $g'$ is as in~\eq{co2eq1}.

Set $C'=C\sm\{0\}$. Define $\iota:\Si\t(0,\iy)\ra\C^m$ by
$\iota(\si,r)=r\si$. Then $\iota$ has image $C'$. By an abuse
of notation, {\it identify} $C'$ with $\Si\t(0,\iy)$ using
$\iota$. The {\it cone metric} on $C'\cong\Si\t(0,\iy)$
is~$g'=\iota^*(g')=\d r^2+r^2g_\sSi$.

For $\al\in\R$, we say that a function $u:C'\ra\R$ is
{\it homogeneous of order} $\al$ if $u\circ t\equiv t^\al u$ for
all $t>0$. Equivalently, $u$ is homogeneous of order $\al$ if
$u(\si,r)\equiv r^\al v(\si)$ for some function~$v:\Si\ra\R$.
\label{co3def1}
\end{dfn}

In \cite[Lem.~2.3]{Joyc3} we study {\it homogeneous harmonic
functions} on~$C'$.

\begin{lem} In the situation of Definition \ref{co3def1},
let\/ $u(\si,r)\equiv r^\al v(\si)$ be a homogeneous function
of order $\al$ on $C'=\Si\t(0,\iy)$, for $v\in C^2(\Si)$. Then
\begin{equation*}
\De u(\si,r)=r^{\al-2}\bigl(\De_\sSi v-\al(\al+m-2)v\bigr),
\end{equation*}
where $\De$, $\De_\sSi$ are the Laplacians on $(C',g')$
and\/ $(\Si,g_\sSi)$. Hence, $u$ is harmonic on $C'$
if and only if\/ $v$ is an eigenfunction of\/ $\De_\sSi$
with eigenvalue~$\al(\al+m-2)$.
\label{co3lem}
\end{lem}

Following \cite[Def.~2.5]{Joyc3}, we define:

\begin{dfn} In the situation of Definition \ref{co3def1},
suppose $m>2$ and define
\e
\D_\sSi=\bigl\{\al\in\R:\text{$\al(\al+m-2)$ is
an eigenvalue of $\De_\sSi$}\bigr\}.
\label{co3eq1}
\e
Then $\D_\sSi$ is a countable, discrete subset of
$\R$. By Lemma \ref{co3lem}, an equivalent definition is that
$\D_\sSi$ is the set of $\al\in\R$ for which there
exists a nonzero homogeneous harmonic function $u$ of order
$\al$ on~$C'$.

Define $m_\sSi:\D_\sSi\ra\N$ by taking
$m_\sSi(\al)$ to be the multiplicity of the eigenvalue
$\al(\al+m-2)$ of $\De_\sSi$, or equivalently the
dimension of the vector space of homogeneous harmonic
functions $u$ of order $\al$ on $C'$. Define
$N_\sSi:\R\ra\Z$ by
\begin{equation*}
N_\sSi(\de)=
-\sum_{\!\!\!\!\al\in\D_\sSi\cap(\de,0)\!\!\!\!}m_\sSi(\al)
\;\>\text{if $\de<0$, and}\;\>
N_\sSi(\de)=
\sum_{\!\!\!\!\al\in\D_\sSi\cap[0,\de]\!\!\!\!}m_\sSi(\al)
\;\>\text{if $\de\ge 0$.}
\end{equation*}
Then $N_\sSi$ is monotone increasing and upper semicontinuous,
and is discontinuous exactly on $\D_\sSi$, increasing by
$m_\sSi(\al)$ at each $\al\in\D_\sSi$. As the
eigenvalues of $\De_\sSi$ are nonnegative, we see that
$\D_\sSi\cap(2-m,0)=\emptyset$ and $N_\sSi\equiv 0$
on~$(2-m,0)$.
\label{co3def2}
\end{dfn}

\subsection{The definition of SL $m$-folds with conical singularities}
\label{co32}

Now we can define {\it conical singularities} of SL $m$-folds,
following~\cite[Def.~3.6]{Joyc3}.

\begin{dfn} Let $(M,J,\om,\Om)$ be an almost Calabi--Yau $m$-fold
for $m>2$, and define $\psi:M\ra(0,\iy)$ as in \eq{co2eq3}. Suppose
$X$ is a compact singular SL $m$-fold in $M$ with singularities at
distinct points $x_1,\ldots,x_n\in X$, and no other singularities.

Fix isomorphisms $\up_i:\C^m\ra T_{x_i}M$ for $i=1,\ldots,n$
such that $\up_i^*(\om)=\om'$ and $\up_i^*(\Om)=\psi(x_i)^m\Om'$,
where $\om',\Om'$ are as in \eq{co2eq1}. Let $C_1,\ldots,C_n$ be SL
cones in $\C^m$ with isolated singularities at 0. For $i=1,\ldots,n$
let $\Si_i=C_i\cap{\cal S}^{2m-1}$, and let $\mu_i\in(2,3)$ with
$(2,\mu_i]\cap\D_\sSii=\emptyset$, where $\D_\sSii$ is defined in
\eq{co3eq1}. Then we say that $X$ has a {\it conical singularity}
at $x_i$, with {\it rate} $\mu_i$ and {\it cone} $C_i$ for
$i=1,\ldots,n$, if the following holds.

By Darboux' Theorem \cite[Th.~3.15]{McSa} there exist embeddings
$\Up_i:B_R\ra M$ for $i=1,\ldots,n$ satisfying $\Up_i(0)=x_i$,
$\d\Up_i\vert_0=\up_i$ and $\Up_i^*(\om)=\om'$, where $B_R$ is
the open ball of radius $R$ about 0 in $\C^m$ for some small $R>0$.
Define $\iota_i:\Si_i\t(0,R)\ra B_R$ by $\iota_i(\si,r)=r\si$
for~$i=1,\ldots,n$.

Define $X'=X\sm\{x_1,\ldots,x_n\}$. Then there should exist a
compact subset $K\subset X'$ such that $X'\sm K$ is a union of
open sets $S_1,\ldots,S_n$ with $S_i\subset\Up_i(B_R)$, whose
closures $\bar S_1,\ldots,\bar S_n$ are disjoint in $X$. For
$i=1,\ldots,n$ and some $R'\in(0,R]$ there should exist a smooth
$\phi_i:\Si_i\t(0,R')\ra B_R$ such that $\Up_i\circ\phi_i:\Si_i
\t(0,R')\ra M$ is a diffeomorphism $\Si_i\t(0,R')\ra S_i$, and
\e
\bmd{\na^k(\phi_i-\iota_i)}=O(r^{\mu_i-1-k})
\quad\text{as $r\ra 0$ for $k=0,1$.}
\label{co3eq2}
\e
Here $\na,\md{\,.\,}$ are computed using the cone metric
$\iota_i^*(g')$ on~$\Si_i\t(0,R')$.
\label{co3def3}
\end{dfn}

The reasoning behind this definition was discussed in
\cite[\S 3.3]{Joyc3}. We suppose $m>2$ for two reasons.
Firstly, the only SL cones in $\C^2$ are finite unions
of SL planes $\R^2$ in $\C^2$ intersecting only at 0.
Thus any SL 2-fold with conical singularities is actually
{\it nonsingular} as an immersed 2-fold, so there is
really no point in studying them.

Secondly, $m=2$ is a special case in the analysis of
\cite[\S 2]{Joyc3}, and it is simpler to exclude it.
Therefore we will suppose $m>2$ throughout the paper.
We will need the following tool \cite[Def.~2.6]{Joyc3},
a smoothed out version of the distance from the singular
set $\{x_1,\ldots,x_n\}$ in~$X$.

\begin{dfn} Let $(M,J,\om,\Om)$ be an almost Calabi--Yau
$m$-fold and $X$ a compact SL $m$-fold in $M$ with conical
singularities at $x_1,\ldots,x_n$, and use the notation of
Definition \ref{co3def3}. Define a {\it radius function}
$\rho$ on $X'$ to be a smooth function $\rho:X'\ra(0,1]$
such that $\rho\equiv 1$ on $K$ and $\rho(y)=d(x_i,y)$ for
$y\in S_i$ close to $x_i$, where $d$ is the metric on $X$.
Radius functions always exist.
\label{co3def4}
\end{dfn}

\subsection{Homology, cohomology and Hodge theory}
\label{co33}

Next we discuss {\it homology} and {\it cohomology} of SL $m$-folds
with conical singularities, following \cite[\S 2.4]{Joyc3}. For a
general reference, see for instance Bredon \cite{Bred}. When $Y$
is a manifold, write $H^k(Y,\R)$ for the $k^{\rm th}$ {\it de Rham
cohomology group} and $H^k_{\rm cs}(Y,\R)$ for the $k^{\rm th}$
{\it compactly-supported de Rham cohomology group} of $Y$. If $Y$
is compact then $H^k(Y,\R)=H^k_{\rm cs}(Y,\R)$. The {\it Betti
numbers} of $Y$ are $b^k(Y)=\dim H^k(Y,\R)$ and~$b^k_{\rm cs}(Y)
=\dim H^k_{\rm cs}(Y,\R)$.

Let $Y$ be a topological space, and $Z\subset Y$ a subspace.
Write $H_k(Y,\R)$ for the $k^{\rm th}$ {\it real singular
homology group} of $Y$, and $H_k(Y;Z,\R)$ for the $k^{\rm th}$
{\it real singular relative homology group} of $(Y;Z)$. When
$Y$ is a manifold and $Z$ a submanifold we define $H_k(Y,\R)$
and $H_k(Y;Z,\R)$ using {\it smooth\/} simplices, as in
\cite[\S V.5]{Bred}. Then the pairing between (singular)
homology and (de Rham) cohomology is defined at the chain
level by integrating $k$-forms over $k$-simplices.

Suppose $X$ is a compact SL $m$-fold in $M$ with conical
singularities $x_1,\ldots,x_n$ and cones $C_1,\ldots,C_n$, and set
$X'=X\sm\{x_1,\ldots,x_n\}$ and $\Si_i=C_i\cap{\cal S}^{2m-1}$ as
in \S\ref{co32}. Then by \cite[\S 2.4]{Joyc3} there is a natural
long exact sequence
\e
\cdots\ra
H^k_{\rm cs}(X',\R)\ra H^k(X',\R)\ra\bigoplus_{i=1}^n
H^k(\Si_i,\R)\ra H^{k+1}_{\rm cs}(X',\R)\ra\cdots,
\label{co3eq3}
\e
and natural isomorphisms
\begin{gather}
H_k\bigl(X;\{x_1,\ldots,x_n\},\R\bigr)^*\!\cong\!
H^k_{\rm cs}(X',\R)\!\cong\!H_{m-k}(X',\R)\!\cong\!H^{m-k}(X',\R)^*
\label{co3eq4}\\
\text{and}\quad
H^k_{\rm cs}(X',\R)\cong H_k(X,\R)^*
\quad\text{for all $k>1$.}
\label{co3eq5}
\end{gather}
The inclusion $\iota:X\ra M$ induces homomorphisms~$\iota_*:
H_k(X,\R)\ra H_k(M,\R)$.

If $(Y,g)$ is a compact Riemannian manifold, then {\it Hodge theory}
shows that each class in $H^k(Y,\R)$ is represented by a unique
$k$-form $\al$ with $\d\al=\d^*\al=0$. Here is an analogue of this
on $X'$ when $k=1$, given in~\cite[Th.~5.4]{Joyc3}.

\begin{thm} Let\/ $(M,J,\om,\Om)$ be an almost Calabi--Yau $m$-fold,
and define $\psi:M\ra(0,\iy)$ as in \eq{co2eq3}. Suppose $X$ is a
compact SL\/ $m$-fold in $M$ with conical singularities at
$x_1,\ldots,x_n$, and let\/ $X',K,R',\Si_i,\Up_i,\phi_i,S_i$
and\/ $\mu_i$ be as in Definition \ref{co3def3}, $\D_\sSii$
as in Definition \ref{co3def2}, and\/ $\rho$ as in Definition
\ref{co3def4}. Define
\e
\begin{split}
Y_\sXp=\bigl\{\al\in C^\iy(T^*X'):\,&\d\al=0,\quad
\d^*(\psi^m\al)=0,\\
&\text{$\md{\na^k\al}=O(\rho^{-1-k})$ for $k\ge 0$}\bigr\}.
\end{split}
\label{co3eq6}
\e
Then $\pi:Y_\sXp\!\ra\!H^1(X',\R)$ given by $\pi:\al\!\mapsto\![\al]$
is an isomorphism. Furthermore:
\begin{itemize}
\item[{\rm(a)}] Fix $\al\in Y_\sXp$. By Hodge theory there exists
a unique $\ga_i\in C^\iy(T^*\Si_i)$ with\/ $\d\ga_i=\d^*\ga_i=0$
for $i=1,\ldots,n$, such that the image of\/ $\pi(\al)$ under the
map $H^1(X',\R)\ra\bigoplus_{i=1}^nH^1(\Si_i,\R)$ of\/ \eq{co3eq3}
is $\bigl([\ga_1],\ldots,[\ga_n]\bigr)$. There exist unique $T_i\in
C^\iy\bigl(\Si_i\t(0,R')\bigr)$ for $i=1,\ldots,n$ such that
\end{itemize}
\ea
(\Up_i\circ\phi_i)^*(\al)&=\pi_i^*(\ga_i)+\d T_i
\quad \text{on $\Si_i\t(0,R')$ for $i=1,\ldots,n,$ and}
\label{co3eq7}\\
\na^kT_i(\si,r)&=O(r^{\nu_i-k})\qquad
\begin{aligned}
&\text{as\/ $r\ra 0$, for all\/ $k\ge 0$ and}\\
&\text{$\nu_i\in(0,\mu_i-2)$ with\/ $(0,\nu_i]\cap\D_\sSii=\emptyset$.}
\end{aligned}
\label{co3eq8}
\ea
\begin{itemize}
\item[{\rm(b)}] Suppose $\ga_i\in C^\iy(T^*\Si_i)$ with\/
$\d\ga_i=\d^*\ga_i=0$ for $i=1,\ldots,n$, and the image
of\/ $\bigl([\ga_1],\ldots,[\ga_n]\bigr)$ under
$\bigoplus_{i=1}^nH^1(\Si_i,\R)\ra H^2_{\rm cs}(X',\R)$ in
\eq{co3eq3} is $[\be]$ for some exact\/ $2$-form $\be$ on $X'$
supported on $K$. Then there exists $\al\in C^\iy(T^*X')$ with\/
$\d\al=\be$, $\d^*(\psi^m\al)=0$ and\/ $\md{\na^k\al}=O(\rho^{-1-k})$
for $k\ge 0$, such that\/ \eq{co3eq7} and\/ \eq{co3eq8} hold
for~$T_i\in C^\iy\bigl(\Si_i\t(0,R')\bigr)$.
\item[{\rm(c)}] Let\/ $f\in C^\iy(X')$ with\/ $\md{\na^kf}=
O(\rho^{-1-k})$ for $k\ge 0$ and\/ $\int_{X'}f\,\d V=0$.
Then there exists a unique exact\/ $1$-form $\al$ on $X'$ with\/
$\d^*(\psi^m\al)=f$ and\/ $\md{\na^k\al}=O(\rho^{-1-k})$ for
$k\ge 0$, such that\/ \eq{co3eq7} and\/ \eq{co3eq8} hold for
$\ga_i=0$ and\/~$T_i\in C^\iy\bigl(\Si_i\t(0,R')\bigr)$.
\end{itemize}
\label{co3thm1}
\end{thm}

\subsection{Lagrangian Neighbourhood Theorems and regularity}
\label{co34}

We recall some symplectic geometry, which can be found in McDuff
and Salamon \cite{McSa}. Let $N$ be a real $m$-manifold. Then its
tangent bundle $T^*N$ has a {\it canonical symplectic form}
$\hat\om$, defined as follows. Let $(x_1,\ldots,x_m)$ be local
coordinates on $N$. Extend them to local coordinates
$(x_1,\ldots,x_m,y_1,\ldots,y_m)$ on $T^*N$ such that
$(x_1,\ldots,y_m)$ represents the 1-form $y_1\d x_1+\cdots+y_m
\d x_m$ in $T_{(x_1,\ldots,x_m)}^*N$. Then~$\hat\om=\d x_1\w\d y_1+
\cdots+\d x_m\w\d y_m$.

Identify $N$ with the zero section in $T^*N$. Then $N$ is a
{\it Lagrangian submanifold\/} of $T^*N$. The {\it Lagrangian
Neighbourhood Theorem} \cite[Th.~3.33]{McSa} shows that any
compact Lagrangian submanifold $N$ in a symplectic manifold
looks locally like the zero section in~$T^*N$.

\begin{thm} Let\/ $(M,\om)$ be a symplectic manifold and\/
$N\subset M$ a compact Lagrangian submanifold. Then there
exists an open tubular neighbourhood\/ $U$ of the zero
section $N$ in $T^*N$, and an embedding $\Phi:U\ra M$ with\/
$\Phi\vert_N=\id:N\ra N$ and\/ $\Phi^*(\om)=\hat\om$, where
$\hat\om$ is the canonical symplectic structure on~$T^*N$.
\label{co3thm2}
\end{thm}

In \cite[\S 4]{Joyc3} we extend Theorem \ref{co3thm2}
to situations involving conical singularities, first to
{\it SL cones},~\cite[Th.~4.3]{Joyc3}.

\begin{thm} Let\/ $C$ be an SL cone in $\C^m$ with isolated
singularity at\/ $0$, and set\/ $\Si=C\cap{\cal S}^{2m-1}$.
Define $\iota:\Si\t(0,\iy)\ra\C^m$ by $\iota(\si,r)=r\si$,
with image $C\sm\{0\}$. For $\si\in\Si$, $\tau\in T_\si^*\Si$,
$r\in(0,\iy)$ and\/ $u\in\R$, let\/ $(\si,r,\tau,u)$ represent
the point\/ $\tau+u\,\d r$ in $T^*_{\smash{(\si,r)}}\bigl(\Si\!
\t\!(0,\iy)\bigr)$. Identify $\Si\t(0,\iy)$ with the zero section
$\tau\!=\!u\!=\!0$ in $T^*\bigl(\Si\t(0,\iy)\bigr)$. Define an
action of\/ $(0,\iy)$ on $T^*\bigl(\Si\!\t\!(0,\iy)\bigr)$ by
\e
t:(\si,r,\tau,u)\longmapsto (\si,tr,t^2\tau,tu)
\quad\text{for $t\in(0,\iy),$}
\label{co3eq9}
\e
so that\/ $t^*(\hat\om)\!=\!t^2\hat\om$, for $\hat\om$ the
canonical symplectic structure on~$T^*\bigl(\Si\!\t\!(0,\iy)\bigr)$.

Then there exists an open neighbourhood\/ $U_\sC$ of\/
$\Si\t(0,\iy)$ in $T^*\bigl(\Si\t(0,\iy)\bigr)$ invariant under
\eq{co3eq9} given by
\e
U_\sC=\bigl\{(\si,r,\tau,u)\in T^*\bigl(\Si\t(0,\iy)\bigr):
\bmd{(\tau,u)}<2\ze r\bigr\}\quad\text{for some $\ze>0,$}
\label{co3eq10}
\e
where $\md{\,.\,}$ is calculated using the cone metric $\iota^*(g')$
on $\Si\t(0,\iy)$, and an embedding $\Phi_\sC:U_\sC\ra\C^m$
with\/ $\Phi_\sC\vert_{\Si\t(0,\iy)}=\iota$, $\Phi_{\sst
C}^*(\om')=\hat\om$ and\/ $\Phi_\sC\circ t=t\,\Phi_\sC$
for all\/ $t>0$, where $t$ acts on $U_\sC$ as in \eq{co3eq9}
and on $\C^m$ by multiplication.
\label{co3thm3}
\end{thm}

In \cite[Th.~4.4]{Joyc3} we construct a particular choice of
$\phi_i$ in Definition~\ref{co3def3}.

\begin{thm} Let\/ $(M,J,\om,\Om)$, $\psi,X,n,x_i,\up_i,C_i,\Si_i,
\mu_i,R,\Up_i$ and\/ $\iota_i$ be as in Definition \ref{co3def3}.
Theorem \ref{co3thm3} gives $\ze>0$, neighbourhoods $U_\sCi$
of\/ $\Si_i\t(0,\iy)$ in $T^*\bigl(\Si_i\t(0,\iy)\bigr)$ and
embeddings $\Phi_\sCi:U_\sCi\ra\C^m$ for~$i=1,\ldots,n$.

Then for sufficiently small\/ $R'\in(0,R]$ there exist unique
closed\/ $1$-forms $\eta_i$ on $\Si_i\t(0,R')$ for $i=1,\ldots,n$
written $\eta_i(\si,r)=\eta_i^1(\si,r)+\eta_i^2(\si,r)\d r$ for
$\eta_i^1(\si,r)\in T_\si^*\Si_i$ and\/ $\eta_i^2(\si,r)\in\R$,
and satisfying $\md{\eta_i(\si,r)}<\ze r$ and\/
$\bmd{\na^k\eta_i}=O(r^{\mu_i-1-k})$ as $r\ra 0$ for $k=0,1,$ 
computing $\na,\md{\,.\,}$ using the cone metric $\iota_i^*(g')$,
such that the following holds.

Define $\phi_i:\Si_i\t(0,R')\ra B_R$ by $\phi_i(\si,r)=\Phi_{\sst
C_i}\bigl(\si,r,\eta_i^1(\si,r),\eta_i^2(\si,r)\bigr)$. Then
$\Up_i\circ\phi_i:\Si_i\t(0,R')\ra M$ is a diffeomorphism
$\Si_i\t(0,R')\ra S_i$ for open sets $S_1,\ldots,S_n$ in $X'$
with\/ $\bar S_1,\ldots,\bar S_n$ disjoint, and\/ $K=X'\sm(S_1
\cup\cdots\cup S_n)$ is compact. Also $\phi_i$ satisfies
\eq{co3eq2}, so that\/ $R',\phi_i,S_i,K$ satisfy
Definition~\ref{co3def3}.
\label{co3thm4}
\end{thm}

In \cite[\S 5]{Joyc3} we study the asymptotic behaviour of the maps
$\phi_i$ of Theorem \ref{co3thm4}, using the elliptic regularity of
the special Lagrangian condition. Combining \cite[Th.~5.1]{Joyc3},
\cite[Lem.~4.5]{Joyc3} and \cite[Th.~5.5]{Joyc3} proves:

\begin{thm} In the situation of Theorem \ref{co3thm4} we have
$\eta_i=\d A_i$ for $i=1,\ldots,n$, where $A_i:\Si_i\t(0,R')\ra\R$
is given by $A_i(\si,r)=\int_0^r\eta_i^2(\si,s)\d s$. Suppose
$\mu_i'\in(2,3)$ with\/ $(2,\mu_i']\cap\D_\sSii=\emptyset$ for
$i=1,\ldots,n$. Then
\e
\begin{gathered}
\bmd{\na^k(\phi_i-\iota_i)}=O(r^{\mu_i'-1-k}),\quad
\bmd{\na^k\eta_i}=O(r^{\mu_i'-1-k})\quad\text{and}\\
\bmd{\na^kA_i}=O(r^{\mu_i'-k})
\quad\text{as $r\ra 0$ for all\/ $k\ge 0$ and\/ $i=1,\ldots,n$.}
\end{gathered}
\label{co3eq11}
\e

Hence $X$ has conical singularities at $x_i$ with cone $C_i$
and rate $\mu_i'$, for all possible rates $\mu_i'$ allowed by
Definition \ref{co3def3}. Therefore, the definition of
conical singularities is essentially independent of the
choice of rate~$\mu_i$.
\label{co3thm5}
\end{thm}

Next we extend Theorem \ref{co3thm2} to SL $m$-folds with
conical singularities \cite[Th.~4.6]{Joyc3}, in a way
compatible with Theorems \ref{co3thm3} and~\ref{co3thm4}.

\begin{thm} Suppose $(M,J,\om,\Om)$ is an almost Calabi--Yau
$m$-fold and\/ $X$ a compact SL\/ $m$-fold in $M$ with conical
singularities at\/ $x_1,\ldots,x_n$. Let the notation $\psi,\up_i,
C_i,\Si_i,\mu_i,R,\Up_i$ and\/ $\iota_i$ be as in Definition
\ref{co3def3}, and let\/ $\ze,U_\sCi,\allowbreak
\Phi_\sCi,\allowbreak R',\allowbreak \eta_i,\allowbreak
\eta_i^1,\eta_i^2,\phi_i,S_i$ and\/ $K$ be as in Theorem~\ref{co3thm4}.

Then making $R'$ smaller if necessary, there exists an open tubular
neighbourhood\/ $U_\sXp\subset T^*X'$ of the zero section
$X'$ in $T^*X'$, such that under $\d(\Up_i\circ\phi_i):T^*\bigl(
\Si_i\t(0,R')\bigr)\ra T^*X'$ for $i=1,\ldots,n$ we have
\e
\bigl(\d(\Up_i\circ\phi_i)\bigr)^*(U_\sXp)=\bigl\{(\si,r,\tau,u)
\in T^*\bigl(\Si_i\t(0,R')\bigr):\bmd{(\tau,u)}<\ze r\bigr\},
\label{co3eq12}
\e
and there exists an embedding $\Phi_\sXp:U_\sXp\ra M$ with\/
$\Phi_\sXp\vert_{X'}=\id:X'\ra X'$ and\/ $\Phi_\sXp^*(\om)=\hat\om$,
where $\hat\om$ is the canonical symplectic structure on $T^*X'$,
such that
\e
\Phi_\sXp\circ\d(\Up_i\circ\phi_i)(\si,r,\tau,u)\equiv\Up_i\circ
\Phi_\sCi\bigl(\si,r,\tau+\eta_i^1(\si,r),u+\eta_i^2(\si,r)\bigr)
\label{co3eq13}
\e
for all\/ $i=1,\ldots,n$ and\/ $(\si,r,\tau,u)\in T^*\bigl(\Si_i\t(0,R')
\bigr)$ with\/ $\bmd{(\tau,u)}<\ze r$. Here $\md{(\tau,u)}$ is computed
using the cone metric $\iota_i^*(g')$ on~$\Si_i\t(0,R')$.
\label{co3thm6}
\end{thm}

In \cite[Th.~4.8]{Joyc3} we extend Theorem \ref{co3thm6} to {\it
families} of almost Calabi--Yau $m$-folds $(M,J^s,\om^s,\Om^s)$
for $s\in\F$. If $\om^s\vert_{X'}$ is not exact then we cannot
deform $X'$ to a Lagrangian $m$-fold in $(M,\om^s)$. Therefore
we replace the condition $\Phi_\sXp^*(\om)=\hat\om$ in Theorem
\ref{co3thm6} by $(\Phi^s_\sXp)^*(\om^s)=\hat\om+\pi^*(\nu^s)$,
where $\nu^s$ is a compactly-supported closed 2-form on~$X'$.

\begin{thm} Let\/ $(M,J,\om,\Om)$ be an almost Calabi--Yau
$m$-fold and\/ $X$ a compact SL\/ $m$-fold in $M$ with conical
singularities at $x_1,\ldots,x_n$. Let the notation
$R,\Up_i,\ze,\Phi_\sCi,R',\eta_i,\eta_i^1,\eta_i^2,
\phi_i,S_i,K$ be as in Theorem \ref{co3thm4}, and let\/
$U_\sXp,\Phi_\sXp$ be as in Theorem \ref{co3thm6}.
Suppose $\bigl\{(M,J^s,\om^s,\Om^s):s\in \F\bigr\}$ is
a smooth family of deformations of $(M,J,\om,\Om)$ with base
space $\F\subset\R^d$. Define $\psi^s:M\ra(0,\iy)$ for
$s\in\F$ as in \eq{co2eq3}, but using~$\om^s,\Om^s$.

Then making $R,R'$ and\/ $U_\sXp$ smaller if necessary,
for some connected open $\F'\subseteq\F$ with\/
$0\in\F'$ and all\/ $s\in \F'$ there exist
\begin{itemize}
\item[{\rm(a)}] isomorphisms $\up_i^s:\C^m\ra T_{x_i}M$ for
$i=1,\ldots,n$ with\/ $\up_i^0=\up_i$, $(\up_i^s)^*(\om^s)=\om'$
and\/~$(\up_i^s)^*(\Om)=\psi^s(x_i)^m\Om'$,
\item[{\rm(b)}] embeddings $\Up_i^s:B_R\ra M$ for $i=1,\ldots,n$
with\/ $\Up_i^0=\Up_i$, $\Up_i^s(0)=x_i$, $\d\Up_i^s\vert_0=\up_i^s$
and\/~$(\Up_i^s)^*(\om^s)=\om'$,
\item[{\rm(c)}] a closed\/ $2$-form $\nu^s\in C^\iy(\La^2T^*X')$
supported in $K\subset X'$ with\/ $\nu^0=0$, and
\item[{\rm(d)}] an embedding $\Phi^s_\sXp\!:\!U_\sXp\!\ra\!M$
with\/ $\Phi_\sXp^0\!=\!\Phi_\sXp$ and\/~$(\Phi^s_\sXp
)^*(\om^s)\!=\!\hat\om\!+\!\pi^*(\nu^s)$,
\end{itemize}
all depending smoothly on $s\in\F'$ with
\e
\Phi_\sXp^s\circ\d(\Up_i\circ\phi_i)(\si,r,\tau,u)\equiv\Up_i^s\circ
\Phi_\sCi\bigl(\si,r,\tau+\eta_i^1(\si,r),u+\eta_i^2(\si,r)\bigr)
\label{co3eq14}
\e
for all\/ $s\in \F'$, $i=1,\ldots,n$ and\/ $(\si,r,\tau,u)\in
T^*\bigl(\Si_i\t(0,R')\bigr)$ with\/~$\bmd{(\tau,u)}<\ze r$.
\label{co3thm7}
\end{thm}

The 2-forms $\nu^s$ in Theorem \ref{co3thm7} define classes
$[\nu^s]$ in $H^2_{\rm cs}(X',\R)$. In \cite[Th.~4.9]{Joyc3}
we investigate these classes, and the freedom to choose~$\nu^s$.

\begin{thm} In the situation of Theorem \ref{co3thm7}, under
the isomorphism \eq{co3eq5}, the class $[\nu^s]\in H^2_{\rm cs}
(X',\R)$ is identified with the map $H_2(X,\R)\ra\R$ given by
$\ga\mapsto\iota_*(\ga)\cdot[\om^s]$, where $\iota:X\ra M$ is
the inclusion, $\iota_*:H_2(X,\R)\ra H_2(M,\R)$ the induced map,
and\/ $[\om^s]\in H^2(M,\R)$. Thus $[\nu^s]$ depends only on
$X,M$ and\/~$[\om^s]\in H^2(M,\R)$.

Let\/ $V\cong H^2_{\rm cs}(X',\R)$ be a vector space of smooth
closed\/ $2$-forms on $X'$ supported in $K$ representing
$H^2_{\rm cs}(X',\R)$. Then making $\F'$ smaller if
necessary, we can choose $\Up_i^s,\nu^s$ and\/ $\Phi_\sXp^s$ in
Theorem \ref{co3thm7} so that\/ $\nu^s\in V$ for all\/ $s\in\F'$.
In particular, if\/ $[\nu^s]=0$ in $H^2_{\rm cs}(X',\R)$ then we can
choose~$\nu^s=0$.
\label{co3thm8}
\end{thm}

\section{Asymptotically Conical SL $m$-folds}
\label{co4}

Let $C$ be an SL cone in $\C^m$ with an isolated singularity
at 0. Section \ref{co3} considered SL $m$-folds with conical
singularities, which are asymptotic to $C$ at 0. We now
discuss {\it Asymptotically Conical\/} SL $m$-folds $L$ in
$\C^m$, which are asymptotic to $C$ at infinity. Here is
the definition.

\begin{dfn} Let $C$ be an SL cone in $\C^m$ with isolated
singularity at 0 for $m>2$, as in Definition \ref{co3def1}, and
let $\Si=C\cap{\cal S}^{2m-1}$, so that $\Si$ is a compact,
nonsingular $(m-1)$-manifold, not necessarily connected.
Let $g_\sSi$ be the metric on $\Si$ induced by the
metric $g'$ on $\C^m$ in \eq{co2eq1}, and $r$ the radius
function on $\C^m$. Define $\iota:\Si\t(0,\iy)\ra\C^m$ by
$\iota(\si,r)=r\si$. Then the image of $\iota$ is $C\sm\{0\}$,
and $\iota^*(g')=r^2g_\sSi+\d r^2$ is the cone metric
on~$C\sm\{0\}$.

Let $L$ be a closed, nonsingular SL $m$-fold in $\C^m$ and
$\la<2$. We call $L$ {\it Asymptotically Conical (AC)}
with {\it rate} $\la$ and {\it cone} $C$ if there exists a
compact subset $K\subset L$ and a diffeomorphism $\vp:\Si\t
(T,\iy)\ra L\sm K$ for some $T>0$, such that
\e
\bmd{\na^k(\vp-\iota)}=O(r^{\la-1-k})
\quad\text{as $r\ra\iy$ for $k=0,1$.}
\label{co4eq1}
\e
Here $\na,\md{\,.\,}$ are computed using the cone metric
$\iota^*(g')$ on~$\Si\t(T,\iy)$.
\label{co4def1}
\end{dfn}

This is very similar to Definition \ref{co3def3}, and in fact
there are strong parallels between the theories of SL $m$-folds
with conical singularities and of Asymptotically Conical SL
$m$-folds. We recall some results from \cite[\S 7]{Joyc3},
including versions of the material in \S\ref{co34}. We
continue to assume $m>2$ throughout.

\subsection{Cohomological invariants of AC SL $m$-folds}
\label{co41}

Let $L$ be an AC SL $m$-fold in $\C^m$ with cone $C$, and let
$\Si=C\cap{\cal S}^{2m-1}$. Using the notation of \S\ref{co33},
as in \eq{co3eq3} there is a long exact sequence
\e
\cdots\ra
H^k_{\rm cs}(L,\R)\ra H^k(L,\R)\ra H^k(\Si,\R)\ra
H^{k+1}_{\rm cs}(L,\R)\ra\cdots.
\label{co4eq2}
\e
Following \cite[Def.~7.2]{Joyc3} we define {\it cohomological
invariants\/} $Y(L),Z(L)$ of~$L$.

\begin{dfn} Let $L$ be an AC SL $m$-fold in $\C^m$ with cone $C$,
and let $\Si=C\cap{\cal S}^{2m-1}$. As $\om',\Im\Om'$ in \eq{co2eq1}
are closed forms with $\om'\vert_L\equiv\Im\Om'\vert_L=0$, they
define classes in the relative de Rham cohomology groups $H^k(\C^m;
L,\R)$ for $k=2,m$. For $k>1$ we have the exact sequence
\begin{equation*}
0=H^{k-1}(\C^m,\R)\ra H^{k-1}(L,\R){\buildrel\cong\over\longra}
H^k(\C^m;L,\R)\ra H^k(\C^m,\R)=0.
\end{equation*}
Let $Y(L)\in H^1(\Si,\R)$ be the image of $[\om']$ in
$H^2(\C^m;L,\R)\cong H^1(L,\R)$ under $H^1(L,\R)\ra H^1(\Si,R)$
in \eq{co4eq2}, and $Z(L)\in H^{m-1}(\Si,\R)$ be the image of
$[\Im\Om']$ in $H^m(\C^m;L,\R)\cong H^{m-1}(L,\R)$ under
$H^{m-1}(L,\R)\ra H^{m-1}(\Si,R)$ in~\eq{co4eq2}.
\label{co4def2}
\end{dfn}

Here are some conditions for $Y(L)$ or $Z(L)$ to be
zero,~\cite[Prop.~7.3]{Joyc3}.

\begin{prop} Let\/ $L$ be an AC SL\/ $m$-fold in $\C^m$ with
cone $C$ and rate $\la$, and let\/ $\Si=C\cap{\cal S}^{2m-1}$.
If\/ $\la<0$ or $b^1(L)=0$ then $Y(L)=0$. If\/ $\la<2-m$ or
$b^0(\Si)=1$ then~$Z(L)=0$.
\label{co4prop}
\end{prop}

Here is a (trivial) lemma on {\it dilations} of AC SL $m$-folds.

\begin{lem} Let\/ $L$ be an AC SL\/ $m$-fold in $\C^m$ with rate
$\la$ and cone $C$, and let\/ $t>0$. Then $tL=\{t\,{\bf x}:{\bf x}
\in L\}$ is also an AC SL\/ $m$-fold in $\C^m$ with rate $\la$
and cone $C$, satisfying $Y(tL)=t^2Y(L)$ and\/~$Z(tL)=t^mZ(L)$.
\label{co4lem}
\end{lem}

\subsection{Lagrangian Neighbourhood Theorems and regularity}
\label{co42}

Next we give versions of \S\ref{co34} for AC SL $m$-folds
rather than SL $m$-folds with conical singularities. Here
\cite[Th.~7.4]{Joyc3} is the analogue of Theorem~\ref{co3thm4}.

\begin{thm} Let\/ $C$ be an SL cone in $\C^m$ with isolated
singularity at\/ $0$, and set\/ $\Si=C\cap{\cal S}^{2m-1}$.
Define $\iota:\Si\t(0,\iy)\ra\C^m$ by $\iota(\si,r)=r\si$. Let\/
$\ze$, $U_\sC\subset T^*\bigl(\Si\t(0,\iy)\bigr)$ and\/
$\Phi_\sC:U_\sC\ra\C^m$ be as in Theorem~\ref{co3thm3}.

Suppose $L$ is an AC SL\/ $m$-fold in $\C^m$ with cone $C$
and rate $\la<2$. Then there exists a compact\/ $K\subset L$
and a diffeomorphism $\vp:\Si\t(T,\iy)\ra L\sm K$ for some
$T>0$ satisfying \eq{co4eq1}, and a closed\/ $1$-form $\chi$ on
$\Si\t(T,\iy)$ written $\chi(\si,r)=\chi^1(\si,r)+\chi^2(\si,r)\d r$
for $\chi^1(\si,r)\in T_\si^*\Si$ and\/ $\chi^2(\si,r)\in\R$,
satisfying
\e
\begin{gathered}
\bmd{\chi(\si,r)}<\ze r,\quad \vp(\si,r)\equiv
\Phi_\sC\bigl(\si,r,\chi^1(\si,r),\chi^2(\si,r)\bigr)\\
\text{and}\quad\bmd{\na^k\chi}=O(r^{\la-1-k})
\quad\text{as $r\ra\iy$ for $k=0,1,$}
\end{gathered}
\label{co4eq3}
\e
computing $\na,\md{\,.\,}$ using the cone metric~$\iota^*(g')$.
\label{co4thm1}
\end{thm}

The next two theorems are analogous to Theorem \ref{co3thm5}.
In \cite[Prop.~7.6]{Joyc3} and \cite[Th.~7.7]{Joyc3} we
decompose $\chi$ in Theorem~\ref{co4thm1}.

\begin{thm} In Theorem \ref{co4thm1} we have $[\chi]=Y(L)$ in
$H^1\bigl(\Si\t(T,\iy),\R\bigr)\cong H^1(\Si,\R)$, where $Y(L)$
is as in Definition \ref{co4def2}. Let\/ $\ga$ be the unique
$1$-form on $\Si$ with\/ $\d\ga=\d^*\ga=0$ and\/ $[\ga]=Y(L)\in
H^1(\Si,\R)$, which exists by Hodge theory. Then $\chi=\pi^*(\ga)
+\d E$, where $\pi:\Si\t(T,\iy)\ra\Si$ is the projection and\/
$E\in C^\iy\bigl(\Si\t(T,\iy)\bigr)$ satisfies
\e
\begin{gathered}
\bmd{\na^k(\vp\!-\!\iota)}=O(r^{\la-1-k}),\;\>
\bmd{\na^k\chi}=O(r^{\la-1-k}),\;\>
\bmd{\na^{k+1}E}=O(r^{\la-1-k})\\
\text{for all $k\ge 0$, and}\quad
\md{E}=\begin{cases}O(r^\la), & \la\ne 0, \\
O\bigl(\md{\log r}\bigr), & \la=0. \end{cases}
\end{gathered}
\label{co4eq4}
\e
Here $\na,\md{\,.\,}$ are computed using the cone
metric $\iota^*(g')$ on~$\Si\t(T,\iy)$.
\label{co4thm2}
\end{thm}

Then \cite[Th.~7.11]{Joyc3} we improve the rate of convergence~$\la$.

\begin{thm} Let\/ $L$ be an AC SL\/ $m$-fold in $\C^m$ with
cone $C$ and rate $\la$. Set\/ $\Si=C\cap{\cal S}^{2m-1}$,
and let\/ $\D_\sSi,N_\sSi$ be as in Definition \ref{co3def2}.
Let\/ $\iota,T,\vp,\chi$ be as in Theorem \ref{co4thm1},
and\/ $Y(L),\ga,E$ as in Theorem \ref{co4thm2}. Then
\begin{itemize}
\item[{\rm(a)}] Suppose $\la,\la'$ lie in the same
connected component of\/ $\R\sm\D_\sSi$. Then
\e
\begin{gathered}
\bmd{\na^k(\vp-\iota)}=O(r^{\la'-1-k}),\quad
\bmd{\na^k\chi}=O(r^{\la'-1-k})\quad\text{and}\\
\bmd{\na^kE}=O(r^{\la'-k})
\quad\text{as $r\ra\iy$ for all\/ $k\ge 0$.}
\end{gathered}
\label{co4eq5}
\e
Hence $L$ is an AC SL\/ $m$-fold with rate $\la'$. In particular,
if\/ $\la\in(2-m,0)$ then $L$ is an AC SL\/ $m$-fold with rate
$\la'$ for all\/~$\la'\in(2-m,0)$.
\item[{\rm(b)}] Suppose $0\le\la<\min\bigl(\D_\sSi\cap(0,\iy)\bigr)$.
Then adding a constant to $E$ if necessary, for all\/ $\la'\in\bigl(
\max(-2,2-m),0\bigr)$ we have
\e
\bmd{\na^kE}=O(r^{\la'-k})\quad\text{as $r\ra\iy$ for all\/ $k\ge 0$.}
\label{co4eq6}
\e
Thus if\/ $Y(L)=0=\ga$ then $L$ is an AC SL\/ $m$-fold with rate
$\la'$, and if\/ $Y(L)\ne 0\ne\ga$ then $L$ is an AC SL\/ $m$-fold
with rate~$0$.
\end{itemize}
\label{co4thm3}
\end{thm}

Here is the analogue of Theorem \ref{co3thm6}, proved
in~\cite[Th.~7.5]{Joyc3}.

\begin{thm} Suppose $L$ is an AC SL\/ $m$-fold in $\C^m$ with
cone $C$. Let\/ $\Si,
\allowbreak
\iota,
\allowbreak
\ze,
\allowbreak
U_\sC,
\allowbreak
\Phi_\sC,
\allowbreak
K,T,\vp,\chi,\chi^1,\chi^2$ be as in Theorem \ref{co4thm1}.
Then making $T,K$ larger if necessary, there exists an open tubular
neighbourhood\/ $U_\sL\subset T^*L$ of the zero section $L$ in $T^*L$,
such that under $\d\vp:T^*\bigl(\Si\t(T,\iy)\bigr)\ra T^*L$ we have
\e
(\d\vp)^*(U_\sL)=\bigl\{(\si,r,\tau,u)\in
T^*\bigl(\Si\t(T,\iy)\bigr):\bmd{(\tau,u)}<\ze r\bigr\},
\label{co4eq7}
\e
and there exists an embedding $\Phi_\sL:U_\sL\ra\C^m$ with\/
$\Phi_\sL\vert_L=\id:L\ra L$ and\/ $\Phi_\sL^*(\om')=\hat\om$,
where $\hat\om$ is the canonical symplectic structure on $T^*L$,
such that
\e
\Phi_\sL\circ\d\vp(\si,r,\tau,u)\equiv
\Phi_\sC\bigl(\si,r,\tau+\chi^1(\si,r),u+\chi^2(\si,r)\bigr)
\label{co4eq8}
\e
for all\/ $(\si,r,\tau,u)\!\in\!T^*\bigl(\Si\!\t\!(T,\iy)\bigr)$
with\/ $\md{(\tau,u)}<\ze r$, computing $\md{\,.\,}$
using~$\iota^*(g')$.
\label{co4thm4}
\end{thm}

In \cite[Th.~7.10]{Joyc3} we study the {\it bounded harmonic
functions\/} on~$L$.

\begin{thm} Suppose $L$ is an AC SL\/ $m$-fold in $\C^m$, with cone
$C$. Let\/ $\Si,T$ and\/ $\vp$ be as in Theorem \ref{co4thm1}. Let\/
$l=b^0(\Si)$, and\/ $\Si^1,\ldots,\Si^l$ be the connected components
of\/ $\Si$. Let\/ $V$ be the vector space of bounded harmonic
functions on $L$. Then $\dim V=l$, and for each\/
${\bf c}=(c^1,\ldots,c^l)\in\R^l$ there exists a unique $v^{\bf c}
\in V$ such that for all\/ $j=1,\ldots,l$, $k\ge 0$ and\/ $\be\in
(2-m,0)$ we have
\e
\na^k\bigl(\vp^*(v^{\bf c})-c^j\,\bigr)=O\bigl(\md{{\bf c}}
r^{\be-k}\bigr) \quad\text{on $\Si^j\t(T,\iy)$ as $r\ra\iy$.}
\label{co4eq9}
\e
Note also that\/ $V=\{v^{\bf c}:{\bf c}\in\R^l\}$
and\/~$v^{(1,\ldots,1)}\equiv 1$.
\label{co4thm5}
\end{thm}

\section{Review of the main results of \cite{Joyc5}}
\label{co5}

Our goal is to generalize the results of \cite{Joyc5} to
more complicated situations. This gave me a problem in
writing this paper, as I want it to make sense on its
own without constant reference to \cite{Joyc5}, but to
control the length I don't want to reproduce large parts
of \cite{Joyc5} as introductory material here.

The solution I have adopted is to reproduce only the three
major theorems from \cite{Joyc5} in this section, with some
supporting definitions and explanations. However, much of
\S\ref{co6}--\S\ref{co8} (for instance, Definitions
\ref{co6def1}--\ref{co6def3} below) has in effect been
copied from \cite[\S 6--\S 7]{Joyc5} and then modified. I
hope this is more economical and readable than reproducing
long definitions from \cite{Joyc5} unchanged, and then
explaining later how to change them.

The other way I save space is that if the proof of a
result in \cite{Joyc5} requires only superficial changes
for the new situations in this paper, then I give only
the result but not the proof here, or else make only
brief comments on how to adapt the proof in~\cite{Joyc5}.

\subsection{An analytic existence result for SL $m$-folds}
\label{co51}

The results of \cite{Joyc5} hinged upon an existence result
\cite[Th.~5.3]{Joyc5} for compact SL $m$-folds proved using
analysis. Here is some notation we will need.

\begin{dfn} Let $(N,g)$ be a Riemannian manifold with
Levi-Civita connection $\na$. For each integer $k\ge 0$,
define $C^k(N)$ to be the Banach space of functions $f$
on $N$ with $k$ continuous derivatives, for which the
norm $\cnm{f}k=\sum_{j=0}^k\sup_N\md{\na^jf}$ is finite.
Let~$C^\iy(N)=\bigcap_{k\ge 0}C^k(N)$.

For $q\ge 1$, define the {\it Lebesgue space} $L^q(N)$ to be
the Banach space of locally integrable functions $f$ on $N$
for which the norm $\lnm{f}q=(\int_N\md{f}^q\,\d V_g)^{1/q}$
is finite. For $q\ge 1$ and $k\ge 0$ an integer, define the 
{\it Sobolev space} $L^q_k(N)$ to be the set of $f\in L^q(N)$ 
such that $f$ is $k$ times weakly differentiable and 
$\md{\na^jf}\in L^q(N)$ for $j\le k$, with norm~$\snm{f}qk=
(\sum_{j=0}^k\int_N\md{\na^jf}^q\,\d V_g)^{1/q}$.
\label{co5def1}
\end{dfn}

The following definition \cite[Def.~5.2]{Joyc5} sets up the
situation we shall consider.

\begin{dfn} Let $(M,J,\om,\Om)$ be an almost Calabi--Yau $m$-fold,
with metric $g$. Let $N$ be a compact, oriented, immersed,
Lagrangian $m$-submanifold in $M$, with immersion $\iota:N\ra M$,
so that $\iota^*(\om)\equiv 0$. Define $h=\iota^*(g)$, so that
$(N,h)$ is a Riemannian manifold. Let $\d V$ be the volume form
on $N$ induced by the metric $h$ and orientation.

Let $\psi:M\ra(0,\iy)$ be the smooth function given in \eq{co2eq3}.
Then $\Om\vert_N$ is a complex $m$-form on $N$, and using \eq{co2eq3}
and the Lagrangian condition we find that $\bmd{\Om\vert_N}=\psi^m$,
calculating $\md{\,.\,}$ using $h$ on $N$. Therefore we may write
\e
\Om\vert_N=\psi^m{\rm e}^{i\th}\,\d V \quad\text{on $N$,}
\label{co5eq1}
\e
for some phase function ${\rm e}^{i\th}$ on $N$. Suppose that
$\cos\th\ge\ha$ on $N$. Then we can choose $\th$ to be a smooth
function $\th:N\ra(-\frac{\pi}{3},\frac{\pi}{3})$. Suppose that
$[\iota^*(\Im\Om)]=0$ in $H^m(N,\R)$. Then
$\int_N\psi^m\sin\th\,\d V=0$, by~\eq{co5eq1}.

Suppose we are given a finite-dimensional vector subspace
$W\subset C^\iy(N)$ with $1\in W$. Define $\pi_\sW:L^2(N)\ra W$
to be the projection onto $W$ using the $L^2$-inner product.

For $r>0$, define $\B_r\subset T^*N$ to be the bundle of 1-forms
$\al$ on $N$ with $\md{\al}<r$. Regard $\B_r$ as a noncompact
$2m$-manifold with natural projection $\pi:\B_r\ra N$, whose fibre
at $x\in N$ is the ball of radius $r$ about 0 in $T_x^*N$. We will
sometimes identify $N$ with the zero section of $\B_r$, and
write~$N\subset\B_r$.

At each $y\in\B_r$ with $\pi(y)=x\in N$, the Levi-Civita connection
$\na$ of $h$ on $T^*N$ defines a splitting $T_y\B_r=H\op V$ into
horizontal and vertical subspaces $H,V$, with $H\cong T_xN$ and
$V\cong T_x^*N$. Write $\hat\om$ for the natural symplectic
structure on $\B_r\subset T^*N$, defined using $T\B_r\cong H\op V$
and $H\cong V^*$. Define a natural Riemannian metric $\hat h$ on
$\B_r$ such that the subbundles $H,V$ are orthogonal, and
$\hat h\vert_H=\pi^*(h)$, $\hat h\vert_V=\pi^*(h^{-1})$.

Let $\hat\na$ be the connection on $T\B_r\cong H\op V$ given by
the lift of the Levi-Civita connection $\na$ of $h$ on $N$ in
the horizontal directions $H$, and by partial differentiation
in the vertical directions $V$, which is well-defined as
$T\B_r$ is naturally trivial along each fibre. Then $\hat\na$
preserves $\hat h,\hat\om$ and the splitting $T\B_r\cong H\op V$.
It is {\it not\/} torsion-free in general, but has torsion
$T(\hat\na)$ depending linearly on the Riemann curvature~$R(h)$.

As $N$ is a Lagrangian submanifold of $M$, by Theorem \ref{co3thm2}
the symplectic manifold $(M,\om)$ is locally isomorphic near $N$ to
$T^*N$ with its canonical symplectic structure. That is, for some
small $r>0$ there exists an immersion $\Phi:\B_r\ra M$ such that
$\Phi^*(\om)=\hat\om$ and $\Phi\vert_N=\iota$. Define an $m$-form
$\be$ on $\B_r$ by~$\be=\Phi^*(\Im\Om)$.

If $\al\in C^\iy(T^*N)$ with $\md{\al}<r$, write $\Ga(\al)$ for the
{\it graph} of $\al$ in $\B_r$. Then $\Phi_*(\Ga(\al))$ is a compact,
immersed submanifold in $M$ diffeomorphic to~$N$.
\label{co5def2}
\end{dfn}

With this notation, here is the existence result~\cite[Th.~5.3]{Joyc5}.

\begin{thm} Let\/ $\ka>1$ and\/ $A_1,\ldots,A_8>0$ be real,
and\/ $m\ge 3$ an integer. Then there exist\/ $\ep,K>0$ depending
only on $\ka,A_1,\ldots,A_8$ and\/ $m$ such that the following
holds.

Suppose $0<t\le\ep$ and Definition \ref{co5def2} holds with\/
$r=A_1t$, and
\begin{itemize}
\item[{\rm(i)}] $\lnm{\psi^m\sin\th}{2m/(m+2)}\le A_2t^{\ka+m/2}$,
$\cnm{\psi^m\sin\th}{0}\le A_2t^{\ka-1}$,
\newline
$\lnm{\d(\psi^m\sin\th)}{2m}\le A_2t^{\ka-3/2}$
and\/~$\lnm{\pi_\sW(\psi^m\sin\th)}{1}\le A_2t^{\ka+m-1}$.
\item[{\rm(ii)}] $\psi\ge A_3$ on $N$.
\item[{\rm(iii)}] $\cnm{\hat\na^k\be}{0}\le A_4t^{-k}$ for
$k=0,1,2$ and\/~$3$.
\item[{\rm(iv)}] The injectivity radius $\de(h)$
satisfies~$\de(h)\ge A_5t$.
\item[{\rm(v)}] The Riemann curvature $R(h)$
satisfies~$\cnm{R(h)}{0}\le A_6t^{-2}$.
\item[{\rm(vi)}] If\/ $v\in L^2_1(N)$ with\/ $\pi_\sW(v)=0$,
then $v\in L^{2m/(m-2)}(N)$ by the Sobolev Embedding Theorem,
and\/~$\lnm{v}{2m/(m-2)}\le A_7\lnm{\d v}{2}$.
\item[{\rm(vii)}] For all\/ $w\in W$ we have
$\lnm{\d^*\d w}{2m/(m+2)}\le\ha A_7^{-1}\lnm{\d w}{2}$.
\newline
For all\/ $w\in W$ with\/ $\int_Nw\,\d V=0$
we have~$\cnm{w}{0}\le A_8t^{1-m/2}\lnm{\d w}{2}$.
\end{itemize}
Here norms are computed using the metric $h$ on $N$ in {\rm(i)},
{\rm(v)}, {\rm(vi)} and\/ {\rm(vii)}, and the metric $\hat h$ on
$\B_{A_1t}$ in {\rm(iii)}. Then there exists $f\in C^\iy(N)$ with\/
$\int_Nf\,\d V=0$, such that\/ $\cnm{\d f}{0}\le Kt^\ka<A_1t$
and\/ $\ti N=\Phi_*\bigl(\Ga(\d f)\bigr)$ is an immersed special
Lagrangian $m$-fold in\/~$(M,J,\om,\Om)$.
\label{co5thm1}
\end{thm}

Its proof in \cite[\S 5]{Joyc5} is long and technical. The basic
idea is to write the equation for $\ti N=\Phi_*\bigl(\Ga(\d f)\bigr)$
to be special Lagrangian as a second-order nonlinear elliptic p.d.e.\
on $f$. Conditions (i)--(vii) ensure that this p.d.e.\ is close
to being linear, in a certain sense. We then solve the p.d.e.\ 
for $f$ by a series method, using facts about the solutions of
second-order linear elliptic~p.d.e.s.

On a first reading, Definition \ref{co5def2} and Theorem
\ref{co5thm1} may look like formidably technical abstract
nonsense. We now try to explain (informally, and
oversimplifying a bit) what the theorem does, and the
reasons behind its design.
\begin{itemize}
\setlength{\parsep}{0pt}
\setlength{\itemsep}{0pt}
\item The theorem says that given a Lagrangian $m$-fold
$N$ in $(M,J,\om,\Om)$ which is {\it close to special
Lagrangian} in a certain sense, we can deform $N$ to a
nearby SL $m$-fold $\ti N$ in $M$ by a small Hamiltonian
deformation.
\item A Lagrangian $m$-fold $N$ has a {\it phase function}
${\rm e}^{i\th}$, and $N$ is special Lagrangian with some
orientation if and only if $\sin\th\equiv 0$. Part (i) of
Theorem \ref{co5thm1} requires four norms of $\psi^m\sin\th$
to be small, so it forces $N$ to be {\it close to special
Lagrangian}.
\item We shall apply Theorem \ref{co5thm1} when $N$ is an
explicitly constructed desingularization of a singular SL
$m$-fold. The construction depends on $t\in(0,\de)$, the
{\it length scale} at which the singularities are resolved.
Thus, we actually construct a 1-{\it parameter family} of
Lagrangian $m$-folds $N^t$ for~$t\in(0,\de)$.

When $t$ is small, $N^t$ is {\it close to special Lagrangian},
but it is also {\it close to singular}, in that the metric
$h^t=g\vert_{N^t}$ on $N^t$ has large Riemann curvature $R(h^t)
=O(t^{-2})$, and small injectivity radius~$\de(h^t)=O(t)$.

So the theorem is set up using a real parameter $t>0$.
Parts (iv), (v) say that $(N,h)$ is {\it not too close
to singular}, and part (iii) that the geometry of $M$
near $N$ is {\it not too close to singular},
in terms of~$t$.
\item When $t$ is small $N$ is close to special Lagrangian
by (i), which is an {\it advantage}, but (iii)--(v) allow
$N$ to be close to singular, which is a {\it disadvantage}.
The proof is a delicate balancing act between these two
influences. The advantages win, and for small $t\le\ep$
we can deform $N$ to an SL $m$-fold~$\ti N$.
\item Roughly speaking, to solve the p.d.e.\ on $f$ we
need an {\it inverse $\De^{-1}$ of the Laplacian} $\De$ on
$N$, which should be {\it bounded independent of\/}~$t$.

However $\De 1=0$, so $\De$ is not invertible. Also,
in our applications $\De$ has a finite number of
{\it small positive eigenvalues} of size $O(t^{m-2})$,
so that if it existed $\De^{-1}$ would be $O(t^{2-m})$,
and not bounded independent of~$t$.

To get round this we introduce a vector space
$W\subset C^\iy(N)$ with $1\in W$, which approximates
the eigenspaces of $\De$ with small eigenvalues.
Part (vi) roughly says that $\De^{-1}=O(1)$ on
$W^\perp$. Part (vii) roughly says that
$\De^{-1}=O(t^{2-m})$ on the subspace of $W$
orthogonal to~1.

To control $W$-components of functions in the proof
although $\De^{-1}$ is large on $W$, we assume in (i) that
the $W$-component $\pi_\sW(\psi^m\sin\th)$ is very small.
\end{itemize}

\subsection{Theorems on desingularizing SL $m$-folds}
\label{co52}

In \cite[\S 6--\S 7]{Joyc5} we apply Theorem \ref{co5thm1}
to construct desingularizations of compact SL $m$-folds
$X$ in $(M,J,\om,\Om)$ with conical singularities. The
first main result is \cite[Th.~6.13]{Joyc5}, which
deals with the simplest case in which there are no
obstructions to desingularization.

\begin{thm} Suppose $(M,J,\om,\Om)$ is an almost Calabi--Yau
$m$-fold and\/ $X$ a compact SL\/ $m$-fold in $M$ with conical
singularities at\/ $x_1,\ldots,x_n$ and cones $C_1,\ldots,C_n$.
Let\/ $L_1,\ldots,L_n$ be Asymptotically Conical SL\/ $m$-folds
in $\C^m$ with cones $C_1,\ldots,C_n$ and rates $\la_1,\ldots,\la_n$.
Suppose $\la_i<0$ for $i=1,\ldots,n$, and\/ $X'=X\sm\{x_1,\ldots,x_n\}$
is connected.

Then there exists $\ep>0$ and a smooth family $\bigl\{
\smash{\ti N^t}:t\in(0,\ep]\bigr\}$ of compact, nonsingular
SL\/ $m$-folds in $(M,J,\om,\Om)$, such that\/ $\smash{\ti N^t}$
is constructed by gluing $tL_i$ into $X$ at\/ $x_i$ for
$i=1,\ldots,n$. In the sense of currents, $\smash{\ti N^t}\ra X$
as~$t\ra 0$.
\label{co5thm2}
\end{thm}

The second main result is \cite[Th.~7.10]{Joyc5}. It strengthens
Theorem \ref{co5thm2} by dropping the assumption that $X'$ is
connected. In doing this we encounter {\it topological
obstructions}, so that desingularizations $\smash{\ti N^t}$
exist only if an equation \eq{co5eq2} holds on topological
invariants of $X$ and~$L_1,\ldots,L_n$.

\begin{thm} Suppose $(M,J,\om,\Om)$ is an almost Calabi--Yau
$m$-fold and\/ $X$ a compact SL\/ $m$-fold in $M$ with conical
singularities at\/ $x_1,\ldots,x_n$ and cones $C_1,\ldots,C_n$.
Define $\psi:M\ra(0,\iy)$ as in \eq{co2eq3}. Let\/ $L_1,\ldots,L_n$
be Asymptotically Conical SL\/ $m$-folds in $\C^m$ with cones
$C_1,\ldots,C_n$ and rates $\la_1,\ldots,\la_n$. Suppose $\la_i<0$
for $i=1,\ldots,n$. Write $X'=X\sm\{x_1,\ldots,x_n\}$
and\/~$\Si_i=C_i\cap{\cal S}^{2m-1}$.

Set\/ $q=b^0(X')$, and let\/ $X_1',\ldots,X_q'$ be the connected
components of\/ $X'$. For $i=1,\ldots,n$ let\/ $l_i=b^0(\Si_i)$,
and let\/ $\Si_i^1,\ldots,\Si_i^{\smash{l_i}}$ be the connected
components of\/ $\Si_i$. Define $k(i,j)=1,\ldots,q$ by $\Up_i
\circ\vp_i\bigl(\Si_i^j\t(0,R')\bigr)\subset X'_{\smash{k(i,j)}}$
for $i=1,\ldots,n$ and $j=1,\ldots,l_i$. Suppose that
\e
\sum_{\substack{1\le i\le n, \; 1\le j\le l_i: \\
k(i,j)=k}}\psi(x_i)^mZ(L_i)\cdot[\Si_i^j\,]=0
\quad\text{for all\/ $k=1,\ldots,q$.}
\label{co5eq2}
\e

Suppose also that the compact\/ $m$-manifold\/ $N$ obtained by
gluing $L_i$ into $X'$ at\/ $x_i$ for $i=1,\ldots,n$ is connected.
A sufficient condition for this to hold is that\/ $X$ and\/ $L_i$
for $i=1,\ldots,n$ are connected.

Then there exists $\ep>0$ and a smooth family $\smash{\bigl\{
\ti N^t:t\in(0,\ep]\bigr\}}$ of compact, nonsingular SL\/
$m$-folds in $(M,J,\om,\Om)$ diffeomorphic to $N$, such
that\/ $\smash{\ti N^t}$ is constructed by gluing $tL_i$
into $X$ at\/ $x_i$ for $i=1,\ldots,n$. In the sense of
currents in Geometric Measure Theory, $\smash{\ti N^t}\ra X$
as~$t\ra 0$.
\label{co5thm3}
\end{thm}

When $X'$ is connected, so that $q=1$, it turns out that
\eq{co5eq2} holds automatically and Theorem \ref{co5thm3}
reduces to Theorem \ref{co5thm2}. Theorems \ref{co5thm2}
and \ref{co5thm3} are proved by the following method,
which will also be used in \S\ref{co6}--\S\ref{co8} below.

We shrink $L_i$ by a small factor $t>0$ to get $tL_i$, which
is also an AC SL $m$-fold in $\C^m$. Using the Lagrangian
neighbourhood results of \S\ref{co34} and \S\ref{co42} we
glue $tL_i$ into $X$ at $x_i$ using a partition of unity,
to get a Lagrangian $m$-fold $N^t$ in $(M,\om)$ for $t\in
(0,\de)$. We also glue the Lagrangian neighbourhoods of
$X'$ and $L_i$ together to get a Lagrangian neighbourhood
$\Phi_\sNt$ for~$N^t$.

We define vector spaces $W^t\subset C^\iy(N^t)$, using
spaces of bounded harmonic functions on $L_i$. We then show
that $N^t,W^t$ and $\Phi^t$ satisfy Definition \ref{co5def2}
and parts (i)--(vii) of Theorem \ref{co5thm1} for all
$t\in(0,\de)$, for some $\ka>1$ and $A_1,\ldots,A_8>0$
independent of $t$. Theorem \ref{co5thm1} then gives $\ep>0$
depending on $\ka,A_1,\ldots,A_8$ such that $N^t$ can be
deformed to an SL $m$-fold $\smash{\ti N^t}$ if~$t\le\ep$.

Here is how the obstruction equation \eq{co5eq2} arises.
For Theorem \ref{co5thm2} we take $W^t=\an{1}$, and then
$\pi_\sWt(\psi^m\sin\th^t)=0$ for topological reasons.
But for Theorem \ref{co5thm3} we have $W^t\cong\R^q$, and
$\pi_\sWt(\psi^m\sin\th^t)$ need not be zero. Calculation
shows that the leading contribution to $\pi_\sWt(\psi^m
\sin\th^t)$ is $O(t^m)$, and is proportional to the left
hand side of \eq{co5eq2}. For part (i) of Theorem
\ref{co5thm2} we need $\lnm{\pi_\sWt(\psi^m\sin\th^t)}{1}\le
A_2t^{\ka+m-1}$, and this holds if and only if \eq{co5eq2} does.

\section{Desingularizing when $Y(L_i)\ne 0$}
\label{co6}

In Theorems \ref{co5thm2} and \ref{co5thm3} we desingularized an
SL $m$-fold $X$ with conical singularities using AC SL $m$-folds
$L_i$ with rates $\la_i$ for $i=1,\ldots,n$, where we assumed that
$\la_i<0$, so that $Y(L_i)=0$ by Proposition \ref{co4prop}. We now
explain how to relax this to allow $\la_i\le 0$ and $Y(L_i)\ne 0$.
As in \cite[\S 6.4]{Joyc6} there are many examples of AC SL
$m$-folds $L$ with rate 0 and $Y(L)\ne 0$, so this significantly
increases the scope of the main result.

Allowing $\la_i=0$ complicates the proofs in two main ways.
Firstly, gluing $tL_i$ into $X$ at $x_i$ to make $N^t$ is no
longer a local matter. Instead, there is a {\it global\/}
condition for $N^t$ to exist as a {\it Lagrangian} $m$-fold, that
$\bigl(Y(L_1),\ldots,Y(L_n)\bigr)$ should lie in the image of the
map $H^1(X',\R)\ra\bigoplus_{i=1}^nH^1(\Si_i,\R)$ in~\eq{co3eq3}.

Secondly, the modifications to the definition of $N^t$ introduce
extra error terms in $\Im\Om\vert_{N^t}$, which contribute $O(t^4)$
to $\lnm{\psi^m\sin\th^t}{2m/(m+2)}$. But for part (i) of Theorem
\ref{co5thm1} to hold we need $\lnm{\psi^m\sin\th^t}{2m/(m+2)}=
O(t^{\ka+m/2})$ for $\ka>1$, so we must assume $m<6$. Therefore
our main result, Theorem \ref{co6thm5}, holds only in dimensions
$m=3,4,5$, rather than all $m\ge 3$ as in Theorem~\ref{co5thm2}.

\subsection{Setting up the problem}
\label{co61}

We shall consider the following situation, the analogue
of~\cite[Def.~6.1]{Joyc5}.

\begin{dfn} Let $(M,J,\om,\Om)$ be an almost Calabi--Yau
$m$-fold with metric $g$, and define $\psi:M\ra(0,\iy)$
as in \eq{co2eq3}. Let $X$ be a compact SL $m$-fold in
$M$ with conical singularities at $x_1,\ldots,x_n$ with
identifications $\up_i$, cones $C_i$ and rates $\mu_i\in(2,3)$,
as in Definition \ref{co3def3}. Let $L_1,\ldots,L_n$ be AC SL
$m$-folds in $\C^m$ as in Definition \ref{co4def1}, where $L_i$
has cone $C_i$. Define $\Si_i=C_i\cap{\cal S}^{2m-1}$ and
$Y(L_i)\in H^1(\Si_i,\R)$ and $Z(L_i)\in H^{m-1}(\Si_i,\R)$
for $i=1,\ldots,n$ as in Definition~\ref{co4def2}.

Set $q=b^0(X')$, so that $X'$ has $q$ connected components,
and number them $X_1',\ldots,X_q'$. For $i=1,\ldots,n$ let
$l_i=b^0(\Si_i)$, so that $\Si_i$ has $l_i$ connected
components, and number them $\Si_i^1,\ldots,\Si_i^{\smash{l_i}}$.
If $\Up_i,\vp_i,S_i$ are as in Definition \ref{co3def3}, then
$\Up_i\circ\vp_i$ is a diffeomorphism $\Si_i\t(0,R')\ra S_i
\subset X'$. For each $j=1,\ldots,l_i$, $\Up_i\circ\vp_i\bigl(
\Si_i^j\t(0,R')\bigr)$ is a connected subset of $X'$, and so
lies in exactly one of the $X_k'$ for $k=1,\ldots,q$. Define
numbers $k(i,j)=1,\ldots,q$ for $i=1,\ldots,n$ and $j=1,\ldots,l_i$.
by $\Up_i\circ\vp_i(\Si_i^j\t(0,R'))\subset X'_{\smash{k(i,j)}}$.
Suppose that:
\begin{itemize}
\setlength{\parsep}{0pt}
\setlength{\itemsep}{0pt}
\item[(i)] The dimension $m$ satisfies $2<m<6$, 
\item[(ii)] The AC SL $m$-fold $L_i$ has rate $\la_i\le 0$
for $i=1,\ldots,n$, and
\item[(iii)] There exists $\varrho\in H^1(X',\R)$ such that
$\bigl(Y(L_1),\ldots,Y(L_n)\bigr)$ is the image of
$\varrho$ under the map $H^1(X',\R)\ra\bigoplus_{i=1}^n
H^1(\Si_i,\R)$ in~\eq{co3eq3}.
\item[(iv)] $\displaystyle\sum_{\substack{1\le i\le n, \;
1\le j\le l_i: \\ k(i,j)=k}}\psi(x_i)^mZ(L_i)\cdot[\Si_i^j\,]=0$
for all $k=1,\ldots,q$.
\item[(v)] Let $N$ be the compact $m$-manifold obtained
by gluing $L_i$ into $X'$ at $x_i$ for $i=1,\ldots,n$ in the
obvious way. Suppose $N$ is connected. A sufficient
condition for this to hold is that $X$ and $L_i$ for
$i=1,\ldots,n$ are connected.
\end{itemize}

We use the following notation:
\begin{itemize}
\setlength{\parsep}{0pt}
\setlength{\itemsep}{0pt}
\item Let $R,B_R,X'$ and $\iota_i,\Up_i$ for $i=1,\ldots,n$ be
as in Definition~\ref{co3def3}.
\item Let $\ze$ and $U_\sCi,\Phi_\sCi$ for
$i=1,\ldots,n$ be as in Theorem~\ref{co3thm3}.
\item Let $R',K$ and $\phi_i,\eta_i,\eta_i^1,\eta_i^2,S_i$
for $i=1,\ldots,n$ be as in Theorem~\ref{co3thm4}.
\item Let $U_\sXp,\Phi_\sXp$ be as in Theorem~\ref{co3thm6}.
\item Let $A_i$ be as in Theorem \ref{co3thm5} for $i=1,\ldots,n$,
so that~$\eta_i=\d A_i$.
\item For $i=1,\ldots,n$ let $\ga_i$ be the unique 1-form on $\Si_i$
with $\d\ga_i=\d^*\ga_i=0$ and $[\ga_i]=Y(L_i)$ in $H^1(\Si_i,\R)$,
which exists by Hodge theory.

Let $\pi_i:\Si_i\t(0,\iy)\ra\Si_i$ be the projection, so that
$\pi_i^*(\ga_i)$ is a 1-form on $\Si_i\t(0,\iy)$ with
$\bmd{\pi_i^*(\ga_i)}=O(r^{-1})$ and $\d\bigl(\pi_i^*(\ga_i)\bigr)
=\d^*\bigl(\pi_i^*(\ga_i)\bigr)=0$, computing $\md{\,.\,}$ and
$\d^*$ using the cone metric $\iota_i^*(g')$ on~$\Si_i\t(0,\iy)$.
\item Let $Y_\sXp$ be as in \eq{co3eq6}, and $\al\in Y_\sXp$
the unique element with $\pi(\al)=\varrho$. Then $\d\al=
\d^*(\psi^m\al)=0$. For $i=1,\ldots,n$ let $T_i\in C^\iy\bigl(
\Si_i\t(0,R')\bigr)$ be as in part (a) of Theorem \ref{co3thm1}.
Choose $\nu_i\in(0,\mu_i-2)$ with $(0,\nu_i]\cap\D_\sSii=\emptyset$.
Then part (a) of Theorem \ref{co3thm1} shows that
\end{itemize}
\ea
(\Up_i\circ\phi_i)^*(\al)=\pi_i^*(\ga_i)+\d T_i
&\quad\text{on $\Si_i\t(0,R')$ for $i=1,\ldots,n$,}
\label{co6eq1}\\
\text{and}\quad
\na^kT_i(\si,r)=O(r^{\nu_i-k})
&\quad\text{as $r\ra 0$, for all $k\ge 0$.}
\label{co6eq2}
\ea
\begin{itemize}
\item Apply Theorem \ref{co4thm1} to $L_i$ with $\ze,U_\sCi,
\Phi_\sCi$ as above, for $i=1,\ldots,n$. Let $T>0$ be as in
the theorem, the same for all $i$. Let the subset $K_i\subset L_i$,
the diffeomorphism $\vp_i:\Si_i\t(T,\iy)\ra L_i\sm K_i$ and the
1-form $\chi_i$ on $\Si_i\t(T,\iy)$ with components $\chi_i^1,
\chi_i^2$ be as in Theorem~\ref{co4thm1}.
\item Let $U_\sLi,\Phi_\sLi$ be as in Theorem \ref{co4thm4}
for~$i=1,\ldots,n$.
\item Let $E_i\in C^\iy\bigl(\Si_i\t(T,\iy)\bigr)$ be as in
Theorem \ref{co4thm2} for $i=1,\ldots,n$. Fix $\la\in
\bigl(\max(-2,2-m),\ha(2-m)\bigr)$. (This interval is nonempty
as $2<m<6$.) Then Theorem \ref{co4thm2} and part (b) of
Theorem \ref{co4thm3} show that
\end{itemize}
\ea
\chi_i=\pi_i^*(\ga_i)+\d E_i
&\quad\text{on $\Si_i\t(T,\iy)$ for $i=1,\ldots,n$}
\label{co6eq3}\\
\text{and}\quad
\na^kE_i(\si,r)=O(r^{\la-k})
&\quad\text{as $r\ra\iy$ for all\/ $k\ge 0$.}
\label{co6eq4}
\ea
By \eq{co6eq1} and \eq{co6eq3} we see that the 1-forms $\al$ on
$X'$ and $\chi_i$ on $L_i$ both have leading term $\pi_i^*(\ga_i)$
in their asymptotic expansion on the cone $C_i$. The construction
below will work by matching up these terms $\pi_i^*(\ga_i)$, so
that we do not have to taper them off to zero.
\label{co6def1}
\end{dfn}

Here are some remarks on conditions (i)--(v) above.
\begin{itemize}
\setlength{\parsep}{0pt}
\setlength{\itemsep}{0pt}
\item[(i)] It will turn out in \S\ref{co62} that defining the
Lagrangian $m$-folds $N^t$ when $Y(L_i)\ne 0$ introduces
$O(t^4)$ error terms in $\lnm{\psi^m\sin\th^t}{2m/(m+2)}$.
However, for part (i) of Theorem \ref{co5thm1} to hold we
need $\lnm{\psi^m\sin\th^t}{2m/(m+2)}=O(t^{\ka+m/2})$ for~$\ka>1$.

Thus we need $\ka+m/2\le 4$ and $\ka>1$, giving $m<6$, which
is why we suppose $m<6$ in (i) above. We also use $2<m<6$ in
choosing $\la\in\bigl(\max(-2,2-m),\ha(2-m)\bigr)$. With some
more work on the definition of $N^t$, the result can probably
be extended to the case~$m=6$.
\item[(ii)] The point of this section is to relax the assumption
$\la_i<0$ in Theorems \ref{co5thm2} and \ref{co5thm3}, and so
allow $Y(L_i)\ne 0$. We suppose $\la_i\le 0$ so that part (b)
of Theorem \ref{co4thm3} applies. Although going from $\la_i<0$
to $\la_i\le 0$ may not seem like much of an improvement, in
fact as in \cite[\S 6.4]{Joyc6} there are many examples of AC
SL $m$-folds $L$ with rate exactly $\la=0$, so allowing $\la=0$
will make our results much more useful.
\item[(iii)] The condition that $\bigl(Y(L_1),\ldots,Y(L_n)\bigr)$
lies in the image of $H^1(X',\R)\ra\bigoplus_{i=1}^nH^1(\Si_i,\R)$
is necessary for the existence of any {\it Lagrangian} $m$-fold
$N^t$ made by gluing $tL_i$ into $X$ at $x_i$ for $i=1,\ldots,n$.
So it is clearly also necessary for the existence of an SL
$m$-fold $\smash{\ti N^t}$ made in the same way.
\item[(iv)] This condition was introduced in \cite[Def.~7.1]{Joyc5},
to deal with analytic problems in desingularizing $X$ when $X'$ is
not connected. As in \cite[\S 7.4]{Joyc5}, when $X'$ is connected
$q=1$ and (iv) holds automatically since~$Z(L_i)\cdot[\Si_i]=0$.
\item[(v)] We assume this for simplicity. If $N$ is not
connected then we can apply Theorem \ref{co6thm5} below to each
component of $N$ separately.
\end{itemize}

Here is the analogue of \cite[Def.s 6.2 \& 6.3]{Joyc5},
constructing $N^t$ for~$t\in(0,\de)$.

\begin{dfn} Choose a smooth, increasing $F:(0,\iy)\ra[0,1]$
with $F(r)\equiv 0$ for $r\in(0,1)$ and $F(r)\equiv 1$ for
$r>2$, and write $F'=\d F/\d r$. Let $\tau$ satisfy
\e
0<\max\Bigl\{\frac{m}{m+1},\frac{m+2}{2\mu_i+m-2},
\frac{m-2}{2\nu_i+m-2}:i=1,\ldots,n\Bigr\}<\tau<1,
\label{co6eq5}
\e
which is possible as $\mu_i>2$, $\nu_i>0$ and $m>2$. For
$i=1,\ldots,n$ and small enough $t>0$, define $P_i^t=
\Up_i(tK_i)$. This is well-defined, Lagrangian in
$(M,\om)$, and diffeomorphic to~$K_i$.

For $i=1,\ldots,n$ and $t>0$ with $tT<t^\tau<2t^\tau<R'$ define
a 1-form $\xi_i^t$ on $\Si_i\t(tT,R')$ by
\e
\begin{split}
\xi_i^t(\si,r)&=\d\bigl[F(t^{-\tau}r)A_i(\si,r)+
t^2(1-F(t^{-\tau}r))E_i(\si,t^{-1}r)\bigr]\\
&\quad+t^2\pi_i^*(\ga_i)+t^2\d\bigl[F(t^{-\tau}r)T_i(\si,r)\bigr]\\
&=F(t^{-\tau}r)\eta_i(\si,r)+
t^{-\tau}F'(t^{-\tau}r)A_i(\si,r)\d r\\
&\quad+t^2(1\!-\!F(t^{-\tau}r))\chi_i(\si,t^{-1}r)-
t^{2-\tau}F'(t^{-\tau}r)E_i(\si,t^{-1}r)\d r\\
&\quad+t^2F(t^{-\tau}r)(\Up_i\circ\phi_i)^*(\al)(\si,r)
+t^{2-\tau}F'(t^{-\tau}r)T_i(\si,r)\d r,
\end{split}
\label{co6eq6}
\e
by \eq{co6eq1} and \eq{co6eq3}. Let $\xi_i^{1,t},\xi_i^{2,t}$ be the
components of $\xi_i^t$ in $T^*\Si$ and $\R$. Then when $r\ge 2t^\tau$
we have $\xi_i^t\equiv\eta_i+t^2(\Up_i\circ\phi_i)^*(\al)$, and when
$r\le t^\tau$ we have $\xi_i^t(\si,r)=t^2\chi_i(\si,t^{-1}r)$. Thus
$\xi_i^t$ is a closed 1-form which interpolates between
$\eta_i+t^2(\Up_i\circ\phi_i)^*(\al)$ near $r=R'$ and
$t^2\chi_i(\si,t^{-1}r)$ near~$r=tT$.

Choose $\de\in(0,1]$ with $\de T\le\de^\tau<2\de^\tau\le R'$, $\de
K_i\!\subset\!B_R\!\subset\!\C^m$ and $\bmd{\xi_i^t(\si,r)}<\ze r$
on $\Si_i\t(tT,R')$ for all $i=1,\ldots,n$ and $t\in(0,\de)$. As in
\cite[Def.~6.2]{Joyc5}, this is possible. Define $\Xi_i^t:\Si_i\t
(tT,R')\ra M$ by
\begin{equation*}
\Xi_i^t(\si,r)=\Up_i\circ\Phi_\sCi\bigl(\si,r,
\xi_i^{1,t}(\si,r),\xi_i^{2,t}(\si,r)\bigr)
\end{equation*}
for $i=1,\ldots,n$ and $t\in(0,\de)$. Making $R'$ smaller if
necessary, this is well-defined, and an embedding. Define
$Q_i^t=\Xi_i^t\bigl(\Si_i\t(tT,R')\bigr)$ for $i=1,\ldots,n$ and
$t\in(0,\de)$. As $\Up_i^*(\om)=\om'$, $\Phi_\sCi^*(\om')=\hat\om$
and $\xi_i^t$ is a closed 1-form we see that $(\Xi_i^t)^*(\om)\equiv 0$.
Thus $Q_i^t$ is {\it Lagrangian} in $(M,\om)$, and is a noncompact
embedded submanifold diffeomorphic to~$\Si_i\t(tT,R')$.

Let $\Ga(t^2\al)$ be the graph of the 1-form $t^2\al$ in $T^*X'$.
Then $\Ga(t^2\al)\cap\pi^*(K)\subset T^*K$ is the graph of
$t^2\al\vert_K$. Recall that $U_\sXp$ is an open neighbourhood
of the zero section in $T^*X'$. By compactness of $K$, making
$\de$ smaller if necessary, we can suppose that $\Ga(t^2\al)
\cap\pi^*(K)\subset U_\sXp$ for all $t\in(0,\de)$. Define
\e
K^t=\Phi_\sXp\bigl(\Ga(t^2\al)\cap\pi^*(K)\bigr)
\quad\text{for $t\in(0,\de)$.}
\label{co6eq7}
\e
Then $K^t$ is a submanifold of $M$ with boundary, diffeomorphic
to $K$. As $\al$ is closed $\Ga(t^2\al)$ is Lagrangian in
$(T^*X',\hat\om)$, and $\Phi_\sXp^*(\om)=\hat\om$, so $K^t$
is {\it Lagrangian} in $(M,\om)$. For $t\in(0,\de)$, define
$N^t$ to be the disjoint union of $K^t$, $P_1^t,\ldots,P_n^t$
and~$Q_1^t,\ldots,Q_n^t$.

Then $N^t$ is {\it Lagrangian} in $(M,\om)$, as $K^t,P_i^t$ and
$Q_i^t$ are. Moreover, $N^t$ is a compact, smooth submanifold
of $M$ {\it without boundary}. To see this, note that
$\xi_i^t(\si,r)=t^2\chi_i(\si,t^{-1}r)$ on $\Si_i\!\t\!
(tT,t^\tau]$, and so
\begin{align*}
\Xi_i^t(\si,r)&=\Up_i\circ\Phi_\sCi\bigl(\si,r,
t^2\chi_i^1(\si,t^{-1}r),t\chi_i^2(\si,t^{-1}r)\bigr)\\
&=\Up_i\bigl(t\,\Phi_\sCi\bigl(\si,t^{-1}r,
\chi_i^1(\si,t^{-1}r),\chi_i^2(\si,t^{-1}r)\bigr)\bigr)
=\Up_i\bigl(t\,\vp_i(\si,t^{-1}r)\bigr)
\end{align*}
on $\Si_i\!\t\!(tT,t^\tau]$, using \eq{co4eq3} and the dilation
equivariance of $\Phi_{\sst C_i}$ in Theorem~\ref{co3thm3}.

Thus the end $\Xi_i^t\bigl(\Si_i\!\t\!(tT,t^\tau]\bigr)$ of $Q_i^t$
is $\Up_i\bigl(t\,\vp_i(\Si_i\!\t\!(T,t^{\tau-1}])\bigr)\subset
\Up_i(tL_i)$, and as $\vp_i\bigl(\Si_i\!\t\!(T,t^{\tau-1}]\bigr)
\subset L_i$ joins smoothly onto $K_i\subset L_i$ we see that
$\Xi_i^t\bigl(\Si_i\!\t\!(tT,t^\tau]\bigr)\subset Q_i^t$ joins
smoothly onto $P_i^t=\Up_i(tK_i)$. Similarly, the boundary
$\pd K^t$ is the disjoint union of pieces $\Si_i$ for
$i=1,\ldots,n$ which join smoothly onto $Q_i^t\cong\Si_i\!\t\!
(tT,R')$ at the $\Si_i\!\t\!\{R'\}$ end, as $K^t$ is the
graph of $t^2\al$ over $K$ and $\xi_i^t\equiv\eta_i+t^2(
\Up_i\circ\phi_i)^*(\al)$ on~$\Si_i\t[2t^\tau,R')$.

As $X',L_i$ are SL $m$-folds they are oriented, and
$N^t$ is made by gluing $X',L_1,\ldots,L_n$ together
in an orientation-preserving way, so $N^t$ is also
oriented. Let $h^t$ be the restriction of $g$ to $N^t$
for $t\in(0,\de)$, so that $(N^t,h^t)$ is a compact
Riemannian manifold, and let $\d V^t$ be the induced
volume form on $N^t$. Then $\Om\vert_{N^t}=\psi^m{\rm
e}^{\smash{i\th^t}}\,\d V^t$ for some phase function
${\rm e}^{\smash{i\th^t}}$ on $N^t$. Write $\ve^t=
\psi^m\sin\th^t$, so that $\Im\Om\vert_{N^t}=\ve^t\,
\d V^t$ for~$t\in(0,\de)$.
\label{co6def2}
\end{dfn}

From \eq{co6eq6} we can explain the reason for condition (iii)
in Definition \ref{co6def1}. The problem is that we need to
define a {\it closed\/} 1-form $\xi_i^t$ on $\Si_i\t(tT,R')$
which interpolates between $t^2\chi_i(\si,t^{-1}r)$ near $r=tT$
and $\eta_i+t^2(\Up_i\circ\phi_i)^*(\al)$ near $r=R'$. This is
possible if and only if $t^2\chi_i(\si,t^{-1}r)$ and $\eta_i+
t^2(\Up_i\circ\phi_i)^*(\al)$ have the {\it same cohomology
class} in~$H^1\bigl(\Si_i\t(tT,R'),\R)\cong H^1(\Si_i,\R)$.

The cohomology class of $\chi_i(\si,t^{-1}r)$ is $Y(L_i)\in
H^1(\Si_i,\R)$ by Theorem \ref{co4thm2}, and $\eta_i=\d A_i$
is exact. Hence the cohomology class of $(\Up_i\circ\phi_i)^*
(\al)$ in $H^1(\Si_i,\R)$ must be $Y(L_i)$. That is, the
cohomology class $[\al]=\varrho\in H^1(X',\R)$ must have
image $Y(L_i)$ in $H^1(\Si_i,\R)$ for all $i=1,\ldots,n$
under the map $H^1(X',\R)\ra\bigoplus_{i=1}^nH^1(\Si_i,\R)$
in \eq{co3eq3}, which gives (iii). This is a {\it necessary
condition} for $N^t$ to exist as a Lagrangian $m$-fold.

Following \cite[Def.~7.2]{Joyc5} we define a vector subspace
$W^t\subset C^\iy(N^t)$. This will be $W$ in Definition
\ref{co5def2}. It is an approximation to the eigenspaces
of the Laplacian $\De$ on $N^t$ with small eigenvalues,
and is a tool to deal with some analytic problems when
$X'$ is not connected.

\begin{dfn} We work in the situation of Definitions \ref{co6def1}
and \ref{co6def2}. For $i=1,\ldots,n$ apply Theorem \ref{co4thm5}
to the AC SL $m$-fold $L_i$ in $\C^m$, using the numbering $\Si_i^j$
chosen in Definition \ref{co6def1} for the connected components of
$\Si_i$. This gives a vector space $V_i$ of bounded harmonic
functions on $L_i$ with $\dim V_i=l_i$. For each ${\bf c}_i=
(c_i^1,\ldots,c_i^{l_i})\in\R^{l_i}$ there exists a unique
$v^{{\bf c}_i}_{\smash{i}}\in V_i$ with
\e
\na^k\bigl(\vp_i^*(v^{{\bf c}_i}_{\smash{i}})-c^j_i\,\bigr)=
O\bigl(\md{{\bf c}_i}r^{\be-k}\bigr)
\quad\text{on $\Si_i^j\t(T,\iy)$ as $r\ra\iy$,}
\label{co6eq8}
\e
for all $i=1,\ldots,n$, $j=1,\ldots,l_i$, $k\ge 0$ and~$\be\in(2-m,0)$.

We shall define a vector subspace $W^t\subset C^\iy(N^t)$
for $t\in(0,\de)$, with an isomorphism $W^t\cong\R^q$. Fix
${\bf d}=(d_1,\ldots,d_q)\in\R^q$, and set $c_i^j=d_{k(i,j)}$
for $i=1,\ldots,n$ and $j=1,\ldots,l_i$. Let ${\bf c}_i=
(c_i^1,\ldots,c_i^{l_i})$. This defines vectors ${\bf c}_i\in
\R^{l_i}$ for $i=1,\ldots,n$, which depend linearly on $\bf d$.
Hence we have harmonic functions $v^{{\bf c}_i}_{\smash{i}}\in
V_i\subset C^\iy(L_i)$, which also depend linearly on~$\bf d$.

Let $F:(0,\iy)\ra[0,1]$ and $\tau\in(0,1)$ be as in Definition
\ref{co6def2}. Make $\de>0$ smaller if necessary so that
$tT<\ha t^\tau$ for all $t\in(0,\de)$. For $t\in(0,\de)$,
define a function $w_{\bf d}^t\in C^\iy(N^t)$ as follows:
\begin{itemize}
\item[(i)] The subset $K^t\subset N^t$ is diffeomorphic to $K$,
and so has $q$ connected components $K_k^t$ diffeomorphic to
$K\cap X_k'$ for $k=1,\ldots,q$. Define $w_{\bf d}^t$ on
$K^t$ by $w_{\bf d}^t\equiv d_k$ on $K^t_k$ for~$k=1,\ldots,q$.
\item[(ii)] Define $w_{\bf d}^t$ on $P_i^t\subset N^t$ by
$(\Up_i\circ t\circ\vp_i)^*(w_{\bf d}^t)\equiv
v^{{\bf c}_i}_{\smash{i}}$ on $K_i$ for~$i=1,\ldots,n$.
\item[(iii)] Define $w_{\bf d}^t$ on $Q_i^t\subset N^t$ by
\e
(\Xi_i^t)^*(w_{\bf d}^t)(\si,r)=\bigl(1-F(2t^{-\tau}r)\bigr)
\vp_i^*(v^{{\bf c}_i}_{\smash{i}})(\si,t^{-1}r)+F(2t^{-\tau}r)c_i^j
\label{co6eq9}
\e
on $\Si_i^j\t(tT,R')$, for $i=1,\ldots,n$ and~$j=1,\ldots,l_i$.
\end{itemize}
It is easy to see that $w_{\bf d}^t$ is smooth over the joins between
$P_i^t,Q_i^t$ and $K^t$, so $w_{\bf d}^t\in C^\iy(N^t)$. Also
$w_{\bf d}^t$ is linear in ${\bf d}$, as $v^{{\bf c}_i}_{\smash{i}}$
is. Thus $W^t=\{w_{\bf d}^t:{\bf d}\in\R^q\}$ is a vector subspace of
$C^\iy(N^t)$ isomorphic to $\R^q$, for all~$t\in(0,\de)$.

If ${\bf d}=(1,\ldots,1)$ then $c_i^j\equiv 1$, so ${\bf c}_i=(1,
\ldots,1)$ for $i=1,\ldots,n$, and thus $v^{{\bf c}_i}_{\smash{i}}
\equiv 1$ for $i=1,\ldots,n$ by Theorem \ref{co4thm5}. Therefore
$w_{\smash{(1,\ldots,1)}}^{\smash{t}}\equiv 1$ by (i)--(iii)
above, and $1\in W^t$ for all $t\in(0,\de)$. This corresponds
to the condition $1\in W$ in Definition \ref{co5def2}. Define
$\pi_\sWt:L^2(N^t)\ra W^t$ to be the projection onto
$W^t$ using the $L^2$-inner product, as for $\pi_{\sst W}$
in Definition~\ref{co5def2}.
\label{co6def3}
\end{dfn}

Note that if $X'$ is connected, so that $q=1$, then~$W^t=\an{1}$.

\subsection{Estimating $\Im\Om\vert_{N^t}$}
\label{co62}

We now prove estimates for $\Im\Om\vert_{N^t}=\ve^t\,\d V^t$,
following \cite[\S 6.2]{Joyc5}. First we compute bounds for
$\ve^t$ at each point, as in~\cite[Prop.~6.4]{Joyc5}.

\begin{prop} In the situation above, making $\de>0$ smaller if
necessary, there exists $C>0$ such that for all\/ $t\in(0,\de)$ we have
\ea
\bmd{(\Xi_i^t)^*(\ve^t)}(\si,r)&\le
\begin{cases}
Cr, & r\in(tT,t^\tau],\\
\begin{aligned}
&\!Ct^{4-4\tau}+Ct^{\tau(\mu_i-2)}+\\
&\!Ct^{(1-\tau)(2-\la)}+Ct^{2+\tau(\nu_i-2)},
\quad   
\end{aligned}
& r\in(t^\tau,2t^\tau),\\
Ct^4r^{-4}, & r\in[2t^\tau,R'),
\end{cases}
\label{co6eq10}\\
\bmd{(\Xi_i^t)^*(\d\ve^t)}(\si,r)&\le
\begin{cases}
C, & r\in(tT,t^\tau],\\
\begin{aligned}
&\!Ct^{4-5\tau}+Ct^{\tau(\mu_i-3)}+\\
&\!Ct^{(1-\tau)(2-\la)-\tau}+Ct^{2+\tau(\nu_i-3)},
\end{aligned}
& r\in(t^\tau,2t^\tau),\\
Ct^4r^{-5}, & r\in[2t^\tau,R'),
\end{cases}
\label{co6eq11}\\
\md{\ve^t}&\le Ct^4,\quad \md{\d\ve^t}\le Ct^4\quad\text{on $K^t$,}
\label{co6eq12}\\
\text{and}\quad
\md{\ve^t}&\le Ct,\quad
\md{\d\ve^t}\le C
\quad\text{on $P_i^t$ for all\/ $i=1,\ldots,n$.}
\label{co6eq13}
\ea
Here $\md{\,.\,}$ is computed using $(\Xi_i^t)^*(h^t)$ in
\eq{co6eq11} and\/ $h^t$ in \eq{co6eq12} and\/~\eq{co6eq13}.
\label{co6prop1}
\end{prop}

\begin{proof} As $\Up_i^*(\Im\Om)$ is a smooth $m$-form on $B_R$
and $\Up_i^*(\Im\Om)\vert_0=\up_i^*(\Im\Om)=\psi(x_i)^m\Im\Om'$
by Definition \ref{co3def3}, we see that $\Up_i^*(\Im\Om)=
\psi(x_i)^m\Im\Om'+O(r)$ on $B_R$, by Taylor's Theorem. Since
$tL_i$ is special Lagrangian in $\C^m$ we have $\Im\Om'
\vert_{tL_i}\equiv 0$. Thus
\e
\bmd{\Up_i^*(\Im\Om)\vert_{tL_i}}=O(r)\quad\text{on $tL_i\cap B_R$,}
\label{co6eq14}
\e
computing $\md{\,.\,}$ using the metric $\Up_i^*(g)$ on $B_R$,
restricted to~$tL_i$.

Now $N^t$ coincides with $\Up_i(tL_i)$ on $P_i^t$ and
$\Xi_i^t\bigl(\Si_i\t(tT,t^\tau]\bigr)$, so
$\ve^t\d V^t=\Im\Om\vert_{\Up_i(tL_i)}$ on these regions.
As $h^t$ is the restriction of $g$ to $N^t$ we have
$\md{\d V^t}=1$, computing $\md{\,.\,}$ using $g$, so
\e
\bmd{\Up_i^*(\ve^t)}=\bmd{\Up_i^*(\Im\Om)\vert_{tL_i}}
\quad\text{on $t\bigl(K\cup\vp_i(\Si_i\t(T,t^{\tau-1}]\bigr)
\subset tL_i\cap B_R$.}
\label{co6eq15}
\e

Combining \eq{co6eq14} and \eq{co6eq15} gives $\bmd{\ve^t}=
O\bigl((\Up_i)_*(r)\bigr)$ on $P_i^t$ and $\Xi_i^t\bigl(
\Si_i\t(tT,t^\tau]\bigr)$. As $(\Up_i)_*(r)=O(t)$ on
$P_i^t$, we see that
\e
\bmd{(\Xi_i^t)^*(\ve^t)}(\si,r)=O(r)
\quad\text{for $r\in(tT,t^\tau]$, and}\quad
\md{\ve^t}=O(t)\quad\text{on $P_i^t$.}
\label{co6eq16}
\e
A similar argument for the derivative $\d\ve^t$ gives
\e
\bmd{(\Xi_i^t)^*(\d\ve^t)}(\si,r)=O(1)
\quad\text{for $r\in(tT,t^\tau]$, and}\quad
\md{\d\ve^t}=O(1)\quad\text{on $P_i^t$.}
\label{co6eq17}
\e

In \cite[Prop.~2.10]{Joyc4} we show that if $N$ is a compact
SL $m$-fold in $M$ and $\ti N$ is a nearby Lagrangian $m$-fold
written as the graph of a $C^1$ small closed 1-form $\al$ on
$N$ using Theorem \ref{co3thm2}, and $\Im\Om\vert_{\ti N}=
\psi^m\sin\th\,\d V$, then
\e
\psi^m\sin\th=-\d^*(\psi^m\al)+O(\ms{\al})+O(\ms{\na\al})
\quad\text{when $\md{\al},\md{\na\al}$ are small.}
\label{co6eq18}
\e
This is extended in \cite[Prop.~6.3]{Joyc4} to the case in
which $X$ is a compact SL $m$-fold with conical singularities,
and $\al$ a small 1-form on~$X'$.

Now on $K^t$ and the annuli $\Xi_i^t\bigl(\Si_i\t[2t^\tau,R')\bigr)$
for $i=1,\ldots,n$, $N^t$ is the graph of a small closed 1-form
$t^2\al$ on $X'$. Abusing notation, identify $K^t$ and $\Xi_i^t
\bigl(\Si_i\t[2t^\tau,R')\bigr)$ with the corresponding regions
in $X'$, so that $\psi,\al,\rho$ make sense on these regions in
$N^t$. Then \cite[Prop.~6.3]{Joyc4} shows that
\e
\ve^t=\psi^m\sin\th^t=-\d^*(\psi^mt^2\al)+O(\rho^{-2}\ms{t^2\al})
+O(\ms{t^2\na\al})
\label{co6eq19}
\e
on $K^t$ and $\Xi_i^t\bigl(\Si_i\t[2t^\tau,R')\bigr)$ when
$\rho^{-1}\md{t^2\al},\md{t^2\na\al}$ are small. Here
$\rho:X'\ra(0,1]$ is a {\it radius function}, as in
Definition~\ref{co3def4}.

Since $\md{\na^k\al}=\rho^{-1-k}$ by \eq{co3eq6} as $\al\in Y_\sXp$,
we see that $\rho^{-1}\md{t^2\al},\md{t^2\na\al}$ are $O(t^{2-2\tau})$
on $K^t$ and $\Xi_i^t\bigl(\Si_i\t[2t^\tau,R')\bigr)$, so \eq{co6eq19}
holds, giving
\e
\ve^t=O(t^4)\;\>\text{on $K^t$, and}\;\>
\bmd{(\Xi_i^t)^*(\ve^t)}=O(t^4r^{-4})\;\>
\text{on $\Xi_i^t\bigl(\Si_i\t[2t^\tau,R')\bigr)$,}
\label{co6eq20}
\e
as $\d^*(\psi^m\al)=0$. By a similar argument for derivatives we obtain
\e
\d\ve^t=O(t^4)\;\>\text{on $K^t$, and}\;\>
\bmd{(\Xi_i^t)^*(\d\ve^t)}=O(t^4r^{-5})\;\>
\text{on $\Xi_i^t\bigl(\Si_i\t[2t^\tau,R')\bigr)$.}
\label{co6eq21}
\e

On the annuli $\Xi_i^t\bigl(\Si_i\t(t^\tau,2t^\tau)\bigr)$ we can
apply the same argument, but we have to be rather more careful.
Here $N^t$ is not the graph of $t^2\al$ over $X'$, but instead
the graph of the 1-form $(\Up_i\circ\phi_i)_*(\xi_i^t-\eta_i)$ on
the corresponding annulus in $X'$. From \eq{co6eq6} we find that
\e
\begin{split}
(\xi_i^t&-\eta_i)(\si,r)=t^2(\Up_i\circ\phi_i)^*(\al)(\si,r)\\
&+\d\bigl[\bigl(1-F(t^{-\tau}t)\bigl)\bigl(t^2E_i(\si,t^{-1}r)
-A_i(\si,r)-t^2T_i(\si,r)\bigr)\bigr].
\end{split}
\label{co6eq22}
\e
So applying \cite[Prop.~6.3]{Joyc4} again as in \eq{co6eq19} gives
\e
\begin{split}
(\Xi_i^t&)^*(\ve^t)(\si,r)=O\bigl(r^{-2}\ms{\xi_i^t-\eta_i}\bigr)
+O\bigl(\ms{\na(\xi_i^t-\eta_i)}\bigr)\\
&-\d^*\bigl(\psi^m\d\bigl[\bigl(1-F(t^{-\tau}t)\bigl)
\bigl(t^2E_i(\si,t^{-1}r)-A_i(\si,r)-t^2T_i(\si,r)\bigr)\bigr]\bigr)
\end{split}
\label{co6eq23}
\e
on $\Si_i\t(t^\tau,2t^\tau)$, provided $r^{-1}\md{\xi_i^t-\eta_i}$
and $\md{\na(\xi_i^t-\eta_i)}$ are small. Here $\md{\,.\,}$ and $\d^*$
are computed using the metric $(\Up_i\circ\phi_i)^*(g)$, and we have
used $\d^*(\psi^m\al)=0$ to eliminate the term~$-\d^*\bigl(\psi^mt^2
(\Up_i\circ\phi_i)^*(\al)\bigr)$.

It is important that there are {\it no linear terms} in
$t^2(\Up_i\circ\phi_i)^*(\al)$ or $t^2\pi_i^*(\ga_i)$ in
\eq{co6eq23}. If we had applied the cruder method used in
\cite[Prop.~6.4]{Joyc5} there would have been such linear
terms, which would have contributed $O(t^{2-2\tau})$ to
$\md{(\Xi_i^t)^*(\ve^t)}$ on $\Si_i\t(t^\tau,2t^\tau)$ and
$O(t^{2+\tau(m/2-1)})$ to $\lnm{\ve^t}{2m/(m+2)}$, which is
too large for part (i) of Theorem \ref{co5thm1} to hold.

Using \eq{co6eq6} to write terms in \eq{co6eq23} in terms of
$\ga_i,A_i,E_i,T_i$ and estimating these using \eq{co3eq11},
\eq{co6eq2} and \eq{co6eq4}, we find that $r^{-1}\md{\xi_i^t
-\eta_i}$ and $\md{\na(\xi_i^t-\eta_i)}$ are $O(t^{2-2\tau})$
on $\Si_i\t(t^\tau,2t^\tau)$, as the $t^2\pi_i^*(\ga_i)$ term
in \eq{co6eq6} dominates, so that \eq{co6eq23} holds for small
$t$ and on $\Si_i\t(t^\tau,2t^\tau)$ we have
\e
(\Xi_i^t)^*(\ve^t)=O(t^{4-4\tau})+O(t^{\tau(\mu_i-2)})+
O(t^{(1-\tau)(2-\la)})+O(t^{2+\tau(\nu_i-2)}).
\label{co6eq24}
\e
By a similar argument for derivatives we find that on
$\Si_i\t(t^\tau,2t^\tau)$ we have
\e
(\Xi_i^t)^*(\d\ve^t)=O(t^{4-5\tau})+O(t^{\tau(\mu_i-3)})+
O(t^{(1-\tau)(2-\la)-\tau})+O(t^{2+\tau(\nu_i-3)}).
\label{co6eq25}
\e
Making $\de>0$ smaller if necessary, equations
\eq{co6eq10}--\eq{co6eq13} now follow from \eq{co6eq16},
\eq{co6eq17}, \eq{co6eq20}, \eq{co6eq21}, \eq{co6eq24}
and \eq{co6eq25}, for some $C>0$ independent of~$t$.
\end{proof}

Now we can estimate norms of $\ve^t$ and $\d\ve^t$, as in part
(i) of Theorem~\ref{co5thm1}.

\begin{prop} There exists $C'>0$ such that for all\/ $t\in(0,\de)$
we have
\ea
&\lnm{\ve^t}{\sst\frac{2m}{m+2}}\!\le\!C't^{\tau(1\!+\!m/2)}
\Bigl(t^{4\!-\!4\tau}\!+\!t^{(1\!-\!\tau)(2\!-\!\la)}\!+\!
\sum_{i=1}^n\bigl(t^{\tau(\mu_i\!-\!2)}\!+\!
t^{2\!+\!\tau(\nu_i\!-\!2)}\bigr)\Bigr),
\label{co6eq26}\\
&\cnm{\ve^t}{0}\!\le\!
C'\Bigl(t^{4-4\tau}\!+\!t^{(1-\tau)(2-\la)}\!+\!
\sum_{i=1}^n\bigl(t^{\tau(\mu_i-2)}\!+\!t^{2+\tau(\nu_i-2)}\bigr)\Bigr),
\;\>\text{and}
\label{co6eq27}\\
&\lnm{\d\ve^t}{2m}\!\le\!C't^{-\tau/2}\Bigl(t^{4-4\tau}\!+\!
t^{(1-\tau)(2-\la)}\!+\!
\sum_{i=1}^n\bigl(t^{\tau(\mu_i-2)}\!+\!t^{2+\tau(\nu_i-2)}\bigr)\Bigr),
\label{co6eq28}
\ea
computing norms with respect to the metric $h^t$ on~$N^t$.
\label{co6prop2}
\end{prop}

\begin{proof} The proof is similar to that of \cite[Prop.~6.5]{Joyc5}.
Using \eq{co6eq10}, \eq{co6eq12} and \eq{co6eq13} we find that $\lnm{
\ve^t}{2m/(m+2)}$ has contributions $O(t^{2+m/2})$ from $P_i^t$,
$O(t^{\tau(2+m/2)})$ from $\Xi_i^t\bigl(\Si_i\t(tT,t^\tau]\bigr)$, and
\begin{equation*}
O(t^{4+\tau(m-6)/2})\!+\!O(t^{\tau(\mu_i-1+m/2)})\!+\!
O(t^{2-\la+\tau(\la-1+m/2)})\!+\!O(t^{2+\tau(\nu_i-1+m/2)})
\end{equation*}
from $\Xi_i^t\bigl(\Si_i\t(t^\tau,2t^\tau)\bigr)$, and $O(t^{4+
\tau(m-6)/2})$ from $\Xi_i^t\bigl(\Si_i\t[2t^\tau,R')\bigr)$ as
$m<6$, and $O(t^4)$ from $K^t$, since $\vol(K^t)=O(1)$.

Now $O(t^{\tau(\mu_i-1+m/2)})$ dominates $O(t^{\tau(2+m/2)})$
and $O(t^{2+m/2})$ as $\mu_i<3$ and $0<\tau<1$, and 
$O(t^{4+\tau(m-6)/2})$ dominates $O(t^4)$ as $m<6$ and $\tau>0$.
The remaining terms give \eq{co6eq26} for some $C'>0$. Equations
\eq{co6eq27}--\eq{co6eq28} follow from \eq{co6eq10}--\eq{co6eq13}
in the same way.                                    
\end{proof}

We also need to estimate $\lnm{\pi_\sWt(\ve^t)}{1}$. Following
\cite[\S 7.2]{Joyc5}, which assumes condition (iv) of Definition
\ref{co6def1}, with no significant changes we prove:

\begin{prop} There exists $C''>0$ such that for all\/
$t\in(0,\de)$ we have $\lnm{\pi_\sWt(\ve^t)}{1}
\le C''t^{(m+1)\tau}$, computing norms using $h^t$ on~$N^t$.
\label{co6prop3}
\end{prop}

For part (i) of Theorem \ref{co5thm1} to hold, we want
$\lnm{\ve^t}{2m/(m+2)}\le A_2t^{\ka+m/2}$, $\cnm{\ve^t}{0}\le
A_2t^{\ka-1}$, $\lnm{\d\ve^t}{2m}\le A_2t^{\ka-3/2}$ and
$\lnm{\pi_\sWt(\ve^t)}{1}\le A_2t^{\ka+m-1}$ for
some $\ka>1$. As $t<1$ we see from \eq{co6eq26} that
$\lnm{\ve^t}{2m/(m+2)}\le A_2t^{\ka+m/2}$ holds with
$A_2\ge 2(n\!+\!1)C'$ provided
\ea
\tau(1\!+\!m/2)\!+\!4\!-\!4\tau &\!\ge\!\ka\!+\!m/2,&\;
\tau(1\!+\!m/2)\!+\!(1\!-\!\tau)(2\!-\!\la)&\!\ge\!\ka\!+\!m/2,
\label{co6eq29}\\             
\tau(1\!+\!m/2)\!+\!\tau(\mu_i\!-\!2)&\!\ge\!\ka\!+\!m/2,&\;
\tau(1\!+\!m/2)\!+\!2\!+\!\tau(\nu_i\!-\!2)&\ge\ka\!+\!m/2,
\label{co6eq30}
\ea
for all $i=1,\ldots,n$. As $\tau\le 1$ we find from
\eq{co6eq26}--\eq{co6eq28} that \eq{co6eq29}--\eq{co6eq30} also
imply $\cnm{\ve^t}{0}\le A_2t^{\ka-1}$ and $\lnm{\d\ve^t}{2m}\le
A_2t^{\ka-3/2}$. Similarly, Proposition \ref{co6prop3} implies
that $\lnm{\pi_\sWt(\ve^t)}{1}\le A_2t^{\ka+m-1}$
holds with $A_2\ge C''$ provided
\e
(m+1)\tau\ge\ka+m-1.
\label{co6eq31}
\e

Now \eq{co6eq29} is equivalent to $\ka\!\le\!1\!+\!(1\!-\!\tau)
(6\!-\!m)/2$ and $\ka\le 1+(1-\tau)(1+m/2-\la)$, which admit a
solution $\ka>1$ as $\tau<1$, $m<6$ and $\la<-1-m/2$. Also, the
conditions on $\tau$ in \eq{co6eq5} ensure that \eq{co6eq30}
for $i=1,\ldots,n$ and \eq{co6eq31} admit a solution $\ka>1$,
and this is the reason for \eq{co6eq5}. Taking this $\ka$ and
$A_2=\max\bigl(2(n+1)C',C''\bigr)$, we have proved the
analogue of~\cite[Th.~6.6]{Joyc5}:

\begin{thm} Making $\de>0$ smaller if necessary, there
exist\/ $A_2>0$ and\/ $\ka>1$ such that the functions
$\ve^t=\psi^m\sin\th^t$ on $N^t$ satisfy
$\lnm{\ve^t}{2m/(m+2)}\le A_2t^{\ka+m/2}$,
$\cnm{\ve^t}{0}\le A_2t^{\ka-1}$,
$\lnm{\d\ve^t}{2m}\le A_2t^{\ka-3/2}$ and\/
$\lnm{\pi_\sWt(\ve^t)}{1}\le A_2t^{\ka+m-1}$
for all\/ $t\in(0,\de)$, as in part\/ {\rm(i)}
of Theorem~\ref{co5thm1}.
\label{co6thm1}
\end{thm}

\subsection{Lagrangian neighbourhoods and bounds on $R(h^t),\de(h^t)$}
\label{co63}

It turns out that \cite[\S 6.3--\S 6.5]{Joyc5} need only minor changes
to apply to the Lagrangian $m$-folds $N^t$ of Definition \ref{co6def2}.
Here is the analogue of~\cite[Def.~6.7]{Joyc5}.

\begin{dfn} Define an open neighbourhood $U_\sNt\subset T^*N^t$
of the zero section $N^t$ in $T^*N^t$ and a smooth map $\Phi_\sNt:
U_\sNt\ra M$ as follows. As $N^t$ is the disjoint union of $K^t$
and $P_i^t,Q_i^t$ for $i=1,\ldots,n$ we shall define $U_\sNt$
and $\Phi_\sNt$ separately over $K^t,P_i^t$ and~$Q_i^t$.

Let $\Pi^t:K^t\ra K$ be the natural projection, recalling
from \eq{co6eq7} that $K^t$ is isomorphic to the graph
$\Ga(t^2\al)$ of $t^2\al$ over $K$. Then $\Pi^t$ is a
diffeomorphism, so $\d\Pi^t:T^*K^t\ra T^*K$ is a vector
bundle isomorphism. Write $\pi$ for both projections
$\pi:T^*N^t\ra N^t$ and $\pi:T^*K\ra K$.

Now let $y\in T^*N^t\cap\pi^*(K^t)$, and define
\e
x=\pi(y)\in K^t,\quad
y'=\d\Pi^t(y)\in T^*K\quad\text{and}\quad
x'=\pi(y')=\Pi^t(x)\in K.
\label{co6eq32}
\e
Define $U_\sNt\cap\pi^*(K^t)$ and $\Phi_\sNt\vert_{U_\sNt
\cap\pi^*(K^t)}$ by
\e
\begin{gathered}
\text{$y\in U_\sNt\cap\pi^*(K^t)$ if and only if
$y'+t^2\al(x')\in U_\sXp$,}\\
\text{and then $\Phi_\sNt(y)=\Phi_\sXp\bigl(y'+t^2\al(x')\bigr)$.}
\end{gathered}
\label{co6eq33}
\e

For $i=1,\ldots,n$, define $U_\sNt\cap\pi^*(P_i^t)$ and
$\Phi_\sNt\vert_{U_\sNt\cap\pi^*(P_i^t)}$ by
\e
\begin{gathered}
U_\sNt\cap\pi^*(P_i^t)=\d(\Up_i\circ t)
\bigl(\{\ga\in T^*K_i:t^{-2}\ga\in U_\sLi\}\bigr)\\
\text{and}\quad \Phi_\sNt\circ\d(\Up_i\circ t)(\ga)=
\Up_i\circ t\circ\Phi_\sLi(t^{-2}\ga).
\end{gathered}
\label{co6eq34}
\e
Here the diffeomorphism $\Up_i\circ t:K_i\ra P_i^t$ induces
$\d(\Up_i\circ t):T^*K_i\ra T^*P_i^t$, and $\ga\mapsto t^{-2}\ga$
is multiplication by $t^{-2}$ in the vector space fibres
of~$T^*K_i\ra K_i$.

Let $F$ be as in Definition \ref{co6def2}. Now $\Xi_i^t:\Si_i\t
(tT,R')\ra Q_i^t$ is a diffeomorphism, and induces an isomorphism
$\d\Xi_i^t:T^*\bigl(\Si_i\t(tT,R')\bigr)\ra T^*Q_i^t$. As in
\eq{co3eq12}--\eq{co3eq13} and \eq{co4eq7}--\eq{co4eq8}, define
$U_\sNt\cap\pi^*(Q_i^t)$ and $\Phi_\sNt\vert_{U_\sNt\cap\pi^*(Q_i^t)}$~by
\begin{gather}
\begin{split}
(\d\Xi_i^t)^*(U_\sNt)=\bigl\{&(\si,r,\varsigma,u)
\in T^*\bigl(\Si_i\t(tT,R')\bigr):\\
&\bmd{(\varsigma,u)+t^2F(t^{-\tau}r)
(\Up_i\circ\phi_i)^*(\al)(\si,r)}<\ze r\bigr\}
\quad\text{and}
\end{split}
\label{co6eq35}\\
\Phi_\sNt\circ\d\Xi_i^t(\si,r,\varsigma,u)\equiv\Up_i\circ\Phi_\sCi
\bigl(\si,r,\varsigma+\xi_i^{1,t}(\si,r),u+\xi_i^{2,t}(\si,r)\bigr)
\label{co6eq36}
\end{gather}
for all $(\si,r,\varsigma,u)\in(\d\Xi_i^t)^*(U_\sNt)$, computing
$\na,\md{\,.\,}$ using~$\iota_i^*(g')$.

Careful consideration shows that $U_\sNt$ is well-defined, and
$\Phi_\sNt$ is well-defined in \eq{co6eq33}, \eq{co6eq34} and
\eq{co6eq36} for small $t$, so making $\de>0$ smaller if necessary
$\Phi_\sNt$ is well-defined for $t\in(0,\de)$. Clearly $\Phi_\sNt$
is smooth on each of $U_\sNt\cap\pi^*(K^t)$, $U_\sNt\cap\pi^*(P_i^t)$
and $U_\sNt\cap\pi^*(Q_i^t)$, but we must still show that $\Phi_\sNt$
is smooth over the joins between them.

The condition $\md{(\varsigma,u)+t^2F(t^{-\tau}r)(\Up_i\circ
\phi_i)^*(\al)(\si,r)}<\ze r$ in \eq{co6eq35} is chosen to ensure
that the definitions of $U_\sNt$ over $K^t,P_i^t$ and $Q_i^t$ join
smoothly together over $\pd K^t$ and $\pd P_i^t$. Therefore $U_\sNt$
is an {\it open tubular neighbourhood\/} of $N^t$ in $T^*N^t$.
Following \cite[Def.~6.4]{Joyc5} we find that $\Phi_\sNt$ is
{\it smooth\/} on $U_\sNt$, that is, the definitions over
$P_i^t,Q_i^t$ and $K^t$ join together smoothly over $\pd P_i^t$
and $\pd K^t$.

We shall show that $\Phi_\sNt^*(\om)=\hat\om$. On $U_\sNt\cap\pi^*(K^t)$
this holds as $\Phi_\sXp^*(\om)=\hat\om$ and $\al$ is closed. On
$U_\sNt\cap\pi^*(P_i^t)$ it follows from $\Up_i^*(\om)=\om'$,
$\Phi_\sLi^*(\om')=\hat\om$, and the fact that the powers of $t$
in \eq{co6eq34} cancel out in their effect on $\Phi_\sNt^*(\om)$.
On $U_\sNt\cap\pi^*(Q_i^t)$ it holds as $\Up_i^*(\om)=\om'$,
$\Phi_\sCi^*(\om')=\hat\om$ and $\xi_i^t$ is closed.

Define an $m$-form $\be^t$ on $U_\sNt$ by $\be^t=\Phi_\sNt^*(\Im\Om)$,
as in Definition~\ref{co5def2}.
\label{co6def4}
\end{dfn}

Using this Lagrangian neighbourhood map $\Phi_\sNt$ in part (iii),
we find that that parts (ii)--(v) of Theorem \ref{co5thm1} hold
for $N^t$ when $t\in(0,\de)$, as in~\cite[Th.~6.8]{Joyc5}.

\begin{thm} Making $\de>0$ smaller if necessary, there exist\/
$A_1,A_3,\ldots,A_6>0$ such that for all\/ $t\in(0,\de)$,
as in parts {\rm(ii)--(v)} of Theorem \ref{co5thm1} we have
\begin{itemize}
\item[{\rm(ii)}] $\psi\ge A_3$ on $N^t$.
\item[{\rm(iii)}] The subset\/ $\B_{A_1t}\subset T^*N^t$
of Definition \ref{co5def2} lies in $U_\sNt$, and
$\cnm{\hat\na^k\be^t}{0}\le A_4t^{-k}$ on $\B_{A_1t}$ for
$k=0,1,2$ and\/~$3$.
\item[{\rm(iv)}] The injectivity radius $\de(h^t)$
satisfies~$\de(h^t)\ge A_5t$.
\item[{\rm(v)}] The Riemann curvature $R(h^t)$
satisfies~$\cnm{R(h^t)}{0}\le A_6t^{-2}$.
\end{itemize}
Here part\/ {\rm(iii)} uses the notation of Definition \ref{co5def2},
and parts {\rm(iv)} and\/ {\rm(v)} refer to the compact Riemannian
manifold\/~$(N^t,h^t)$.
\label{co6thm2}
\end{thm}

\begin{proof} This was proved in \cite[Th.~6.8]{Joyc5}, but
with $N^t,\Phi_\sNt$ defined more simply. The changes
for the new $N^t,\Phi_\sNt$ are very minor. The main
difference is that in \cite[\S 6.3]{Joyc5} $N^t,\Phi_\sNt$
are independent of $t$ over $K$, but here $N^t,\Phi_\sNt$ do
depend on $t$ over $K^t$. Identifying $K^t$ with $K$ using
$\Pi^t$ in Definition \ref{co6def4} we have $\na^j(h^t-g\vert_K)
=O(t^2)$ for $j\ge 0$. Therefore the contributions to
$\cnm{\hat\na^k\be^t}{0},\de(h^t)$ and $\cnm{R(h^t)}{0}$ from
$K^t$ are all $O(1)$ for small $t$, as in~\cite[\S 6.3]{Joyc5}.
\end{proof}

\subsection{Sobolev embedding inequalities on $N^t$}
\label{co64}

We now prove that parts (vi) and (vii) of Theorem \ref{co5thm1}
hold for $N^t,W^t$. Here is the analogue of \cite[Th.~7.8]{Joyc5},
which gives part (vi) of Theorem \ref{co5thm1} for~$N^t,W^t$.

\begin{thm} Making $\de>0$ smaller if necessary, there exists
$A_7>0$ such that for all\/ $t\in(0,\de)$, if\/ $v\in L^2_1(N^t)$
with\/ $\int_{N^t}vw\,\d V^t=0$ for all\/ $w\in W^t$ then $v\in
L^{2m/(m-2)}(N^t)$ and\/~$\lnm{v}{2m/(m-2)}\le A_7\lnm{\d v}{2}$.
\label{co6thm3}
\end{thm}

\begin{proof} This was proved in \cite[Th.~6.12]{Joyc5} for
$X'$ connected and \cite[Th.~7.8]{Joyc5} for general $X'$,
but with $N^t,\Phi_\sNt$ defined more simply. Here is how
to modify the proofs. First consider the case of
\cite[\S 6.4]{Joyc5}, when $X'$ is connected.

In \cite[\S 6.4]{Joyc5} we define a partition of unity
function $F^t$ on $N^t$. The new definition is as follows.
Choose $a,b\in\R$ with $0<a<b<\tau$. Making $\de>0$ smaller
if necessary we have $2t^\tau<t^b<t^a<\min(1,R')$ for all
$t\in(0,\de)$. Let $G:(0,\iy)\ra[0,1]$ be a smooth,
decreasing function with $G(s)=1$ for $s\in(0,a]$ and $G(s)=0$
for $s\in[b,\iy)$. For $t\in(0,\de)$, define a function
$F^t:N^t\ra[0,1]$ by
\e
F^t(x)=\begin{cases} 1, & x\in K^t, \\
G\bigl((\log r)/(\log t)\bigr), & x=\Xi_i^t(\si,r)\in Q_i^t,
\quad i=1,\ldots,n, \\
0, & x\in P_i^t,\quad i=1,\ldots,n.
\end{cases}
\label{co6eq37}
\e

The main point we have to deal with is that in
\cite[\S 6.4]{Joyc5} the function $F^t$ was supported on
$N^t\cap X'$, and so if $v\in C^1(N^t)$ then we could treat
$F^tv$ as a compactly-supported function on $X'$. But for
the new $N^t$ this no longer holds. Instead, on the support
of $F^t$ we have $N^t=\Phi_\sXp\bigl(\Ga(t^2\al)\bigr)$,
where $\Ga(t^2\al)$ is the graph of $t^2\al$ in~$T^*X'$.

Identifying $\Ga(t^2\al)$ with $X'$ using $\pi:T^*X'\ra X'$
and $N^t$ with $\Ga(t^2\al)$ using $\Phi_\sXp$ defines an
identification between $N^t$ and $X'$ on the support of
$F^t$. Thus, if $v\in C^1(N^t)$ we can regard $F^tv$ as a
compactly-supported function on $X'$. Then \cite[Prop.~6.11]{Joyc5}
gives
\e
\lnm{F^tv}{2m/(m-2)}\le D_2\bigl(\blnm{\,\d(F^tv)}{2}+
\bmd{\ts\int_{X'}F^tv\,\d V_g\,}\,\bigr),
\label{co6eq38}
\e
computing norms using $g$ on $X'$, where $D_2>0$ is
independent of~$v,t$.

For the proof of \cite[Th.~6.12]{Joyc5} to work for the new
$N^t$, we need \eq{co6eq38} to hold computed using $h^t$ on
$N^t$. Since $h^t$ is the restriction of $h^t$ to
$\Ga(t^2\al)$ and $\al=O(\rho^{-1})$, $\na\al=O(\rho^{-2})$
we find that
\e
h^t=g\vert_{X'}+O(t^2\rho^{-2})=g\vert_{X'}+O(t^{2-2\tau})
\label{co6eq39}
\e
on the support of $F^t$, identifying $N^t$ and $X'$ there as
above. Thus $h^t$ and $g\vert_{X'}$ are uniformly equivalent
on the support of $F^t$ for small $t$, so increasing $D_2$
we see that \eq{co6eq38} holds with $\lnm{\,.\,}{2m/(m-2)}$,
$\lnm{\,.\,}{2}$  computed using~$h^t$.

From \eq{co6eq39} we have $\d V_g=\bigl(1+O(t^{2-2\tau})\bigr)
\d V^t$, so that
\e
\begin{split}
\ts\bmd{\int_{X'}F^tv\,\d V_g\,}&\le\ts\bmd{\int_{N^t}F^tv\,\d V^t\,}
+Ct^{2-2\tau}\int_{N^t}\md{F^tv}\,\d V^t\\
&\le\ts\bmd{\int_{N^t}F^tv\,\d V^t\,}+C't^{2-2\tau}
\lnm{F^tv}{2m/(m-2)},
\end{split}
\label{co6eq40}
\e
for small $t$ and $C,C'>0$ independent of $t$. Here norms
are computed using $h^t$, and in the second line we use
$\lnm{F^tv}{1}\le\vol(N^t)^{(m+2)/2m}\lnm{F^t}{2m/(m-2)}$,
and that $\vol(N^t)$ is bounded independently of~$t$.

Substituting \eq{co6eq40} into \eq{co6eq38} gives
\begin{equation*}
(1-C'D_2t^{2-2\tau})\lnm{F^tv}{2m/(m-2)}\le D_2\bigl(
\blnm{\,\d(F^tv)}{2}+\bmd{\ts\int_{N^t}F^tv\,\d V^t\,}\,\bigr),
\end{equation*}
computing norms using $h^t$. For small $t$ we have
$C'D_2t^{2-2\tau}\le\ha$, and so
\begin{equation*}
\lnm{F^tv}{2m/(m-2)}\le 2D_2\bigl(
\blnm{\,\d(F^tv)}{2}+\bmd{\ts\int_{N^t}F^tv\,\d V^t\,}\,\bigr).
\end{equation*}
This is the first line of \cite[eq.~(79)]{Joyc5}, with $2D_2$ in
place of $D_2$. The rest of the proof of \cite[Th.~6.12]{Joyc5}
then follows. This proves the theorem when $X'$ is connected.
For general $X'$ we follow the proof of \cite[Th.~7.8]{Joyc5},
but with $F^t$ defined as in \eq{co6eq37}. The modifications
above extend easily, and the theorem follows.
\end{proof}

The proof of \cite[Th.~7.9]{Joyc5} also holds for the new
$N^t,W^t$ with no significant changes. This proves part (vii)
of Theorem \ref{co5thm1} for~$N^t,W^t$.

\begin{thm} Making $\de>0$ smaller if necessary, for all\/
$t\in(0,\de)$ and\/ $w\in W^t$ we have $\lnm{\d^*\d w}{2m/(m+2)}
\le\ha A_7^{-1}\lnm{\d w}{2}$, where $A_7>0$ is as in Theorem
\ref{co6thm3}. Also there exists $A_8>0$ such that for all\/
$t\in(0,\de)$ and\/ $w\in W^t$ with\/ $\int_{N^t}w\,\d V^t=0$
we have~$\cnm{w}{0}\le A_8t^{1-m/2}\lnm{\d w}{2}$.
\label{co6thm4}
\end{thm}

\subsection{The main result when $Y(L_i)\ne 0$ and $\la_i=0$}
\label{co65}

Here is the first of the main results of this paper, the
analogue of Theorem \ref{co5thm3} but allowing $\la_i=0$
and~$Y(L_i)\ne 0$.

\begin{thm} Let\/ $(M,J,\om,\Om)$ be an almost Calabi--Yau $m$-fold
for $2\!<\!m\!<\nobreak\!6$, and\/ $X$ a compact SL\/ $m$-fold
in $M$ with conical singularities at\/ $x_1,\ldots,x_n$ and cones
$C_1,\ldots,C_n$. Define $\psi:M\ra(0,\iy)$ as in \eq{co2eq3}. Let\/
$L_1,\ldots,L_n$ be Asymptotically Conical SL\/ $m$-folds in $\C^m$
with cones $C_1,\ldots,C_n$ and rates $\la_1,\ldots,\la_n$. Write
$X'=X\sm\{x_1,\ldots,x_n\}$ and\/ $\Si_i=C_i\cap{\cal S}^{2m-1}$.
Suppose that\/ $\la_i\le 0$ for $i=1,\ldots,n$, and that there
exists $\varrho\in H^1(X',\R)$ such that\/ $\bigl(Y(L_1),\ldots,
Y(L_n)\bigr)$ is the image of\/ $\varrho$ under the map $H^1(X',\R)
\ra\bigoplus_{i=1}^nH^1(\Si_i,\R)$ in~\eq{co3eq3}.

Set\/ $q=b^0(X')$, and let\/ $X_1',\ldots,X_q'$ be the connected
components of\/ $X'$. For $i=1,\ldots,n$ let\/ $l_i=b^0(\Si_i)$,
and let\/ $\Si_i^1,\ldots,\Si_i^{\smash{l_i}}$ be the connected
components of\/ $\Si_i$. Define $k(i,j)=1,\ldots,q$ by $\Up_i
\circ\vp_i\bigl(\Si_i^j\t(0,R')\bigr)\subset X'_{\smash{k(i,j)}}$
for $i=1,\ldots,n$ and\/ $j=1,\ldots,l_i$. Suppose that
\e
\sum_{\substack{1\le i\le n, \; 1\le j\le l_i: \\
k(i,j)=k}}\psi(x_i)^mZ(L_i)\cdot[\Si_i^j\,]=0
\quad\text{for all\/ $k=1,\ldots,q$.}
\label{co6eq41}
\e

Suppose also that the compact\/ $m$-manifold\/ $N$ obtained
by gluing $L_i$ into $X'$ at\/ $x_i$ for $i=1,\ldots,n$ is
connected. A sufficient condition for this to hold is that\/
$X$ and\/ $L_i$ for $i=1,\ldots,n$ are connected.

Then there exists $\ep>0$ and a smooth family $\smash{\bigl\{
\ti N^t:t\in(0,\ep]\bigr\}}$ of compact, nonsingular SL\/
$m$-folds in $(M,J,\om,\Om)$ diffeomorphic to $N$, such
that\/ $\smash{\ti N^t}$ is constructed by gluing $tL_i$
into $X$ at\/ $x_i$ for $i=1,\ldots,n$. In the sense of
currents in Geometric Measure Theory, $\smash{\ti N^t}\ra
X$ as~$t\ra 0$.
\label{co6thm5}
\end{thm}

\begin{proof} The hypotheses of the theorem imply that
Definition \ref{co6def1} holds. Let $\de>0$ and $N^t,W^t$ for
$t\in(0,\de)$ be as in Definitions \ref{co6def2} and \ref{co6def3},
and make $\de>0$ smaller if necessary so that Theorems \ref{co6thm1},
\ref{co6thm2}, \ref{co6thm3} and \ref{co6thm4} apply. Theorem
\ref{co6thm1} gives constants $\ka>1$ and $A_2>0$ such that
part (i) of Theorem \ref{co5thm1} holds for $N^t$ for all
$t\in(0,\de)$, replacing $N,W,\th$ by $N^t,W^t,\th^t$ respectively.

Let the Lagrangian neighbourhood $\Phi_\sNt:U_\sNt\ra M$ and the
$m$-form $\be^t$ on $U_\sNt$ be as in Definition \ref{co6def4}.
Then Theorem \ref{co6thm2} gives constants $A_1,A_3,\ldots,A_6>0$
such that parts (ii)--(v) of Theorem \ref{co5thm1} hold for $N^t$
for all $t\in(0,\de)$, replacing $N,\be,h$ by $N^t,\be^t,h^t$
respectively. Theorems \ref{co6thm3} and \ref{co6thm4} give
$A_7,A_8>0$ such that parts (vi) and (vii) of Theorem \ref{co5thm1}
hold for $N^t$ for all $t\in(0,\de)$, replacing $N,W$ by $N^t,W^t$
respectively.

We have not yet shown that the inequality $\cos\th^t\ge\ha$ in
Definition \ref{co5def2} holds. From parts (i) and (ii) of Theorem
\ref{co5thm1} we see that $\md{\sin\th^t}\le A_2A_3^{-m}t^{\ka-1}$
on $N^t$. Thus for small $t\in(0,\de)$ we have $\md{\sin\th^t}\le
\frac{\sqrt{3}}{2}$ as $\ka>1$, so that $\md{\cos\th^t}\ge\ha$.
As $N^t$ is approximately special Lagrangian on $K^t$ we have
${\rm e}^{i\th^t}\approx 1$ on $K^t$, so $\cos\th^t\ge\ha$ on
$N^t$ as $\cos\th^t$ is continuous and $N^t$ connected.

Let $\ep,K>0$ be as given in Theorem \ref{co5thm1} depending
on $\ka,A_1,\ldots,A_8$ and $m$, and make $\ep>0$ smaller if
necessary to ensure that $\ep<\de$ and $\cos\th^t\ge\ha$ on
$N^t$ for $t\le\ep$. Then Theorem \ref{co5thm1} shows that for
all $t\in(0,\ep]$ we can deform $N^t$ to a nearby compact,
nonsingular SL $m$-fold $\smash{\ti N^t}$, given by
$\smash{\ti N^t}=(\Phi_\sNt)_*\bigl(\Ga(\d f^t)\bigr)$ for
some $f^t\in C^\iy(N^t)$ with~$\cnm{\d f^t}{0}\le Kt^\ka<A_1t$.

Since $N^t$ and $\Phi_\sNt$ depend smoothly on $t$, we see that
$f^t$ is the locally unique solution of a nonlinear elliptic
p.d.e.\ on $N^t$ depending smoothly on $t$. It quickly follows
from general theory that $f^t$ depends smoothly on $t$, and so
$\smash{\ti N^t}$ does. We show that $\smash{\ti N^t}\ra X$ as
currents as $t\ra 0$ as in~\cite[\S 6.5]{Joyc5}.
\end{proof}

When $m=6$, the proof of Theorem \ref{co6thm1} shows that the
bounds in part (i) of Theorem \ref{co5thm1} hold with $\ka=1$
rather than $\ka>1$. Thus the proof {\it only just\/} fails
when $m=6$, and with some improvements to the material of
\S\ref{co61}--\S\ref{co62} Theorem \ref{co6thm5} can probably
be proved for $m=6$. However, the author does not know how to
extend Theorem \ref{co6thm5} to the case~$m>6$.

If $X'$ is connected, so that $q=1$, then $k(i,j)\equiv 1$
and \eq{co6eq41} becomes
\begin{equation*}
\sum_{i=1}^n\psi(x_i)^mZ(L_i)\cdot\sum_{j=1}^{l_i}[\Si_i^j\,]=0.
\end{equation*}
But $\sum_{j=1}^{l_i}[\Si_i^j\,]=[\Si_i]$, and $Z(L_i)\cdot[\Si_i]=0$
as $Z(L_i)$ is the image of a class in $H^{m-1}(L_i,\R)$ by Definition
\ref{co4def2}, and $\Si_i=\pd L_i$, so $[\Si_i]$ maps to zero in
$H_{m-1}(L_i,\R)$. Therefore \eq{co6eq41} holds automatically when
$X'$ is connected. Also, $N$ is automatically connected. Thus
Theorem \ref{co6thm5} simplifies in this case to give an analogue
of Theorem~\ref{co5thm2}:

\begin{thm} Let\/ $(M,J,\om,\Om)$ be an almost Calabi--Yau $m$-fold
for $2\!<\!m\!<\nobreak\!6$, and\/ $X$ a compact SL\/ $m$-fold
in $M$ with conical singularities at\/ $x_1,\ldots,x_n$ and cones
$C_1,\ldots,C_n$. Let\/ $L_1,\ldots,L_n$ be Asymptotically
Conical SL\/ $m$-folds in $\C^m$ with cones $C_1,\ldots,C_n$
and rates $\la_1,\ldots,\la_n$. Suppose that\/ $\la_i\le 0$
for\/ $i=1,\ldots,n$, that\/ $X'=X\sm\{x_1,\ldots,x_n\}$ is
connected, and that there exists $\varrho\in H^1(X',\R)$
such that\/ $\bigl(Y(L_1),\ldots,Y(L_n)\bigr)$ is the image
of\/ $\varrho$ under the map $H^1(X',\R)\ra\bigoplus_{i=1}^n
H^1(\Si_i,\R)$ in \eq{co3eq3}, where~$\Si_i=C_i\cap{\cal S}^{2m-1}$.

Then there exists $\ep>0$ and a smooth family $\bigl\{
\smash{\ti N^t}:t\in(0,\ep]\bigr\}$ of compact, nonsingular SL\/
$m$-folds in $(M,J,\om,\Om)$, such that\/ $\smash{\ti N^t}$ is
constructed by gluing $tL_i$ into $X$ at\/ $x_i$ for
$i=1,\ldots,n$. In the sense of currents, $\smash{\ti N^t}\ra
X$ as~$t\ra 0$.
\label{co6thm6}
\end{thm}

\section{Desingularizing in families when $Y(L_i)=0$}
\label{co7}

Next we consider a different generalization of Theorems
\ref{co5thm2} and \ref{co5thm3}, to {\it families} of
almost Calabi--Yau $m$-folds $(M,J^s,\om^s,\Om^s)$. Let
$\F\subset\R^d$ be open with $0\in\F$, and $\bigl\{(M,J^s,
\om^s,\Om^s):s\in\F\bigr\}$ be a smooth family of almost
Calabi--Yau $m$-folds with~$(M,J^0,\om^0,\Om^0)=(M,J,\om,\Om)$.

Suppose $X$ is an SL $m$-fold in $(M,J,\om,\Om)$ with
conical singularities at $x_1,\ldots,x_n$ and cones $C_i$,
and that $L_1,\ldots,L_n$ are AC SL $m$-folds in $\C^m$ with
cones $C_i$ and rates $\la_i$. In the rest of the paper we
shall construct special Lagrangian desingularizations
$\smash{\ti N^{s,t}}$ of $X$ not just in $(M,J,\om,\Om)$, as
in \cite{Joyc5} and \S\ref{co6}, but in $(M,J^s,\om^s,\Om^s)$
for small~$s\in\F$.

In this section, as in Theorems \ref{co5thm2} and \ref{co5thm3}
but not as in \S\ref{co6}, we shall assume that $\la_i<0$, so
that $Y(L_i)=0$ by Proposition \ref{co4prop}. Then \S\ref{co8}
combines the new material of \S\ref{co6} and this section to
study desingularization in families when~$Y(L_i)\ne 0$.

The advantage of desingularizing in families $(M,J^s,\om^s,\Om^s)$,
rather than a single almost Calabi--Yau $m$-fold $(M,J,\om,\Om)$, is
that by varying the cohomology classes $[\om^s]$ and $[\Im\Om^s]$ we
can overcome obstructions to the existence of SL desingularizations
$\smash{\ti N^t}$ in $(M,J,\om,\Om)$. In this section, by varying
$[\Im\Om^s]$ we show how to relax equation \eq{co5eq2} in
Theorem~\ref{co5thm3}.

\subsection{Setting up the problem}
\label{co71}

The next three definitions are analogues of Definitions
\ref{co6def1}--\ref{co6def3}. We consider the following
situation, modifying \cite[Def.~6.1 \& Def.~7.1]{Joyc5}
to the families case.

\begin{dfn} Let $(M,J,\om,\Om)$ be an almost Calabi--Yau
$m$-fold with metric $g$, and define $\psi:M\ra(0,\iy)$
as in \eq{co2eq3}. Let $X$ be a compact SL $m$-fold in
$M$ with conical singularities at $x_1,\ldots,x_n$ with
identifications $\up_i$, cones $C_i$ and rates $\mu_i$.
Define $\Si_i=C_i\cap{\cal S}^{2m-1}$ for $i=1,\ldots,n$.
Let $L_1,\ldots,L_n$ be AC SL $m$-folds in $\C^m$, where
$L_i$ has cone $C_i$ and rate $\la_i$ for $i=1,\ldots,n$.
As in Definition \ref{co2def5}, let $\bigl\{(M,J^s,\om^s,
\Om^s):s\in\F\bigr\}$ be a {\it smooth family of deformations}
of $(M,J,\om,\Om)$, with base space $\F\subset\R^d$. Suppose that
\begin{itemize}
\setlength{\parsep}{0pt}
\setlength{\itemsep}{0pt}
\item[(i)] $L_i$ has rate $\la_i<0$ for $i=1,\ldots,n$, so
that $Y(L_i)=0$ by Proposition~\ref{co4prop}.
\item[(ii)] $[\om^s]\cdot\iota_*(\ga)=0$ for all $s\in\F$ and
$\ga\in H_2(X,\R)$, where $\iota_*:H_2(X,\R)\ra H_2(M,\R)$ is
the natural inclusion.
\item[(iii)] Let $N$ be the compact $m$-manifold obtained
by gluing $L_i$ into $X'$ at $x_i$ for $i=1,\ldots,n$ in the
obvious way. Suppose $N$ is connected. A sufficient
condition for this to hold is that $X$ and $L_i$ for
$i=1,\ldots,n$ are connected.
\end{itemize}

We use the following notation:
\begin{itemize}
\setlength{\parsep}{0pt}
\setlength{\itemsep}{0pt}
\item Let $R,B_R,X'$ and $\iota_i,\Up_i$ for $i=1,\ldots,n$ be
as in Definition~\ref{co3def3}.
\item Let $\ze$ and $U_\sCi,\Phi_\sCi$ for
$i=1,\ldots,n$ be as in Theorem~\ref{co3thm3}.
\item Let $R'$ and $\phi_i,\eta_i,\eta_i^1,\eta_i^2$ for
$i=1,\ldots,n$ be as in Theorem~\ref{co3thm4}.
\item Let $q=b^0(X')$, and write the connected components of
$X'$ as~$X_1',\ldots,X_q'$.
\item For $i=1,\ldots,n$ let $l_i=b^0(\Si_i)$, and write the
connected components of $\Si_i$ as~$\Si_i^1,\ldots,\Si_i^{\smash{l_i}}$.
\item Then $\Up_i\circ\vp_i\bigl(\Si_i^j\t(0,R')\bigr)$ is a
connected subset of $X'$, and so lies in exactly one $X_k'$.
Define $k(i,j)$ for $i=1,\ldots,n$ and $j=1,\ldots,l_i$
by~$\Up_i\circ\vp_i\bigl(\Si_i^j\t(0,R')\bigr)\subset
X'_{\smash{k(i,j)}}$.
\item Let $U_\sXp,\Phi_\sXp$ be as in Theorem~\ref{co3thm6}.
\item Let $A_i$ be as in Theorem \ref{co3thm5} for $i=1,\ldots,n$,
so that~$\eta_i=\d A_i$.
\item Apply Theorem \ref{co4thm1} to $L_i$ with $\ze,U_\sCi,
\Phi_\sCi$ as above, for $i=1,\ldots,n$. Let $T>0$ be as in
the theorem, the same for all $i$. Let the subset $K_i\subset L_i$,
the diffeomorphism $\vp_i:\Si_i\t(T,\iy)\ra L_i\sm K_i$ and the
1-form $\chi_i$ on $\Si_i\t(T,\iy)$ with components $\chi_i^1,
\chi_i^2$ be as in Theorem~\ref{co4thm1}.
\item Let $U_\sLi,\Phi_\sLi$ be as in Theorem \ref{co4thm4}
for~$i=1,\ldots,n$.
\item Let $E_i\in C^\iy\bigl(\Si_i\t(T,\iy)\bigr)$ be as in
Theorem \ref{co4thm2} for~$i=1,\ldots,n$.

Note that as $Y(L_i)=0$ by assumption we have $\ga_i=0$ in
Theorem \ref{co4thm2}, and therefore $\chi_i=\d E_i$
for~$i=1,\ldots,n$.
\item Apply Theorems \ref{co3thm7} and \ref{co3thm8} to $X$
and $\bigl\{(M,J^s,\om^s,\Om^s):s\in\F\bigr\}$. By Theorem
\ref{co3thm8}, part (ii) above implies that $[\nu^s]=0$, so
we can take~$\nu^s\equiv 0$.

Then Theorem \ref{co3thm7} gives an open $\F'\subseteq\F$
with $0\in\F'$ and $\psi^s,\up_i^s,\Up_i^s,\Phi_\sXp^s$ for
$s\in\F'$ depending smoothly on $s$, satisfying \eq{co3eq14} and
\end{itemize}
\e
\begin{gathered}
\up_i^0=\up_i,\;\>
\Up_i^0=\Up_i,\;\>
\Phi_\sXp^0=\Phi_\sXp,\;\>
(\up_i^s)^*(\Om)=\psi^s(x_i)^m\Om',\\
\Up_i^s(0)=x_i,\;\>
\d\Up_i^s\vert_0=\up_i^s,\;\>
(\Up_i^s)^*(\om^s)=\om',\;\>\text{and}\;\>
(\Phi^s_\sXp)^*(\om^s)=\hat\om.
\end{gathered}
\label{co7eq1}
\e
\label{co7def1}
\end{dfn}

Instead of just defining $N^t$ for $t\in(0,\de)$ as
in Definition \ref{co6def2}, we define a larger family
$\bigl\{N^{s,t}:s\in\F'$, $t\in(0,\de)\bigr\}$, where
$N^{s,t}$ is Lagrangian in~$(M,\om^s)$.

\begin{dfn} In the situation of Definition \ref{co7def1}, choose a
smooth, increasing function $F:(0,\iy)\ra[0,1]$ with $F(r)\equiv 0$
for $r\in(0,1)$ and $F(r)\equiv 1$ for $r>2$. Write $F'$ for
$\d F/\d r$. Let $\tau\in(0,1)$ satisfy
\e
0<\max_{i=1,\ldots,n}\Bigl\{\frac{m}{m+1},
\frac{m+2}{2\mu_i+m-2}\Bigr\}<\tau<1,
\label{co7eq2}
\e
which is possible as $\mu_i>2$ implies~$(m+2)/(2\mu_i+m-2)<1$.

For $i=1,\ldots,n$, $s\in\F'$ and small $t>0$, define $P_i^{s,t}=
\Up_i^s(tK_i)$. This is well-defined if $tK_i\subset B_R\subset\C^m$,
and is a compact submanifold of $M$ with boundary, diffeomorphic to
$K_i$. As $K_i$ is Lagrangian in $(\C^m,\om')$ and $(\Up_i^s)^*(\om^s)
=\om'$, we see that $P_i^{s,t}$ is {\it Lagrangian} in~$(M,\om^s)$.

For $i=1,\ldots,n$ and $t>0$ with $tT<t^\tau<2t^\tau<R'$,
define a 1-form $\xi_i^t$ on $\Si_i\t(tT,R')$ by 
\e
\begin{split}
\xi_i^t(\si,r)&=\d\bigl[F(t^{-\tau}r)A_i(\si,r)+
t^2(1-F(t^{-\tau}r))E_i(\si,t^{-1}r)\bigr]\\
&=F(t^{-\tau}r)\eta_i(\si,r)+
t^{-\tau}F'(t^{-\tau}r)A_i(\si,r)\d r\\
&\quad+t^2(1-F(t^{-\tau}r))\chi_i(\si,t^{-1}r)-
t^{2-\tau}F'(t^{-\tau}r)E_i(\si,t^{-1}r)\d r.
\end{split}
\label{co7eq3}
\e
Let $\xi_i^{1,t},\xi_i^{2,t}$ be the components of $\xi_i^t$
in $T^*\Si$ and $\R$, as for $\eta_i,\chi_i$ in Theorems
\ref{co3thm4} and \ref{co4thm1}. Note that when $r\ge 2t^\tau$
we have $F(t^{-\tau}r)\equiv 1$ so that $\xi_i^t(\si,r)=\eta_i
(\si,r)$, and when $r\le t^\tau$ we have $F(t^{-\tau}r)\equiv 0$,
so that $\xi_i^t(\si,r)=t^2\chi_i(\si,t^{-1}r)$. Thus $\xi_i^t$
is an exact 1-form which interpolates between $\eta_i(\si,r)$
near $r=R'$ and $t^2\chi_i(\si,t^{-1}r)$ near~$r=tT$.

Choose $\de\in(0,1]$ with $\de T\le\de^\tau<2\de^\tau\le R'$
and $\de K_i\!\subset\!B_R\!\subset\!\C^m$ and
\e
\bmd{\xi_i^t(\si,r)}<\ze r\quad\text{on $\Si_i\t(tT,R')$ for
all $i=1,\ldots,n$ and $t\in(0,\de)$.}
\label{co7eq4}
\e
This is possible as in \cite[Def.~6.2]{Joyc5}. For $i=1,\ldots,n$,
$s\in\F'$ and $t\in(0,\de)$ define $\Xi_i^{s,t}:\Si_i\t(tT,R')\ra M$ by
\e
\Xi_i^{s,t}(\si,r)=\Up_i^s\circ\Phi_\sCi\bigl(\si,r,
\xi_i^{1,t}(\si,r),\xi_i^{2,t}(\si,r)\bigr).
\label{co7eq5}
\e
Define $Q_i^{s,t}=\Xi_i^{s,t}\bigl(\Si_i\t(tT,R')\bigr)$ for
$i=1,\ldots,n$, $s\in\F'$ and $t\in(0,\de)$. As $(\Up_i^s)^*(\om^s)
=\om'$, $\Phi_\sCi^*(\om')=\hat\om$ and $\xi_i^t$ is closed we
see that $(\Xi_i^{s,t})^*(\om^s)\equiv 0$. Thus $Q_i^{s,t}$ is
{\it Lagrangian} in $(M,\om^s)$, and is a noncompact embedded
submanifold diffeomorphic to~$\Si_i\t(tT,R')$.

Recall that $K\subset X'$, and $X'$ is embedded as the zero section
in $U_\sXp\subset T^*X'$. Thus $K\subset U_\sXp$, and $\Phi_\sXp^s$
maps $U_\sXp\ra M$. For $s\in\F'$, define $K^s=\Phi_\sXp^s(K)$. Then
$K^s$ is Lagrangian in $(M,\om^s)$, as $(\Phi_\sXp^s)^*(\om^s)=
\hat\om$, and $K\subset X'$ is Lagrangian in $(U_\sXp,\hat\om)$.
For $s\in\F'$ and $t\in(0,\de)$, define $N^{s,t}$ to be the
disjoint union of $K^s$, $P_1^{s,t},\ldots,P_n^{s,t}$ and
$Q_1^{s,t},\ldots,Q_n^{s,t}$. Then $N^{s,t}$ is {\it Lagrangian}
in $(M,\om^s)$, as $K^s,P_i^{s,t}$ and $Q_i^{s,t}$ are.

Moreover, $N^{s,t}$ is a compact, smooth submanifold of $M$
{\it without boundary}. The proof of this follows
\cite[Def.~6.2]{Joyc5}, with simple changes. In particular,
$\pd P_i^{s,t}$ joins smoothly onto the $\Si_i\t\{tT\}$ end of
$Q_i^{s,t}$ as $\xi_i^t(\si,r)\equiv t^2\chi_i(\si,t^{-1}r)$
near $r=tT$ on $\Si_i\t(tT,R')$, and both $P_i^{s,t},Q_i^{s,t}$
are defined using $\Up_i^s$. Similarly, the $\Si_i\t\{R'\}$ end
of $Q_i^{s,t}$ joins smoothly onto the appropriate component of
$\pd K^s$ as $\xi_i^t\equiv\eta_i$ near $r=R'$ on $\Si_i\t(tT,R')$,
and because of the compatibility \eq{co3eq14} between $\Up_i^s$,
used to define $Q_i^s$, and $\Phi_\sXp^s$, used to define~$K^s$.

Note that $N^{s,t}$ depends smoothly on $s,t$, since $\Up_i^s,
\Phi_\sXp^s$ depend smoothly on $s$. Also, when $s=0$ we have
$\Up_i^0=\Up_i$ and $\Phi_\sXp^0=\Phi_\sXp$ by \eq{co7eq1}, so
we see that $P_i^{0,t}=P_i^t$, $Q_i^{0,t}=Q_i^t$, $K^0=K$ and
$N^{0,t}=N^t$, where $P_i^t,Q_i^t$ and $N^t$ are as
in~\cite[Def.~6.2]{Joyc5}.

Let $h^{s,t}$ be the restriction of $g$ to $N^{s,t}$ for $s\in\F'$
and $t\in(0,\de)$, so that $(N^{s,t},h^{s,t})$ is a compact
Riemannian manifold, which is naturally oriented. Let $\d V^{s,t}$
be the volume form on $N^{s,t}$. As in \eq{co5eq1} we may write
$\Om\vert_{N^{s,t}}=(\psi^s)^m{\rm e}^{\smash{i\th^{s,t}}}\,\d V^{s,t}$
for some phase function ${\rm e}^{\smash{i\th^{s,t}}}$ on $N^{s,t}$.
Write $\ve^{s,t}=(\psi^s)^m\sin\th^{s,t}$, so that
$\Im\Om\vert_{N^{s,t}}=\ve^{s,t}\,\d V^{s,t}$. Then $\ve^{s,t}$
depends smoothly on $s$, as $\Om^s$, $\Up_i^s$ and $\Phi_\sXp^s$~do.
\label{co7def2}
\end{dfn}

Note that Definition \ref{co7def2} does not include the terms in
$\ga_i,T_i$ and $\al$ added in Definition \ref{co6def2}, which
were there to allow $Y(L_i)\ne 0$. Here is the analogue of
Definition \ref{co6def3}, defining vector spaces $W^{s,t}$
which will be $W$ in Definition~\ref{co5def2}.

\begin{dfn} In the situation of Definitions \ref{co7def1} and
\ref{co7def2} define vector spaces $V_i$ of bounded harmonic
functions on $L_i$ and $v^{{\bf c}_i}_{\smash{i}}\in V_i$ for
${\bf c}_i\in\R^{l_i}$ satisfying \eq{co6eq8} for $i=1,\ldots,n$,
as in Definition \ref{co6def3}. We shall define a vector
subspace $W^{s,t}\subset C^\iy(N^{s,t})$ for $s\in\F'$
and $t\in(0,\de)$, with an isomorphism~$W^{s,t}\cong\R^q$.

Fix ${\bf d}=(d_1,\ldots,d_q)\in\R^q$, and as in Definition
\ref{co6def3} define ${\bf c}_i\in\R^{l_i}$ for $i=1,\ldots,n$
depending linearly on $\bf d$. Let $F:(0,\iy)\ra[0,1]$ and
$\tau\in(0,1)$ be as in Definition \ref{co7def2}. Make $\de>0$
smaller if necessary so that $tT<\ha t^\tau$ for all $t\in(0,\de)$.
For $s\in\F'$ and $t\in(0,\de)$, define $w_{\bf d}^{s,t}\in
C^\iy(N^{s,t})$ as follows:
\begin{itemize}
\item[(i)] The subset $K^s\subset N^{s,t}$ has $q$ connected
components $\Phi_\sXp^s(K\cap X_k')$. Define $w_{\bf d}^{s,t}
\equiv d_k$ on $\Phi_\sXp^s(K\cap X_k')$ for~$k=1,\ldots,q$.
\item[(ii)] Define $w_{\bf d}^{s,t}$ on $P_i^{s,t}\subset N^{s,t}$
by $(\Up_i^s\circ t\circ\vp_i)^*(w_{\bf d}^{s,t})\equiv
v^{{\bf c}_i}_{\smash{i}}$ on~$K_i$.
\item[(iii)] Define $w_{\bf d}^{s,t}$ on $Q_i^{s,t}\subset N^{s,t}$ by
\end{itemize}
\e
(\Xi_i^{s,t})^*(w_{\bf d}^{s,t})(\si,r)=\bigl(1-F(2t^{-\tau}r)\bigr)
\vp_i^*(v^{{\bf c}_i}_{\smash{i}})(\si,t^{-1}r)+F(2t^{-\tau}r)c_i^j
\label{co7eq6}
\e
\begin{itemize}
\item[]on $\Si_i^j\t(tT,R')$, for $i=1,\ldots,n$ and~$j=1,\ldots,l_i$.
\end{itemize}

The argument of Definition \ref{co6def3} shows that
$w_{\bf d}^{s,t}$ is smooth, and linear in $\bf d$. Define
$W^{s,t}=\{w_{\bf d}^{s,t}:{\bf d}\in\R^q\}$. Then $W^{s,t}
\subset C^\iy(N^{s,t})$ is a vector subspace isomorphic to
$\R^q$, and $1\in W^{s,t}$. Define $\pi_\sWst:L^2(N^{s,t})
\ra W^{s,t}$ to be projection onto $W^{s,t}$ using the
$L^2$-inner product.
\label{co7def3}
\end{dfn}

\subsection{Estimating $\Im\Om^s\vert_{N^{s,t}}$}
\label{co72}

We now estimate $\Im\Om\vert_{N^{s,t}}$ as in \S\ref{co62},
to show that $N^{s,t},W^{s,t}$ satisfy part (i) of Theorem
\ref{co5thm1}. Following Proposition \ref{co6prop1} we
bound $\ve^{s,t}$ at each point in~$N^{s,t}$.

\begin{prop} In the situation above, making $\F'$ and\/ $\de>0$
smaller if necessary, there exists $C>0$ such that for all\/
$s\in\F'$ and\/ $t\in(0,\de)$ we have
\ea
\md{\ve^{s,t}}&\le C\md{s},\quad \md{\d\ve^{s,t}}\le C\md{s}
\quad\text{on $K^s$,}
\label{co7eq7}\\
\bmd{(\Xi_i^{s,t})^*(\ve^{s,t})}(\si,r)&\le
\begin{cases}
Cr, & r\in(tT,t^\tau],\\
Ct^{\tau(\mu_i-2)}+Ct^{(1-\tau)(2-\la_i)},
\quad 
& r\in(t^\tau,2t^\tau),\\
Cr\md{s}, & r\in[2t^\tau,R'),
\end{cases}
\label{co7eq8}\\
\bmd{(\Xi_i^{s,t})^*(\d\ve^{s,t})}(\si,r)&\le
\begin{cases}
C, & r\in(tT,t^\tau],\\
Ct^{\tau(\mu_i-3)}+Ct^{(1-\tau)(2-\la_i)-\tau}, & r\in(t^\tau,2t^\tau),\\
C\md{s}, & r\in[2t^\tau,R'),
\end{cases}
\label{co7eq9}\\
\text{and}\quad
\md{\ve^{s,t}}&\le Ct,\quad
\md{\d\ve^{s,t}}\le C
\quad\text{on $P_i^{s,t}$ for all\/ $i=1,\ldots,n$.}
\label{co7eq10}
\ea
Here $\md{\,.\,}$ is computed using $h^{s,t}$
or~$(\Xi_i^{s,t})^*(h^{s,t})$.
\label{co7prop1}
\end{prop}

\begin{proof} When $s=0$ we have $N^{0,t}=N^t$, where $N^t$ is as
in \cite[Def.~6.2]{Joyc5}. Therefore \cite[Prop.~6.4]{Joyc5} proves
\eq{co7eq7}--\eq{co7eq10} when $s=0$. Now $K^s=\Phi_\sXp^s(K)$ is
independent of $t$, and so are $h^{s,t}\vert_{K^s}$ and $\ve^{s,t}
\vert_{K^s}$. Also $\ve^{0,t}\vert_{K^0}\equiv 0$ as $K^0=K$ is
special Lagrangian in $(M,J^0,\om^0,\Om^0)$, and $\ve^{s,t}$
depends smoothly on $s$. Therefore by Taylor's Theorem and
compactness of $K^s\cong K$ we see that
\e
\md{\ve^{s,t}}=O\bigl(\md{s}\bigr),\quad
\md{\d\ve^{s,t}}=O\bigl(\md{s}\bigr)
\quad\text{on $K^s$, for small $s\in\F'$.}
\label{co7eq11}
\e

Consider the $m$-form $\psi^s(x_i)^{-m}(\Xi_i^{s,t})^*
(\ve^{s,t}\d V^{s,t})-\psi(x_i)^{-m}(\Xi_i^{0,t})^*
(\ve^{0,t}\d V^{0,t})$ on $\Si_i\t(tT,R')$. By \eq{co7eq5}
this is the pull-back of $\psi^s(x_i)^{-m}(\Up_i^s)^*
(\Im\Om^s)-\psi(x_i)^{-m}\Up_i^*(\Im\Om)$ on $B_R$ under the
map $(\si,r)\mapsto\Phi_\sCi\bigl(\si,r,\xi_i^{1,t}(\si,r),
\xi_i^{2,t}(\si,r)\bigr)$. But since $(\up_i^s)^*(\Om)=
\psi^s(x_i)^m\Om'$, by \eq{co7eq1} we see that
\e
\psi^s(x_i)^{-m}(\Up_i^s)^*(\Im\Om^s)-
\psi(x_i)^{-m}\Up_i^*(\Im\Om)=O\bigl(r\md{s}\bigr)
\quad\text{on $B_R$.}
\label{co7eq12}
\e
Combining all these facts we see that
\e
\bmd{(\Xi_i^{s,t})^*(\ve^{s,t})}=\bigl(1+O(\md{s})\bigr)
\cdot\bmd{(\Xi_i^{0,t})^*(\ve^{0,t})}+O\bigl(r\md{s}\bigr)
\quad\text{on $\Si_i\t(tT,R')$.}
\label{co7eq13}
\e
A similar proof shows that
\e
\bmd{(\Xi_i^{s,t})^*(\d\ve^{s,t})}=\bigl(1+O(\md{s})\bigr)
\cdot\bmd{(\Xi_i^{0,t})^*(\d\ve^{0,t})}+O\bigl(\md{s}\bigr)
\quad\text{on $\Si_i\t(tT,R')$.}
\label{co7eq14}
\e

Since $\bmd{(\Up_i^s)^*(\Im\Om^s)\vert_{tL_i}}=O(r)$ on $tL_i
\cap B_R$ for small $s\in\F'$, as in \eq{co6eq14}, the proofs
of \eq{co6eq16} and \eq{co6eq17} show that for $i=1,\ldots,n$
we have
\e
\md{\ve^{s,t}}=O(t),\quad
\md{\d\ve^{s,t}}=O(1)
\quad\text{on $P_i^{s,t}$, for small $s\in\F'$.}
\label{co7eq15}
\e
Finally, making $\F'$ and $\de>0$ smaller if necessary, for
some $C>0$ we see that equation \eq{co7eq7} follows from
\eq{co7eq11}, equation \eq{co7eq8} from \eq{co7eq8} with $s=0$
and \eq{co7eq13}, equation \eq{co7eq9} from \eq{co7eq9} with
$s=0$ and \eq{co7eq14}, and equation \eq{co7eq10}
from~\eq{co7eq15}.
\end{proof}

Proposition \ref{co6prop2} then immediately generalizes to give:

\begin{prop} For some $C'>0$ and all\/ $s\in\F'$ and\/
$t\in(0,\de)$ we have
\ea
\lnm{\ve^{s,t}}{2m/(m+2)}&\le C'\md{s}+C't^{\tau(1+m/2)}
\sum_{i=1}^n\bigl(t^{\tau(\mu_i-2)}+t^{(1-\tau)(2-\la_i)}\bigr),
\label{co7eq16}\\
\cnm{\ve^{s,t}}{0}&\le C'\md{s}+C'
\sum_{i=1}^n\bigl(t^{\tau(\mu_i-2)}+t^{(1-\tau)(2-\la_i)}\bigr),
\label{co7eq17}\\
\text{and}\quad\lnm{\d\ve^{s,t}}{2m}&\le C'\md{s}+C't^{-\tau/2}
\sum_{i=1}^n\bigl(t^{\tau(\mu_i-2)}+t^{(1-\tau)(2-\la_i)}\bigr),
\label{co7eq18}
\ea
computing norms with respect to the metric $h^{s,t}$ on~$N^{s,t}$.
\label{co7prop2}
\end{prop}

\subsection{Bounding $\lnm{\pi_\sWst(\ve^{s,t})}{1}$ under
conditions on $s,t$}
\label{co73}

Next we estimate $\lnm{\pi_\sWst(\ve^{s,t})}{1}$, using the method
of \cite[\S 7.2]{Joyc5}. First we bound $\int_{N^{s,t}}w_{\bf d}^{s,t}
\ve^{s,t}\,\d V^{s,t}$ for all ${\bf d}\in\R^q$. As $w_{\bf d}^{s,t}
\equiv d_k$ on $\Phi_\sXp^s(K\cap X_k')$ and $\Xi_i^{s,t}(\Si_i^j\t
[t^\tau,R'))$ when $k(i,j)=k$, we see that
\e
\begin{gathered}
\int_{N^{s,t}}\!\!\!
w_{\bf d}^{s,t}\ve^{s,t}\,\d V^{s,t}=
\sum_{i=1}^n\biggl(\,\int_{P_i^{s,t}}\!\!\!\!
w_{\bf d}^{s,t}\ve^{s,t}
\,\d V^{s,t}+\int_{\Xi_i^{s,t}(\Si_i\t(tT,t^\tau))}
\!\!\!\!\!\!\!\!\!\!\!\!\!\!\!\!\!\!\!\!
w_{\bf d}^{\smash{s,t}}\ve^{s,t}\,\d V^{s,t}\biggr)\\
+\sum_{k=1}^qd_k\biggl(\,\int_{\Phi_{X'}^s(K\cap X_k')}
\!\!\!\!\!\!\!\!\!\!\!\!\!\!\!\!
\ve^{s,t}\,\d V^{s,t}
+\!\!\!\!\sum_{\substack{1\le i\le n, \; 1\le j\le l_i: \\
k(i,j)=k}}\int_{\Xi_i^{s,t}(\Si_i^j\t[t^\tau,R'))}
\!\!\!\!\!\!\!\!\!\!\!\!\!\!\!\!
\ve^{s,t}\,\d V^{s,t}\biggr).
\end{gathered}
\label{co7eq19}
\e

Following the proof of \cite[Prop.~7.4]{Joyc5} with trivial
modifications, we find:

\begin{prop} For all\/ $s\in\F'$, $t\in(0,\de)$, ${\bf d}\in\R^q$
and\/ $i=1,\ldots,n$ we have
\e
\int_{P_i^{s,t}}w_{\bf d}^{s,t}\ve^{s,t}\,\d V^{s,t}+
\int_{\Xi_i^{s,t}(\Si_i\t(tT,t^\tau))}w_{\bf d}^{s,t}\ve^{s,t}\,
\d V^{s,t}=O\bigl(\md{{\bf d}}t^{(m+1)\tau}\bigr).
\label{co7eq20}
\e
\label{co7prop3}
\end{prop}

Observe that the closure $\,\smash{\ov{\!X'_k\!}}\,$ of $X_k'$ in $M$
is an $m$-chain in $M$ {\it without boundary}, and thus defines an
integral homology class $[\,\smash{\ov{\!X'_k\!}}\,]\in H_m(M,\Z)$,
with $[X]=\sum_{k=1}^q[\,\smash{\ov{\!X'_k\!}}\,]$. Using this we
can state the analogue of~\cite[Prop.~7.5]{Joyc5}.

\begin{prop} Making $\F'$ smaller if necessary, for all\/ $s\in\F'$,
$t\in(0,\de)$, ${\bf d}\in\R^q$ and\/ $k=1,\ldots,q$ we have
\e
\begin{gathered}
\int_{\Phi_{X'}^s(K\cap X_k')}
\!\!\!\!\!\!\!\!\!\!\!\!\ve^{s,t}\,\d V^{s,t}
+\!\!\sum_{\substack{1\le i\le n, \; 1\le j\le l_i: \\
k(i,j)=k}}\int_{\Xi_i^{s,t}(\Si_i^j\t
[t^\tau,R'))}\!\!\!\!\!\!\!\!\!\!\!\!\!\ve^{s,t}\,\d V^{s,t}=
\bigl[\Im\Om^s\bigr]\cdot\bigl[\,\smash{\ov{\!X'_k\!}}\,]\\
-t^m\!\!\!\!\!\!\!\!\sum_{\substack{1\le i\le n, \; 1\le j\le l_i: \\
k(i,j)=k}}\!\!\!\!\!\!\!\!
\psi(x_i)^mZ(L_i)\cdot[\Si_i^j\,]
+O\bigl(t^{(m+1)\tau}\bigr)+O\bigl(\md{s}t^m\bigr).
\end{gathered}
\label{co7eq21}
\e
Here $[\Im\Om^s]\in H^m(M,\R)$, $[\,\smash{\ov{\!X'_k\!}}\,]\in
H_m(M,\Z)$, $Z(L_i)\in H^{m-1}(\Si_i,\R)$ is as in \S\ref{co41},
and\/~$[\Si_i^j\,]\in H_{m-1}(\Si_i,\Z)$.
\label{co7prop4}
\end{prop}

\begin{proof} As $\ve^{s,t}\d V^{s,t}=\Im\Om^s\vert_{N^{s,t}}$,
the left hand side of \eq{co7eq21} is the integral of $\Im\Om^s$
over the $m$-chain 
\begin{equation*}
Z_k=\Phi_\sXp^s(K\cap X_k')
+\sum_{\substack{1\le i\le n, \; 1\le j\le l_i: \\
k(i,j)=k}}\Xi_i^{s,t}\bigl(\Si_i^j\t[t^\tau,R')\bigr)
\end{equation*}
for $k=1,\ldots,q$, which is a closed subset of $N^{s,t}$,
with boundary $(m\!-\!1)$-chain
\e
\pd Z_k=
-\sum_{\substack{1\le i\le n, \; 1\le j\le l_i: \\
k(i,j)=k}}\Xi_i^{s,t}\bigl(\Si_i^j\t\{t^\tau\}\bigr).
\label{co7eq22}
\e

For each $i=1,\ldots,n$ and $j=1,\ldots,l_i$, define an $m$-chain
$A_i^j$ in $B_R$ to be the image of $\Si_i^j\t[0,1]$
under the map $\Si_i^j\t[0,1]\ra B_R$ given by
\begin{equation*}
(\si,r)\longmapsto r\Phi_\sCi\bigl(\si,t^\tau,t^2\chi_i^1
(\si,t^{\tau-1}),t^2\chi_i^2(\si,t^{\tau-1})\bigr).
\end{equation*}
As $\Up_i^s\circ\Phi_\sCi\bigl(\si,t^\tau,t^2\chi_i^1(\si,t^{\tau-1}),
t^2\chi_i^2(\si,t^{\tau-1})\bigr)\equiv\Xi_i^{s,t}(\si,t^\tau)$ for
$\si\in\Si_i$ by Definition \ref{co7def2}, we see that
\e
\pd\bigl(\Up_i^s(A_i^j)\bigr)=\Xi_i^{s,t}\bigl(\Si_i^j\t\{t^\tau\}\bigr),
\label{co7eq23}
\e
regarding $\Up_i^s(A_i^j)$ as an $m$-chain in~$M$.

Define another $m$-chain $Z_k'$ for $k=1,\ldots,q$ to be
\begin{equation*}
Z_k'=\!\ov{\!X_k'\!}\,-\sum_{\substack{1\le i\le n, \;
1\le j\le l_i: \\ k(i,j)=k}}\Up_i^s(A_i^j).
\end{equation*}
As $\,\ov{\!X_k'\!}\,$ is an $m$-chain without boundary, we see
from \eq{co7eq22} and \eq{co7eq23} that $\pd Z_k'=\pd Z_k$, and
in fact it is easy to see that $Z_k'$ and $Z_k$ are homologous
in $M$. Since $\Im\Om^s$ is a closed $m$-form on $M$, this
implies that $\int_{Z_k'}\Im\Om^s=\int_{Z_k}\Im\Om^s$.

From Theorem \ref{co3thm7} we have $(\up_i^s)^*(\Om^s)=\psi^s
(x_i)^m\Om'$, where $\Om'$ is as in \eq{co2eq1}. Thus as
$(\Up_i^s)^*(\Om^s)$ is smooth on $B_R$, Taylor's theorem gives
\e
(\Up_i^s)^*(\Om^s)=\psi^s(x_i)^m\Om'+O(r)\quad\text{on $B_R$.}
\label{co7eq24}
\e
Making $\F'$ smaller if necessary, this holds {\it uniformly in}
$s$ for~$s\in\F'$.

Now $A_i^j$ is an $m$-chain in $B_R\subset\C^m$ with boundary in
the AC SL $m$-fold $tL_i$, and $[\pd A_i^j]\in H_{m-1}(tL_i,\R)$
is the image of $[\Si_i^j]\in H_{m-1}(\Si_i,\R)$ under the map
$H_{m-1}(\Si_i,\R)\ra H_{m-1}(L_i,\R)$ dual to the map $H^{m-1}
(L_i,\R)\ra H^{m-1}(\Si_i,\R)$ of \eq{co4eq2}. It then follows
easily from Definition \ref{co4def2} and Lemma \ref{co4lem} that
\e
\int_{A_i^j}\Im\Om'=Z(tL_i)\cdot[\Si_i^j]=t^mZ(L_i)\cdot[\Si_i^j].
\label{co7eq25}
\e

But as $r=O(t^\tau)$ on $A_i^j$ and $\vol(A_i^j)=O(t^{m\tau})$
we see from \eq{co7eq24} that
\e
\int_{A_i^j}\bigl((\Up_i^s)^*(\Im\Om^s)-\psi^s(x_i)^m\Im\Om'\bigr)
=O(t^{(m+1)\tau}),
\label{co7eq26}
\e
which again holds uniformly in $s\in\F'$. As the left hand side
of \eq{co7eq21} is $\int_{Z_k}\Im\Om^s=\int_{Z_k'}\Im\Om^s$,
equation \eq{co7eq21} follows from $\int_{X_k'}\Im\Om^s=
[\Im\Om^s]\cdot[\,\smash{\ov{\!X'_k\!}}\,]$, equations
\eq{co7eq25}, \eq{co7eq26} and $\psi^s(x_i)=\psi(x_i)+O(\md{s})$.
This completes the proof.
\end{proof}

Now we can estimate $\lnm{\pi_\sWst(\ve^{s,t})}{1}$ under
conditions on~$s,t$.

\begin{prop} Making $\F',\de$ smaller if necessary, there
exists $C''>0$ such that for all\/ $s\in\F'$ and\/ $t\in(0,\de)$
satisfying
\e
\bigl[\Im\Om^s\bigr]\cdot\bigl[\,\smash{\ov{\!X'_k\!}}\,]=t^m
\!\!\!\!\!\!\!\!\sum_{\substack{1\le i\le n, \; 1\le j\le l_i: \\
k(i,j)=k}}\!\!\!\!\!\!\!
\psi(x_i)^mZ(L_i)\cdot[\Si_i^j\,]
\quad\text{for $k=1,\ldots,q,$}
\label{co7eq27}
\e
we have~$\lnm{\pi_\sWst(\ve^{s,t})}{1}\le C''t^{(m+1)\tau}+C''\md{s}t^m$.
\label{co7prop5}
\end{prop}

\begin{proof} Let $s\in\F'$ and $t\in(0,\de)$ satisfy \eq{co7eq27}.
Combining equations \eq{co7eq19}, \eq{co7eq20} \eq{co7eq21} and
\eq{co7eq27} gives
\e
\int_{N^{s,t}}w_{\bf d}^{s,t}\ve^{s,t}\,\d V^{s,t}=
O\bigl(\md{{\bf d}}t^{(m+1)\tau}\bigr)+
O\bigl(\md{{\bf d}}\md{s}t^m\bigr)
\quad\text{for all ${\bf d}\in\R^q$.}
\label{co7eq28}
\e
One can show from Definition \ref{co7def3} that
$\lnm{w_{\bf d}^{s,t}}{2}\ge C\md{{\bf d}}$ for some $C>0$
and all ${\bf d},s,t$. This and \eq{co7eq28} imply that
$\lnm{\pi_\sWst(\ve^{s,t})}{2}=O(t^{(m+1)\tau})
+O(\md{s}t^m)$. But $\lnm{\pi_\sWst(\ve^{s,t})}{1}\le
\vol(N^{s,t})^{1/2}\lnm{\pi_\sWst(\ve^{s,t})}{2}$, and
$\vol(N^{s,t})=O(1)$. Therefore $\lnm{\pi_\sWst(\ve^{s,t})}{1}
=O(t^{(m+1)\tau})+O(\md{s}t^m)$. So making $\F',\de$ smaller if
necessary, there exists $C''>0$ such that the proposition holds.
\end{proof}

If we do not impose equation \eq{co7eq27} then the $t^m$ term in
\eq{co7eq21} means that we expect $\lnm{\pi_\sWst(\ve^{s,t})}{1}
=O(t^m)$, which is too big for part (i) of Theorem \ref{co5thm1}
to hold for $N^{s,t},W^{s,t}$. We can now prove the analogue of
Theorem~\ref{co6thm1}.

\begin{thm} In the situation of Definitions
\ref{co7def1}--\ref{co7def3}, making $\F'$ and\/ $\de>0$
smaller if necessary, there exist\/ $A_2>0$ and\/ $\ka>1$
such that for all\/ $s\in\F'$ and\/ $t\in(0,\de)$ with\/
$\md{s}\le t^{\ka+m/2}$ and
\e
\bigl[\Im\Om^s\bigr]\cdot\bigl[\,\smash{\ov{\!X'_k\!}}\,]=t^m
\!\!\!\!\!\!\!\!\sum_{\substack{1\le i\le n, \; 1\le j\le l_i: \\
k(i,j)=k}}\!\!\!\!\!\!\!
\psi(x_i)^mZ(L_i)\cdot[\Si_i^j\,]
\quad\text{for $k=1,\ldots,q,$}
\label{co7eq29}
\e
the functions $\ve^{s,t}=(\psi^s)^m\sin\th^{s,t}$ on $N^{s,t}$
satisfy $\lnm{\ve^{s,t}}{2m/(m+2)}\le A_2t^{\ka+m/2}$,
$\cnm{\ve^{s,t}}{0}\le A_2t^{\ka-1}$, $\lnm{\d\ve^{s,t}}{2m}\le
A_2t^{\ka-3/2}$ and\/ $\lnm{\pi_\sWst(\ve^{s,t})}{1}\le A_2
t^{\ka+m-1}$, as in part\/ {\rm(i)} of Theorem~\ref{co5thm1}.
\label{co7thm1}
\end{thm}

\begin{proof} Make $\F',\de$ smaller if necessary so that
Propositions \ref{co7prop1}--\ref{co7prop5} hold. Let $C',C''$
be as in Propositions \ref{co7prop2} and \ref{co7prop5}, and
set $A_2\!=\!\max\bigl((2n\!+\!1)C',2C''\bigr)\!>\!0$. Let
$s\in\F'$ and $t\in(0,\de)$ satisfy $\md{s}\le t^{\ka+m/2}$
and \eq{co7eq29}. We seek $\ka>1$ so that the four bounds on
$\ve^{s,t}$ hold.

As $t<1$, $\md{s}\le t^{\ka+m/2}$ and $A_2\ge(2n+1)C'$, equations
\eq{co7eq16}--\eq{co7eq18} imply that $\lnm{\ve^{s,t}}{2m/(m+2)}
\!\le\!A_2t^{\ka+m/2}$, $\cnm{\ve^{s,t}}{0}\!\le\!A_2t^{\ka-1}$
and $\lnm{\d\ve^{s,t}}{2m}\!\le\!A_2t^{\ka-3/2}$ if
\ea
\tau(1\!+\!m/2)\!+\!\tau(\mu_i\!-\!2)&\!\ge\!\ka\!+\!m/2,&\;\>
\tau(1\!+\!m/2)\!+\!(1\!-\!\tau)(2\!-\!\la_i)&\!\ge\!\ka\!+\!m/2,
\label{co7eq30}\\
\tau(\mu_i\!-\!2)&\!\ge\!\ka\!-\!1,&\;\>
(1\!-\!\tau)(2\!-\!\la_i)&\!\ge\!\ka\!-\!1,
\label{co7eq31}\\
-\tau/2\!+\!\tau(\mu_i\!-\!2)&\ge\ka\!-\!3/2,
&\text{and}\;\>
-\tau/2\!+\!(1\!-\!\tau)(2\!-\!\la_i)&\ge\ka\!-\!3/2
\label{co7eq32}
\ea
for all $i=1,\ldots,n$. As $A_2\ge 2C''$ and $t<1$, by
Proposition \ref{co7prop5} we have $\lnm{\pi_\sWst(\ve^{s,t})}{1}
\le A_2t^{\ka+m-1}$ provided
\e
(m+1)\tau\ge\ka+m-1.
\label{co7eq33}
\e

Since $0<\tau<1$ and $\mu_i>2$ the first equation of \eq{co7eq30}
admits a solution $\ka>1$ provided $\tau>(2+m)/(2\mu_i-2+m)$,
which holds by \eq{co7eq2}. Since $\la_i<0$ by Definition
\ref{co7def1}, we can suppose $\la_i<\ha(2-m)$ by part (a) of
Theorem \ref{co4thm3} applied to $L_i$, and then the second
equation of \eq{co7eq30} admits a solution $\ka>1$. Also
\eq{co7eq30} implies \eq{co7eq31} and \eq{co7eq32} as
$\tau\le 1$. Equation \eq{co7eq33} admits a solution $\ka>1$
as $\tau>\frac{m}{m+1}$ by \eq{co7eq2}. Thus (having decreased
$\la_i$ if necessary) we can choose $\ka>1$ satisfying
\eq{co7eq30}--\eq{co7eq33}, and the theorem is proved.
\end{proof}

Here we have restricted to pairs $s\in\F'$ and $t\in(0,\de)$ with
$\md{s}\le t^{\ka+m/2}$. This is so that the contributions $C'\md{s}$
in \eq{co7eq16}--\eq{co7eq18} can be absorbed into $A_2t^{\ka+m/2}$,
$A_2t^{\ka-1}$ and $A_2t^{\ka-3/2}$. Basically, Theorem \ref{co7thm1}
says that if $s$ is not too large compared to $t$, and \eq{co7eq29}
holds, then part (i) of Theorem \ref{co5thm1} holds for~$N^{s,t},W^{s,t}$.

\subsection{Parts (ii)--(vii) of Theorem \ref{co5thm1}}
\label{co74}

Next we carry out the programme of \S\ref{co63} and
\S\ref{co64} for families. Here is the analogue of
Definition \ref{co6def4} and~\cite[Def.~6.7]{Joyc5}.

\begin{dfn} For each $s\in\F'$ and $t\in(0,\de)$, we define
an open neighbourhood $U_\sNst\subset T^*N^{s,t}$ of the zero
section $N^{s,t}$ in $T^*N^{s,t}$, and a smooth map $\Phi_\sNst:
U_\sNst\ra M$. Let $\pi:T^*N^{s,t}\ra N^{s,t}$ be the projection.
Define
\e
\begin{gathered}
U_\sNst\cap\pi^*(K^s)=
\d(\Phi_\sXp^s\vert_K)\bigl(U_\sXp\cap\pi^*(K)\bigr)
\quad\text{and}\\
\Phi_\sNst\vert_{U_\sNst\cap\pi^*(K^s)}\circ\d(\Phi_\sXp^s
\vert_K)=\Phi_\sXp^s\vert_{U_\sXp\cap\pi^*(K)},
\end{gathered}
\label{co7eq34}
\e
where $\Phi_\sXp^s\vert_K:K\ra K^s$ is a diffeomorphism and
$\d(\Phi_\sXp^s\vert_K):T^*K\ra T^*K^s$ the induced isomorphism.

Following \eq{co6eq34}, define $U_\sNst\cap\pi^*(P_i^{s,t})$
and $\Phi_\sNst\vert_{U_\sNst\cap\pi^*(P_i^{s,t})}$ by
\e
\begin{gathered}
U_\sNst\cap\pi^*(P_i^{s,t})=\d(\Up_i^s\circ t)
\bigl(\{\al\in T^*K_i:t^{-2}\al\in U_\sLi\}\bigr)\\
\text{and}\quad \Phi_\sNst\circ\d(\Up_i^s\circ t)(\al)=
\Up_i^s\circ t\circ\Phi_\sLi(t^{-2}\al).
\end{gathered}
\label{co7eq35}
\e
Modifying \eq{co6eq35} and \eq{co6eq36}, define $U_\sNst\cap\pi^*
(Q_i^{s,t})$ and $\Phi_\sNst\vert_{U_\sNst\cap\pi^*(Q_i^{s,t})}$ by
\begin{gather}
(\d\Xi_i^{s,t})^*(U_\sNst)=\bigl\{(\si,r,\varsigma,u)
\in T^*\bigl(\Si_i\t(tT,R')\bigr):\bmd{(\varsigma,u)}<\ze r\bigr\}
\quad\text{and}
\label{co7eq36}\\
\Phi_\sNst\circ\d\Xi_i^{s,t}(\si,r,\varsigma,u)\equiv\Up_i^s\circ\Phi_\sCi
\bigl(\si,r,\varsigma+\xi_i^{1,t}(\si,r),u+\xi_i^{2,t}(\si,r)\bigr).
\label{co7eq37}
\end{gather}

The proof in \cite[Def.~6.7]{Joyc5} shows that making $\F',\de$
smaller if necessary, $U_\sNst$ and $\Phi_\sNst$ are well-defined
for all $s\in\F'$ and $t\in(0,\de)$, and $U_\sNst$ is an open
tubular neighbourhood of $N^{s,t}$ in $T^*N^{s,t}$. Since
$(\Up_i^s)^*(\om^s)=\om'$ and $(\Phi_\sXp^s)^*(\om^s)=\hat\om$ we
also find that $\Phi_\sNst^*(\om^s)=\hat\om$. Define an $m$-form
$\be^{s,t}$ on $U_\sNst$ by~$\be^{s,t}=\Phi_\sNst^*(\Im\Om^s)$.
\label{co7def4}
\end{dfn}

Here is the analogue of Theorem~\ref{co6thm2}.

\begin{thm} Making $\F',\de$ smaller if necessary,
there exist\/ $A_1,A_3,\ldots,A_6>0$ such that for all\/
$s\in\F'$ and\/ $t\in(0,\de)$, as in {\rm(ii)--(v)} of
Theorem \ref{co5thm1} we have
\begin{itemize}
\item[{\rm(ii)}] $\psi^s\ge A_3$ on $N^{s,t}$.
\item[{\rm(iii)}] The subset\/ $\B_{A_1t}\!\subset\!T^*N^{s,t}$
of Definition \ref{co5def2} lies in $U_\sNst$, and
$\cnm{\hat\na^k\be^{s,t}}{0}
\allowbreak
\!\le\!A_4t^{-k}$ on $\B_{A_1t}$ for
$k=0,1,2$ and\/~$3$.
\item[{\rm(iv)}] The injectivity radius $\de(h^{s,t})$
satisfies~$\de(h^{s,t})\ge A_5t$.
\item[{\rm(v)}] The Riemann curvature $R(h^{s,t})$
satisfies~$\cnm{R(h^{s,t})}{0}\le A_6t^{-2}$.
\end{itemize}
Here part\/ {\rm(iii)} uses the notation of Definition \ref{co5def2},
and parts {\rm(iv)} and\/ {\rm(v)} refer to the compact Riemannian
manifold\/~$(N^{s,t},h^{s,t})$.
\label{co7thm2}
\end{thm}

Note that here we do not require $s$ to be small compared to $t$,
as we did in Theorem \ref{co7thm1}, but instead parts (ii)--(v)
hold uniformly over $s\in\F'$. The proof is a straightforward
extension of that of Theorem \ref{co6thm2} to families. For any
given $s\in\F'$ the previous proof gives $\de,A_1,A_3,\ldots,
A_6>0$, which can be chosen to depend continuously on $s$. Then
making $\F'$ smaller, we can take $\de,A_1,A_3,\ldots,A_6$ to
be independent of~$s$.

Making $\F'$ smaller if necessary, we find that all of \cite[\S 6.4 \&
\S 7.3]{Joyc5} can be extended to $(M,J^s,\om^s,\Om^s)$ for $s\in\F'$,
and applies uniformly for $s\in\F'$. Thus with only minor modifications
we prove versions of Theorems \ref{co6thm3} and~\ref{co6thm4}:

\begin{thm} Making $\F'$ and\/ $\de>0$ smaller if necessary, there
exists $A_7>0$ such that for all\/ $s\in\F'$ and\/ $t\in(0,\de)$,
if\/ $v\in L^2_1(N^{s,t})$ with\/ $\int_{N^{s,t}}vw\,\d V^{s,t}=0$
for all\/ $w\in W^{s,t}$ then $v\in L^{2m/(m-2)}(N^{s,t})$
and\/~$\lnm{v}{2m/(m-2)}\le A_7\lnm{\d v}{2}$.
\label{co7thm3}
\end{thm}

\begin{thm} Making $\F'$ and\/ $\de>0$ smaller if necessary, for
all\/ $s\in\F'$, $t\in(0,\de)$ and\/ $w\in W^{s,t}$ we have
$\lnm{\d^*\d w}{2m/(m+2)}\le\ha A_7^{-1}\lnm{\d w}{2}$, where
$A_7>0$ is as in Theorem \ref{co7thm3}. Also there exists
$A_8>0$ such that for all\/ $s\in\F'$, $t\in(0,\de)$ and\/
$w\in W^{s,t}$ with\/ $\int_{N^{s,t}}w\,\d V^{s,t}=0$ we
have~$\cnm{w}{0}\le A_8t^{1-m/2}\lnm{\d w}{2}$.
\label{co7thm4}
\end{thm}

\subsection{The main result for families when $Y(L_i)=0$}
\label{co75}

Here is our first main result for desingularization in families,
the analogue of Theorem \ref{co5thm3}, when we do not assume
that $X'$ is connected.

\begin{thm} Suppose $(M,J,\om,\Om)$ is an almost Calabi--Yau
$m$-fold and\/ $X$ a compact SL\/ $m$-fold in $M$ with conical
singularities at\/ $x_1,\ldots,x_n$ and cones $C_1,\ldots,C_n$.
Define $\psi:M\ra(0,\iy)$ as in \eq{co2eq3}. Let\/ $L_1,\ldots,L_n$
be Asymptotically Conical SL\/ $m$-folds in $\C^m$ with cones
$C_1,\ldots,C_n$ and rates $\la_1,\ldots,\la_n$. Suppose
$\la_i<0$ for $i=1,\ldots,n$. Write $X'=X\sm\{x_1,\ldots,x_n\}$
and\/~$\Si_i=C_i\cap{\cal S}^{2m-1}$.

Set\/ $q=b^0(X')$, and let\/ $X_1',\ldots,X_q'$ be the connected
components of\/ $X'$. For $i=1,\ldots,n$ let\/ $l_i=b^0(\Si_i)$,
and let\/ $\Si_i^1,\ldots,\Si_i^{\smash{l_i}}$ be the connected
components of\/ $\Si_i$. Define $k(i,j)=1,\ldots,q$ by $\Up_i
\circ\vp_i\bigl(\Si_i^j\t(0,R')\bigr)\subset X'_{\smash{k(i,j)}}$
for $i=1,\ldots,n$ and\/ $j=1,\ldots,l_i$. Suppose the compact\/
$m$-manifold\/ $N$ obtained by gluing $L_i$ into $X'$ at\/ $x_i$
for $i=1,\ldots,n$ is connected. A sufficient condition for this
to hold is that\/ $X$ and\/ $L_i$ for $i=1,\ldots,n$ are connected.

Suppose $\bigl\{(M,J^s,\om^s,\Om^s):s\in\F\bigr\}$ is a smooth
family of deformations of\/ $(M,J,\om,\Om)$, with base space
$\F\subset\R^d$. Let\/ $\iota_*:H_2(X,\R)\ra H_2(M,\R)$ be
the natural inclusion. Suppose that
\e
[\om^s]\cdot\iota_*(\ga)=0
\quad\text{for all\/ $s\in\F$ and\/ $\ga\in H_2(X,\R)$.}
\label{co7eq38}
\e
Define $\G\subseteq\F\t(0,1)$ to be the subset of\/
$(s,t)\in\F\t(0,1)$ with
\e
\bigl[\Im\Om^s\bigr]\cdot\bigl[\,\smash{\ov{\!X'_k\!}}\,]=t^m
\!\!\!\!\!\!\!\!\sum_{\substack{1\le i\le n, \; 1\le j\le l_i: \\
k(i,j)=k}}\!\!\!\!\!\!\!
\psi(x_i)^mZ(L_i)\cdot[\Si_i^j\,]
\quad\text{for $k=1,\ldots,q$.}
\label{co7eq39}
\e

Then there exist\/ $\ep\in(0,1)$ and\/ $\ka>1$ and a smooth family
\e
\bigl\{\smash{\ti N^{s,t}}:(s,t)\in\G,\quad t\in(0,\ep],
\quad \md{s}\le t^{\ka+m/2}\bigr\},
\label{co7eq40}
\e
such that\/ $\smash{\ti N^{s,t}}$ is a compact, nonsingular
SL\/ $m$-fold in $(M,J^s,\om^s,\Om^s)$ diffeomorphic to $N$,
which is constructed by gluing $tL_i$ into $X$ at\/ $x_i$ for
$i=1,\ldots,n$. In the sense of currents in Geometric Measure
Theory, $\smash{\ti N^{s,t}}\ra X$ as~$s,t\ra 0$.
\label{co7thm5}
\end{thm}

\begin{proof} The hypotheses of the theorem imply that
conditions (i)--(iii) of Definition \ref{co7def1} hold.
Let $\F',\de$ and $N^{s,t}$ for $s\in\F'$ and $t\in(0,\de)$
be as in Definition \ref{co7def2} and $W^{s,t}$ as in
Definition \ref{co7def3}, and make $\F'$ and $\de>0$
smaller if necessary so that Theorems \ref{co7thm1},
\ref{co7thm2}, \ref{co7thm3} and \ref{co7thm4} apply.
Theorems \ref{co7thm1}, \ref{co7thm2}, \ref{co7thm3} and
\ref{co7thm4} now give $\ka>1$ and $A_1,\ldots,A_8>0$.
We show that $\cos\th^{s,t}\ge\ha$ on $N^{s,t}$ for
small $s,t$ as in the proof of Theorem~\ref{co6thm5}.

Let $\ep,K>0$ be as given in Theorem \ref{co5thm1} depending
on $\ka,A_1,\ldots,A_8$ and $m$, and make $\ep>0$ smaller if
necessary to ensure that $\ep<\de$, and that if $s\in\F$ with
$\md{s}\le\ep^{\ka+m/2}$ then $s\in\F'$ and $\cos\th^{s,t}\ge\ha$
on $N^{s,t}$ for $t\in(0,\ep]$. Let $\G$ be as in the theorem,
and suppose $(s,t)\in\G$ with $t\in(0,\ep]$ and $\md{s}\le
t^{\ka+m/2}$. Then~$s\in\F'$, as $\md{s}\le t^{\ka+m/2}\le
\ep^{\ka+m/2}$.

Theorem \ref{co7thm1} shows that part (i) of Theorem
\ref{co5thm1} holds for $N^{s,t},W^{s,t}$, as $\md{s}\le
t^{\ka+m/2}$ and \eq{co7eq29} holds by choice of $s,t$.
Theorems \ref{co7thm2}, \ref{co7thm3} and \ref{co7thm4}
show that parts (ii)--(v), (vi) and (vii) of Theorem
\ref{co5thm1} hold for $N^{s,t},W^{s,t}$. Thus as
$t\le\ep$, Theorem \ref{co5thm1} shows that there
exists a nearby compact, nonsingular SL $m$-fold
$\smash{\ti N^{s,t}}$ in $(M,J^s,\om^s,\Om^s)$, as in
\eq{co7eq40}. The remaining conclusions follow as for
Theorem~\ref{co6thm5}.
\end{proof}

Putting $d=0$ and $\F=\{0\}=\R^d$, Theorem \ref{co7thm5}
reduces to Theorem \ref{co5thm3}. Equation \eq{co7eq38} is
a necessary condition for the existence of any Lagrangian
$m$-fold $\smash{\ti N^{s,t}}$ in $(M,J^s,\om^s,\Om^s)$
desingularizing $X$. However, \eq{co7eq39} cannot be
justified in the same way, and actually it's sufficient for
$[\Im\Om^s]\cdot[X]=0$ to hold exactly (this follows from
the sum of \eq{co7eq39} over $k=1,\ldots,q$) and for 
\eq{co7eq39} to hold only approximately, that is, up to
terms of order~$O(t^{\ka+m-1})$.

If $X'$ is connected, so that $q=1$, then as for Theorem
\ref{co6thm6} the right hand side of \eq{co7eq39} is zero
automatically. Then the theorem simplifies to give:

\begin{thm} Let\/ $(M,J,\om,\Om)$ be an almost Calabi--Yau
$m$-fold and\/ $X$ a compact SL\/ $m$-fold in $M$ with conical
singularities at\/ $x_1,\ldots,x_n$ and cones $C_1,\ldots,C_n$.
Let\/ $L_1,\ldots,L_n$ be Asymptotically Conical SL\/ $m$-folds
in $\C^m$ with cones $C_1,\ldots,C_n$ and rates $\la_1,\ldots,
\la_n$. Suppose $\la_i<0$ for $i=1,\ldots,n$, and\/ $X'=X\sm
\{x_1,\ldots,x_n\}$ is connected.

Suppose $\bigl\{(M,J^s,\om^s,\Om^s):s\in\F\bigr\}$ is a smooth
family of deformations of\/ $(M,J,\om,\Om)$, with base space
$\F\subset\R^d$. Let\/ $\iota_*:H_2(X,\R)\ra H_2(M,\R)$ be
the natural inclusion. Suppose that
\e
\begin{split}
[\om^s]\cdot\iota_*(\ga)&=0
\quad\text{for all\/ $s\in\F$ and\/ $\ga\in H_2(X,\R)$, and}\\
[\Im\Om^s]\cdot[X]&=0
\quad\text{for all\/ $s\in\F$, where $[X]\in H_m(M,\R)$.}
\end{split}
\label{co7eq41}
\e

Then there exist\/ $\ep>0$ and\/ $\ka>1$ and a smooth family
\e
\bigl\{\smash{\ti N^{s,t}}:s\in\F,\quad t\in(0,\ep],\quad
\md{s}\le t^{\ka+m/2}\bigr\},
\label{co7eq42}
\e
such that\/ $\smash{\ti N^{s,t}}$ is a compact, nonsingular
SL\/ $m$-fold in $(M,J^s,\om^s,\Om^s)$, which is constructed
by gluing $tL_i$ into $X$ at\/ $x_i$ for $i=1,\ldots,n$. In
the sense of currents in Geometric Measure Theory,
$\smash{\ti N^{s,t}}\ra X$ as~$s,t\ra 0$.
\label{co7thm6}
\end{thm}

Putting $d=0$ and $\F=\{0\}=\R^d$, Theorem \ref{co7thm6}
reduces to Theorem \ref{co5thm2}. When $s=0$ the SL
$m$-fold $\smash{\ti N^{0,t}}$ in Theorem \ref{co7thm6}
coincides with the SL $m$-fold $\smash{\ti N^t}$ in
$(M,J,\om,\Om)$ constructed in Theorem \ref{co5thm2}. Now
given the compact, nonsingular SL $m$-folds $\smash{\ti N^t}$
in Theorem \ref{co5thm2}, we can apply Theorem \ref{co2thm2}
to prove the existence of SL $m$-folds $\smash{\ti N^{s,t}}$
in $(M,J^s,\om^s,\Om^s)$ for $s$ small compared to~$t$.

Thus, much of Theorem \ref{co7thm6} follows from Theorems
\ref{co2thm2} and \ref{co5thm2}. The new feature is that
when $t$ is small $\md{s}\le t^{\ka+m/2}$ is sufficient
for the existence of $\smash{\ti N^{s,t}}$, whereas
Theorems \ref{co2thm2} and \ref{co5thm2} give no
quantitative restrictions on~$s,t$.

However, Theorem \ref{co7thm5} does {\it not\/} follow from
Theorems \ref{co2thm2} and \ref{co5thm3}. This is because
Theorem \ref{co7thm5} can prove the existence of desingularizations
$\smash{\ti N^{s,t}}$ in $(M,J^s,\om^s,\Om^s)$ for $s\ne 0$ in
cases when there do not exist desingularizations $\smash{\ti N^t}$
in $(M,J,\om,\Om)$ because \eq{co5eq2} does not hold, so Theorem
\ref{co5thm3} does not apply. Effectively, by deforming
$(M,J,\om,\Om)$ in a family $\F$ we can overcome obstructions
to desingularizing $X$ in the single almost Calabi--Yau
$m$-fold~$(M,J,\om,\Om)$.

\section{Desingularizing in families when $Y(L_i)\ne 0$}
\label{co8}

Finally we extend the material of \S\ref{co6} to families of almost
Calabi--Yau $m$-folds $(M,J^s,\om^s,\Om^s)$ for $s\in\F$. There are
{\it topological obstructions} to defining a Lagrangian $m$-fold
$N^{s,t}$ in $(M,\om^s)$ by gluing $tL_i$ into $X$ at $x_i$, as
$[\om^s\vert_{N^{s,t}}]$ may be nonzero in $H^2(N^{s,t},\R)$.
Because of this, we can only define $N^{s,t}$ as a Lagrangian
$m$-fold for $(s,t)$ in a subset $\G$ of $\F\t(0,1)$, satisfying
an equation \eq{co8eq9} involving $[\om^s]$, $t$ and $Y(L_i)$
for~$i=1,\ldots,n$.

\subsection{Setting up the problem}
\label{co81}

We shall consider the following situation, combining Definitions
\ref{co6def1} and~\ref{co7def1}.

\begin{dfn} Let $(M,J,\om,\Om)$ be an almost Calabi--Yau
$m$-fold, and $X$ a compact SL $m$-fold in $M$ with conical
singularities at $x_1,\ldots,x_n$ with identifications
$\up_i$, cones $C_i$ and rates $\mu_i$. Let $L_1,\ldots,L_n$
be AC SL $m$-folds in $\C^m$, where $L_i$ has cone $C_i$ and
rate $\la_i$ for $i=1,\ldots,n$. As in Definition \ref{co2def5},
let $\bigl\{(M,J^s,\om^s,\Om^s):s\in\F\bigr\}$ be a {\it smooth
family of deformations} of $(M,J,\om,\Om)$, with base
space~$\F\subset\R^d$.

Set $q=b^0(X')$, so that $X'$ has $q$ connected components,
and number them $X_1',\ldots,X_q'$. For $i=1,\ldots,n$ let
$l_i=b^0(\Si_i)$, so that $\Si_i$ has $l_i$ connected
components, and number them $\Si_i^1,\ldots,\Si_i^{\smash{l_i}}$.
If $\Up_i,\vp_i,S_i$ are as in Definition \ref{co3def3}, then
$\Up_i\circ\vp_i$ is a diffeomorphism $\Si_i\t(0,R')\ra S_i
\subset X'$. For each $j=1,\ldots,l_i$, $\Up_i\circ\vp_i\bigl(
\Si_i^j\t(0,R')\bigr)$ is a connected subset of $X'$, and so
lies in exactly one of the $X_k'$ for $k=1,\ldots,q$. Define
numbers $k(i,j)=1,\ldots,q$ for $i=1,\ldots,n$ and $j=1,\ldots,l_i$.
by $\Up_i\circ\vp_i(\Si_i^j\t(0,R'))\subset X'_{\smash{k(i,j)}}$.
Suppose that:
\begin{itemize}
\setlength{\parsep}{0pt}
\setlength{\itemsep}{0pt}
\item[(i)] The dimension $m$ satisfies $2<m<6$, and
\item[(ii)] The AC SL $m$-fold $L_i$ has rate $\la_i\le 0$
for~$i=1,\ldots,n$.
\item[(iii)] $[\Im\Om^s]\cdot[\,\ov{\!X_k\!}\,]=0$ for all
$s\in\F$ and $k=1,\ldots,q$, where~$[\,\ov{\!X_k\!}\,]\in
H_m(M,\R)$.
\item[(iv)] $\displaystyle\sum_{\substack{1\le i\le n, \;
1\le j\le l_i: \\ k(i,j)=k}}\psi(x_i)^mZ(L_i)\cdot[\Si_i^j\,]=0$
for all $k=1,\ldots,q$.
\item[(v)] Let $N$ be the compact $m$-manifold obtained
by gluing $L_i$ into $X'$ at $x_i$ for $i=1,\ldots,n$ in the
obvious way. Suppose $N$ is connected. A sufficient
condition for this to hold is that $X$ and $L_i$ for
$i=1,\ldots,n$ are connected.
\end{itemize}

Let $g,\psi,R,B_R,\Si_i,\iota_i,\Up_i,\ze,U_\sCi,\Phi_{\sCi},
R',K,\phi_i,\eta_i,\eta_i^1,\eta_i^2,S_i,U_{\sXp},\Phi_{\sXp},
\allowbreak
A_i$,
\allowbreak
$\ga_i$,
\allowbreak
$\pi_i,
\allowbreak
T,
\allowbreak
K_i,
\vp_i,
\chi_i,\chi_i^1,\chi_i^2,U_\sLi,\Phi_\sLi,E_i$ and $\la$
be as in Definition \ref{co6def1}. Define
\begin{itemize}
\setlength{\parsep}{0pt}
\setlength{\itemsep}{0pt}
\item Let $V\cong H^2_{\rm cs}(X',\R)$ be a vector space
of smooth closed $2$-forms on $X'$ supported in $K$
representing~$H^2_{\rm cs}(X',\R)$.
\item Let $Y(L_i)\in H^1(\Si_i,\R)$ be as in Definition
\ref{co4def2}. Then by Hodge theory there exists a unique
$\ga_i\in C^\iy(T^*\Si_i)$ with $\d\ga_i=\d^*\ga_i=0$
and $[\ga_i]=Y(L_i)\in H^1(\Si_i,\R)$. Let $\pi_i:\Si_i\t(0,R')
\ra\Si_i$ be the projection.
\item Define $\varpi\in H^2_{\rm cs}(X',\R)$ to be the
image of $\bigl(Y(L_1),\ldots,Y(L_n)\bigr)$ under the
map $\bigoplus_{i=1}^nH^1(\Si_i,\R)\ra H^2_{\rm cs}(X',\R)$
in \eq{co3eq3}. Let $\be\in V$ be the unique element
with~$[\be]=\varpi$.
\item Choose $\nu_i\in(0,\mu_i-2)$ with $(0,\nu_i]\cap\D_\sSii
=\emptyset$ for~$i=1,\ldots,n$.
\item Apply part (b) of Theorem \ref{co3thm1} with $\ga_i$
and $\be$ as above. This gives $\al\in C^\iy(T^*X')$ with
$\d\al=\be$, $\d^*(\psi^m\al)=0$ and $\md{\na^j\al}=O(\rho^{-1-j})$
for $k\ge 0$, and $T_i\in C^\iy\bigl(\Si_i\t(0,R')\bigr)$
for $i=1,\ldots,n$ with
\end{itemize}
\ea
(\Up_i\circ\phi_i)^*(\al)&=\pi_i^*(\ga_i)+\d T_i
\quad \text{on $\Si_i\t(0,R')$ for $i=1,\ldots,n$, and}
\label{co8eq1}\\
\na^jT_i(\si,r)&=O(r^{\nu_i-j})\quad
\text{as $r\ra 0$, for all $j\ge 0$.}
\label{co8eq2}
\ea
\begin{itemize}
\item Apply Theorems \ref{co3thm7} and \ref{co3thm8} to $X$ and
$\bigl\{(M,J^s,\om^s,\Om^s):s\in\F\bigr\}$ with $V$ as above. This
gives an open set $\F'\subseteq\F$ with $0\in\F'$ and $\psi^s,
\up_i^s,\Up_i^s,\nu^s,\Phi_\sXp^s$ for $s\in\F'$ with $\nu^s\in V$
depending smoothly on $s$, satisfying \eq{co3eq14} and
\end{itemize}
\e
\begin{gathered}
\up_i^0=\up_i,\;\>
\Up_i^0=\Up_i,\;\>
\nu^0=0,\;\>
\Phi_\sXp^0=\Phi_\sXp,\;\>
(\up_i^s)^*(\Om)=\psi^s(x_i)^m\Om',\\
\Up_i^s(0)=x_i,\;\>
\d\Up_i^s\vert_0=\up_i^s,\;\>
(\Up_i^s)^*(\om^s)=\om',\;\>
(\Phi^s_\sXp)^*(\om^s)=\hat\om+\pi^*(\nu^s).
\end{gathered}
\label{co8eq3}
\e
\begin{itemize}
\item Define a function $f^s\in C^\iy(X')$ for $s\in\F'$ by
$(\Phi_\sXp^s)^*(\Im\Om^s)\vert_{X'}=f^s\,\d V$. As in \eq{co7eq12}, we
can show that $\na^jf^s=O\bigl(\rho^{1+j}\md{s}\bigr)$ for all $j\ge 0$,
where $\rho:X'\ra(0,1]$ is a {\it radius function}, as in Definition
\ref{co3def4}. Part (iii) above implies that $\int_{X'_k}f^s\,\d V=0$
for all $s\in\F'$ and~$k=1,\ldots,q$.
\item Using $\na^jf^s=O\bigl(\rho^{1+j}\md{s}\bigr)$ and
$\int_{X'_k}f^s\,\d V=0$ we can apply part (c) of Theorem
\ref{co3thm1} to $f^s$ on $X'_k$, for $s\in\F'$ and $k=1,\ldots,q$.
This gives exact 1-forms on $X'_k$ and functions on $\Si_i^j\t(0,R')$
for~$k(i,j)=k$.

Put these 1-forms and functions together for all $k$ to
give an exact 1-form $\varrho^s$ on $X'$ with
$\d^*(\psi^m\varrho^s_k)=f^s$ and functions $Z_i^s\in C^\iy
\bigl(\Si_i\t(0,R')\bigr)$ with
$(\Up_i\circ\phi_i)^*(\varrho^s)=\d Z_i^s$ on $\Si_i\t(0,R')$
for $i=1,\ldots,n$. As $\na^jf^s=O\bigl(\rho^{1+j}\md{s}\bigr)$,
we find that
\end{itemize}
\e
\md{\na^j\varrho^s}=O\bigl(\rho^{-1-j}\md{s}\bigr)
\quad\text{and}\quad
\na^jZ_i^s(\si,r)=O\bigl(r^{\nu_i-j}\md{s}\bigr)\quad
\text{for all $j\ge 0$.}
\label{co8eq4}
\e
\label{co8def1}
\end{dfn}

Observe that parts (iii) and (iv) above are equivalent to
requiring that both sides of \eq{co7eq39} are zero. To
prove a more general result, we would prefer to replace
(iii) and (iv) with the single condition \eq{co7eq39}.
However, part (iii) is necessary for the existence of
$\varrho^s$ and $Z_i^s$ above, and this together with
\eq{co7eq39} implies part (iv). We will discuss this
further after Theorem~\ref{co8thm3}.

Next we define $m$-submanifolds $N^{s,t}$ in $M$ for $s\in\F'$ and
$t\in(0,\de)$, adapting Definitions \ref{co6def2} and \ref{co7def2}.
It will turn out that $N^{s,t}$ is only {\it Lagrangian} in
$(M,\om^s)$ if $[\om^s]\cdot\iota_*(\ga)=t^2\varpi\cdot\ga$ for all
$\ga\in H_2(X,\R)$. However, we still define $N^{s,t}$ for all $s,t$,
as it will be useful that $N^{s,t}$ depends smoothly on~$s,t$.

\begin{dfn} We work in the situation of Definition \ref{co8def1}.
For $i=1,\ldots,n$, $s\in\F'$ and small $t>0$, define $P_i^{s,t}=
\Up_i^s(tK_i)$. Then $P_i^{s,t}$ is Lagrangian in $(M,\om^s)$, as
in Definition \ref{co7def2}. Let $F$ and $\tau$ be as in Definition
\ref{co6def2}. Modifying \eq{co6eq6}, for $i=1,\ldots,n$ and $t>0$
with $tT<t^\tau<2t^\tau<R'$, define a closed 1-form $\xi_i^{s,t}$
on $\Si_i\t(tT,R')$ by
\e
\begin{split}
\xi_i^{s,t}(\si,r)&=\d\bigl[F(t^{-\tau}r)A_i(\si,r)+
t^2(1-F(t^{-\tau}r))E_i(\si,t^{-1}r)\bigr]\\
&\quad+t^2\pi_i^*(\ga_i)+t^2\d\bigl[F(t^{-\tau}r)T_i(\si,r)\bigr]
+\d\bigl[F(t^{-\tau}r)Z_i^s(\si,r)\bigr]\\
&=F(t^{-\tau}r)\eta_i(\si,r)+
t^{-\tau}F'(t^{-\tau}r)A_i(\si,r)\d r\\
&\quad+t^2(1\!-\!F(t^{-\tau}r))\chi_i(\si,t^{-1}r)-
t^{2-\tau}F'(t^{-\tau}r)E_i(\si,t^{-1}r)\d r\\
&\quad+t^2F(t^{-\tau}r)(\Up_i\circ\phi_i)^*(\al)(\si,r)
+t^{2-\tau}F'(t^{-\tau}r)T_i(\si,r)\d r\\
&\quad
+F(t^{-\tau}r)(\Up_i\circ\phi_i)^*(\varrho^s)(\si,r)
+t^{-\tau}F'(t^{-\tau}r)Z_i^s(\si,r)\d r.
\end{split}
\label{co8eq5}
\e

Let $\xi_i^{1,s,t},\xi_i^{2,s,t}$ be the components of $\xi_i^{s,t}$ in
$T^*\Si$ and $\R$. Then when $r\ge 2t^\tau$ we have $\xi_i^{s,t}\equiv
\eta_i+t^2(\Up_i\circ\phi_i)^*(\al)+(\Up_i\circ\phi_i)^*(\varrho^s)$, and
when $r\le t^\tau$ we have $\xi_i^{s,t}(\si,r)=t^2\chi_i(\si,t^{-1}r)$.
Choose $\de\in(0,1]$ with $\de T\le\de^\tau<2\de^\tau\le R'$,
$\de K_i\!\subset\!B_R\!\subset\!\C^m$ and $\bmd{\xi_i^{s,t}
(\si,r)}<\ze r$ on $\Si_i\t(tT,R')$ for all $i=1,\ldots,n$,
$s\in\F'$ and $t\in(0,\de)$. Making $\F'$ smaller if necessary, this
is possible. Following \eq{co7eq5}, for $i=1,\ldots,n$, $s\in\F'$
and $t\in(0,\de)$ define $\Xi_i^{s,t}:\Si_i\t(tT,R')\ra M$ by
\e
\Xi_i^{s,t}(\si,r)=\Up_i^s\circ\Phi_\sCi\bigl(\si,r,
\xi_i^{1,s,t}(\si,r),\xi_i^{2,s,t}(\si,r)\bigr).
\label{co8eq6}
\e
Define $Q_i^{s,t}=\Xi_i^{s,t}\bigl(\Si_i\t(tT,R')\bigr)$ for
$i=1,\ldots,n$, $s\in\F'$ and $t\in(0,\de)$. As in Definition
\ref{co7def2}, we find that $Q_i^{s,t}$ is Lagrangian in~$(M,\om^s)$.

Modifying Definition \ref{co6def2}, let $\Ga(t^2\al+\varrho^s)$
be the graph of the 1-form $t^2\al+\varrho^s$ in $T^*X'$. Then
$\Ga(t^2\al+\varrho^s)\cap\pi^*(K)\subset T^*K$ is the graph of
$(t^2\al+\varrho^s)\vert_K$. By compactness of $K$, making $\F'$
and $\de>0$ smaller if necessary, we can suppose that $\Ga(t^2\al+
\varrho^s)\cap\pi^*(K)\subset U_\sXp$ for all $s\in\F'$ and
$t\in(0,\de)$. Define
\e
K^{s,t}=\Phi_\sXp^s\bigl(\Ga(t^2\al+\varrho^s)\cap\pi^*(K)\bigr)
\quad\text{for $s\in\F'$ and $t\in(0,\de)$.}
\label{co8eq7}
\e
Then $K^{s,t}$ is a submanifold of $M$ with boundary,
diffeomorphic to~$K$.

We calculate when $K^{s,t}$ is Lagrangian in $(M,\om^s)$. Write
$\iota^{s,t}$ for the natural diffeomorphism $K\ra K^{s,t}$ given
by $\Phi_\sXp^s\circ(t^2\al+\varrho^s)\vert_K$. Now
$(\Phi^s_\sXp)^*(\om^s)=\hat\om+\pi^*(\nu^s)$ by \eq{co8eq3}.
Pushing down by $\pi:\Ga(t^2\al+\varrho^s)\ra X'$ gives
\e
(\iota^{s,t})^*(\om^s)=\pi_*\bigl(\hat\om\vert_{\Ga(t^2\al+\varrho^s)}
\bigr)+\nu^s=-\d(t^2\al+\varrho^s)+\nu^s=-t^2\be+\nu^s,
\label{co8eq8}
\e
by some standard symplectic geometry, and as $\d\al=\be$ and
$\varrho^s$ is exact.

Hence $K^{s,t}$ is Lagrangian in $(M,\om^s)$ if and only if
$\nu^s=t^2\be$. But $\nu^s,\be\in V$ and $V\cong H^2_{\rm cs}
(X',\R)$ by Definition \ref{co8def1}, so $K^{s,t}$ is Lagrangian
if and only if $[\nu^s]=t^2\varpi$ in $H^2_{\rm cs}(X',\R)$, as
$[\be]=\varpi$. Now Theorem \ref{co3thm8} identifies $[\nu^s]$
as the unique class in $H^2_{\rm cs}(X',\R)$ with
$[\nu^s]\cdot\ga=[\om^s]\cdot\iota_*(\ga)$ for all $\ga\in
H_2(X,\R)$, where $\iota_*:H_2(X,\R)\ra H_2(M,\R)$ is induced
by the inclusion $\iota:X\ra M$. Therefore $K^{s,t}$ is
Lagrangian in $(M,\om^s)$ if and only if
\e
[\om^s]\cdot\iota_*(\ga)=t^2\varpi\cdot\ga
\quad\text{for all $\ga\in H_2(X,\R)$.}
\label{co8eq9}
\e

Define $\G\subseteq\F'\t(0,\de)$ to be the subset of $(s,t)$
satisfying \eq{co8eq9}. Then $K^{s,t}$ is Lagrangian in $(M,\om^s)$
if and only if $(s,t)\in\G$. For $s\in\F'$ and $t\in(0,\de)$,
define $N^{s,t}$ to be the disjoint union of $K^{s,t}$, $P_1^{s,t},
\ldots,P_n^{s,t}$ and $Q_1^{s,t},\ldots,Q_n^{s,t}$. Then $N^{s,t}$
is a compact, smooth $m$-submanifold of $M$ {\it without boundary},
which depends smoothly on $s,t$. The proof of this follows
Definitions \ref{co6def2} and \ref{co7def2}, with simple changes.
Also $N^{s,t}$ is {\it Lagrangian}\/ in $(M,\om^s)$ if and only
if~$(s,t)\in\G$.

Let $h^{s,t}$, $\d V^{s,t}$, ${\rm e}^{\smash{i\th^{s,t}}}$ and
$\ve^{s,t}=(\psi^s)^m\sin\th^{s,t}$ be as in Definition~\ref{co7def2}.
\label{co8def2}
\end{dfn}

Here equation \eq{co8eq9} is a weakening of the combination
of part (iii) of Definition \ref{co6def1} and part (ii)
of Definition \ref{co7def1}. For by exactness in \eq{co3eq3}
and the definition of $\varpi$, we see that part (iii) of
Definition \ref{co6def1} admits a solution
$\varrho$ if and only if $\varpi=0$. But if $\varpi=0$
then \eq{co8eq9} becomes part (ii) of Definition \ref{co7def1}.
Thus part (iii) of Definition \ref{co6def1} and part (ii) of
Definition \ref{co7def1} together imply \eq{co8eq9}, but not
vice versa.

It is not difficult to see that \eq{co8eq9} is a
{\it necessary condition} for the existence of a {\it Lagrangian}
$m$-fold $N^{s,t}$ in $(M,\om^s)$ made by gluing $tL_i$ into $X$
at $x_i$ for $i=1,\ldots,n$, just as part (iii) of Definition
\ref{co6def1} and part (ii) of Definition \ref{co7def1} were
necessary for the existence of Lagrangian $m$-folds $N^t$ and
$N^{s,t}$ in \S\ref{co6} and \S\ref{co7}. The definition shows that
\eq{co8eq9} is also a {\it sufficient condition} for small~$s,t$.

In \S\ref{co72} we saw that $\Im\Om^s\vert_{N^{s,t}}$ is
$O\bigl(\md{s}\bigr)$ on $K^s$. The condition $\md{s}\le
t^{\ka+m/2}$ in Theorem \ref{co7thm6} is there to control
these $O\bigl(\md{s}\bigr)$ error terms. However, in
this section we will need to allow $\md{s}=O(t^2)$, so
we cannot require $\md{s}\le t^{\ka+m/2}$. Therefore, an
$O\bigl(\md{s}\bigr)$ contribution to $\Im\Om^s\vert_{N^{s,t}}$
on $K^{s,t}$ is unacceptably large in this section.

This is the reason for using the exact 1-form $\varrho^s$
and functions $Z_i^s$ to construct $N^{s,t}$. They will
have the effect of cancelling out the $O\bigl(\md{s}\bigr)$
contributions, so that $\Im\Om^s\vert_{N^{s,t}}=O\bigl(\ms{s}
\bigr)$ on $K^{s,t}$. We could have used $\varrho^s$ and
$Z_i^s$ in the same way in \S\ref{co7}, and so relaxed the
condition $\md{s}\le t^{\ka+m/2}$ in Theorem~\ref{co7thm6}.

Here is the analogue of Definitions \ref{co6def3} and~\ref{co7def3}.

\begin{dfn} In the situation of Definitions \ref{co8def1} and
\ref{co8def2} define vector spaces $W^{s,t}\subset C^\iy(N^{s,t})$
for $s\in\F'$ and $t\in(0,\de)$, elements $w_{\bf d}^{s,t}\in
W^{s,t}$ for ${\bf d}\in\R^q$, and the projection $\pi_\sWst$,
exactly as in Definition \ref{co7def3}, except that in (i) we
replace $K^s$ by $K^{s,t}$, and define $w_{\bf d}^{s,t}\equiv d_k$
on the $k^{\rm th}$ connected component of $K^{s,t}$ diffeomorphic
to~$K\cap X'_k$.
\label{co8def3}
\end{dfn}

\subsection{Estimating $\Im\Om^s\vert_{N^{s,t}}$}
\label{co82}

We now bound $\ve^{s,t}$ on $N^{s,t}$, which yields estimates
of $\Im\Om^s\vert_{N^{s,t}}$ as $\Im\Om^s\vert_{N^{s,t}}=
\ve^{s,t}\,\d V^{s,t}$. As this still makes sense if $N^{s,t}$
is not Lagrangian, we do it for all $s\in\F'$ and $t\in(0,\de)$,
rather than for $(s,t)\in\G$. Then we can exploit the smooth
dependence of $N^{s,t}$ on $s,t$, including at $(0,0)$ in a
certain sense. Here is the analogue of Propositions
\ref{co6prop1} and~\ref{co7prop1}.

\begin{prop} In the situation above, making $\F'$ and\/ $\de>0$
smaller if necessary, there exists $C>0$ such that for all\/
$s\in\F'$ and\/ $t\in(0,\de)$ we have
\ea
\md{\ve^{s,t}}&\le C\ms{s}+Ct^4,\quad
\md{\d\ve^{s,t}}\le C\ms{s}+Ct^4 \quad\text{on $K^{s,t}$,}
\label{co8eq10}\\
\bmd{(\Xi_i^{s,t})^*(\ve^{s,t})}(\si,r)&\le\!
\begin{cases}
Cr, & r\!\in\!(tT,t^\tau],\\
\begin{aligned}
&\!C\md{s}t^{\tau(\nu_i-2)}\!+\!
Ct^{4-4\tau}\!+\!Ct^{\tau(\mu_i-2)}\\
&\!+Ct^{(1-\tau)(2-\la)}+Ct^{2+\tau(\nu_i-2)},
\end{aligned}
& r\!\in\!(t^\tau,2t^\tau),\\
C\ms{s}r^{-4}+Ct^4r^{-4}, & r\!\in\![2t^\tau,R'),
\end{cases}
\label{co8eq11}\\
\!\bmd{(\Xi_i^{s,t})^*(\d\ve^{s,t})}(\si,r)&\le\!
\begin{cases}
C, & r\!\in\!(tT,t^\tau],\\
\begin{aligned}
&\!C\md{s}t^{\tau(\nu_i-3)}\!+\!
Ct^{4-5\tau}\!+\!Ct^{\tau(\mu_i-3)}\\
&\!+\!Ct^{(1\!-\!\tau)(2\!-\!\la)\!-\!\tau}\!+\!Ct^{2+\tau(\nu_i-3)},
\end{aligned}
& r\!\in\!(t^\tau,2t^\tau),\\
C\ms{s}r^{-5}+Ct^4r^{-5}, & r\!\in\![2t^\tau,R'),
\end{cases}
\label{co8eq12}\\
\text{and}\quad
\md{\ve^{s,t}}&\le Ct,\quad
\md{\d\ve^{s,t}}\le C
\quad\text{on $P_i^{s,t}$ for all\/ $i=1,\ldots,n$.}
\label{co8eq13}
\ea
Here $\md{\,.\,}$ is computed using $h^{s,t}$
or~$(\Xi_i^{s,t})^*(h^{s,t})$.
\label{co8prop1}
\end{prop}

\begin{proof} The proof combines those of Propositions
\ref{co6prop1} and \ref{co7prop1}, so we will be brief. Identify
$K^{s,t}$ and $\Xi_i^{s,t}\bigl(\Si_i\t[2t^\tau,R')\bigr)$ with
the corresponding regions in $X'$ in the natural way, so that we
can regard $f^s$ as a function, $\al,\varrho^s$ as 1-forms and
$h$ as a metric on these regions in $N^{s,t}$ rather than $X'$.
Generalizing the proof of \eq{co6eq18} and \eq{co6eq19} we find
that on $K^{s,t}$ and $\Xi_i^{s,t}\bigl(\Si_i\t[2t^\tau,R')\bigr)$
we have
\e
\begin{split}
\ve^{s,t}&=f^s\!-\!\d_{s,t}^*\bigl((\psi^s)^m(t^2\al\!+\!\varrho^s)\bigr)
\!+\!O\bigl(\rho^{-2}\ms{t^2\al\!+\!\varrho^s}\bigr)
\!+\!O\bigl(\ms{t^2\na\al\!+\!\na\varrho^s}\bigr)\\
&=f^s-\d_{s,t}^*\bigl((\psi^s)^m(t^2\al+\varrho^s)\bigr)
+O\bigl(\rho^{-4}\ms{s}+\rho^{-4}t^4\bigr)
\end{split}
\label{co8eq14}
\e
when $\rho^{-1}\md{t^2\al+\varrho^s},\md{t^2\na\al+\na\varrho^s}$
are small. Here $\d^*_{s,t}$ is $\d^*$ calculated using $h^{s,t}$,
and we use $\md{\na^j\al}=O(\rho^{-1-j})$ and $\md{\na^j\varrho^s}
=O(\rho^{-1-j}\md{s})$ in the second line.

Now $h^{s,t}=h+O\bigl(t^2\rho^{-2}+\md{s}\rho^{-2}\bigr)$ and
$\na(h^{s,t}-h)=O\bigl(t^2\rho^{-3}+\md{s}\rho^{-3}\bigr)$ in
these regions, and $\psi^s=\psi+O\bigl(\md{s}\bigr)$, $\na(
\psi^s-\psi)=O\bigl(\rho^{-1}\md{s}\bigr)$, so we find that
\e
\begin{split}
\d_{s,t}^*\bigl((\psi^s)^m(t^2\al+\varrho^s)\bigr)&=
\d^*\bigl(\psi^m(t^2\al+\varrho^s)\bigr)
+O\bigl(\rho^{-3}(\md{s}+t^2)\md{t^2\al+\varrho^s}\bigr)\\
&\qquad\qquad
+O\bigl(\rho^{-2}(\md{s}+t^2)\md{\na(t^2\al+\varrho^s)}\bigr)\\
&=f^s+O\bigl(\rho^{-4}\ms{s}+\rho^{-4}t^4\bigr),
\end{split}
\label{co8eq15}
\e
estimating as in \eq{co8eq14}. Here $\d^*$ is calculated w.r.t.\
$h$, and we use $\d^*(\psi^m\al)=0$ and $\d^*(\psi^m\varrho^s)=f^s$
on $X'$. Combining \eq{co8eq14} and \eq{co8eq15} gives
\e
\ve^{s,t}=O\bigl(\rho^{-4}\ms{s}+\rho^{-4}t^4\bigr)
\quad\text{on $K^{s,t}$ and $\Xi_i^{s,t}\bigl(\Si_i
\t[2t^\tau,R')\bigr)$.}
\label{co8eq16}
\e
A similar proof for derivatives shows that
\e
\md{\d\ve^{s,t}}=O\bigl(\rho^{-5}\ms{s}+\rho^{-5}t^4\bigr)
\quad\text{on $K^{s,t}$ and $\Xi_i^{s,t}\bigl(\Si_i
\t[2t^\tau,R')\bigr)$.}
\label{co8eq17}
\e

As $\rho^{-1}=O(1)$ on $K^{s,t}$, equations \eq{co8eq16} and
\eq{co8eq17} imply \eq{co8eq10} for some $C>0$, making
$\F',\de$ smaller if necessary. Also, \eq{co8eq16} and
\eq{co8eq17} prove the bottom lines of \eq{co8eq11}
and \eq{co8eq12} for some $C>0$. On $P_i^{s,t}$ and
$\Xi_i^{s,t}\bigl(\Si_i\t(tT,t^\tau]\bigr)$ the definition
of $N^{s,t}$ coincides with Definition \ref{co7def2}. Thus,
equation \eq{co8eq13} and the top lines of \eq{co8eq11}
and \eq{co8eq12} follow from Proposition~\ref{co7prop1}.

The middle lines of \eq{co8eq11} and \eq{co8eq12} are
proved as in Proposition \ref{co6prop1}, adapted to the
families case as in Proposition \ref{co7prop1}. The additional
terms $C\md{s}t^{\tau(\nu_i-2)}$ and $C\md{s}t^{\tau(\nu_i-3)}$
which do not appear in \eq{co6eq10} and \eq{co6eq11} are there
to bound terms in $\ve^{s,t}$ coming from terms in $Z_i^s$ in
\eq{co8eq5}, in particular $t^{-\tau}F'(t^{-\tau}r)Z_i^s(\si,r)
\d r$ and its derivatives. This completes the proof.
\end{proof}

The purpose of including the $Z_i^s$ and $\varrho^s$ terms
in Definition \ref{co8def2} was to arrange for the $f^s$
terms in \eq{co8eq14} and \eq{co8eq15} to cancel, so
that the estimate \eq{co8eq16} has no $O\bigl(\md{s}\bigr)$
terms, but only $O\bigl(\ms{s}\bigr)$ terms. Propositions
\ref{co6prop2} and \ref{co7prop2} immediately generalize to give:

\begin{prop} For some $C'>0$ and all\/ $s\in\F'$ and\/
$t\in(0,\de)$ we have
\ea
\begin{split}
\lnm{\ve^{s,t}}{2m/(m+2)}&\le C't^{\tau(1+m/2)}
\Bigl(\ms{s}t^{-4\tau}+
t^{4-4\tau}+t^{(1-\tau)(2-\la)}\\
&\quad +\sum_{i=1}^n\bigl(t^{\tau(\mu_i-2)}+
t^{2+\tau(\nu_i-2)}+\md{s}t^{\tau(\nu_i-2)}\bigr)\Bigr),
\end{split}
\label{co8eq18}\\
\begin{split}
\cnm{\ve^{s,t}}{0}&\le
C'\Bigl(\ms{s}t^{-4\tau}+
t^{4-4\tau}+t^{(1-\tau)(2-\la)}\\
&\quad +\sum_{i=1}^n\bigl(t^{\tau(\mu_i-2)}+t^{2+\tau(\nu_i-2)}
+\md{s}t^{\tau(\nu_i-2)}\bigr)\Bigr),
\end{split}
\label{co8eq19}\\
\begin{split}
\text{and}\quad
\lnm{\d\ve^{s,t}}{2m}&\le C't^{-\tau/2}\Bigl(\ms{s}t^{-4\tau}+
t^{4-4\tau}+t^{(1-\tau)(2-\la)}\\
&\quad+\sum_{i=1}^n\bigl(t^{\tau(\mu_i-2)}+t^{2+\tau(\nu_i-2)}
+\md{s}t^{\tau(\nu_i-2)}\bigr)\Bigr),
\end{split}
\label{co8eq20}
\ea
computing norms with respect to the metric $h^{s,t}$ on~$N^{s,t}$.
\label{co8prop2}
\end{prop}

Note that putting $s=0$ in Propositions \ref{co8prop1} and
\ref{co8prop2} gives Propositions \ref{co6prop1} and \ref{co6prop2}.
Parts (iii) and (iv) of Definition \ref{co8def1} imply that
\eq{co7eq27} holds, and therefore Proposition \ref{co7prop5}
generalizes to give:

\begin{prop} Making $\F',\de$ smaller if necessary, there exists
$C''>0$ such that for all\/ $s\in\F'$ and\/ $t\in(0,\de)$ we
have~$\lnm{\pi_\sWst(\ve^{s,t})}{1}\le C''t^{(m+1)\tau}+C''\md{s}t^m$.
\label{co8prop3}
\end{prop}

Here is the analogue of Theorems \ref{co6thm1} and~\ref{co7thm1}.

\begin{thm} Making $\F'$ and\/ $\de>0$ smaller if necessary, there
exist\/ $A_2>0$, $\ka>1$ and\/ $\vartheta\in(0,2)$ such that for all\/
$s\in\F'$ and\/ $t\in(0,\de)$ with\/ $\md{s}\le t^\vartheta$, the
functions $\ve^{s,t}=(\psi^s)^m\sin\th^{s,t}$ on $N^{s,t}$ satisfy
$\lnm{\ve^{s,t}}{2m/(m+2)}\le A_2t^{\ka+m/2}$, $\cnm{\ve^{s,t}}{0}\le
A_2t^{\ka-1}$, $\lnm{\d\ve^{s,t}}{2m}\le A_2t^{\ka-3/2}$ and\/
$\lnm{\pi_\sWst(\ve^{s,t})}{1}\le A_2t^{\ka+m-1}$, as in part\/
{\rm(i)} of Theorem~\ref{co5thm1}.
\label{co8thm1}
\end{thm}

\begin{proof} Choose $\ka>1$ as in the proof of Theorem
\ref{co6thm1}, but requiring that strict inequality hold in
\eq{co6eq29} and \eq{co6eq30}. This is clearly possible. Then
all the terms in \eq{co8eq18}, \eq{co8eq19}, \eq{co8eq20} not
involving $\md{s}$ are bounded by multiples of $t^{\ka+m/2}$,
$t^{\ka-1}$, $t^{\ka-3/2}$ respectively, as in the proof of
Theorem~\ref{co6thm1}.

For the term $C't^{\tau(1+m/2)}\cdot\ms{s}t^{-4\tau}$ in
\eq{co8eq18} to be bounded by a multiple of $t^{\ka+m/2}$ when
$\md{s}\le t^\vartheta$, it is enough that $t^{\tau(1+m/2)}\cdot
t^{2\vartheta}\cdot t^{-4\tau}\le t^{\ka+m/2}$. Since $t\in(0,1)$,
this holds if $\tau(1+m/2)+2\vartheta-4\tau\ge\ka+m/2$. In the
same way, for all the terms in \eq{co8eq18}, \eq{co8eq19},
\eq{co8eq20} involving $\md{s}$ to be bounded by multiples of
$t^{\ka+m/2}$, $t^{\ka-1}$, $t^{\ka-3/2}$ respectively, it is
enough that for $i=1,\ldots,n$ we have
\ea
\tau(1\!+\!m/2)\!+\!2\vartheta\!-\!4\tau&\ge\!\ka\!+\!m/2,&\;
\tau(1\!+\!m/2)\!+\!\vartheta\!+\!\tau(\nu_i\!-\!2)&\ge\!\ka\!+\!m/2,
\label{co8eq21}\\
2\vartheta-4\tau &\ge\ka-1,&\;
\vartheta+\tau(\nu_i-2)&\ge\ka-1,
\label{co8eq22}\\
-\tau/2+2\vartheta-4\tau &\ge\ka-3/2,&\;
-\tau/2+\vartheta+\tau(\nu_i-2)&\ge\ka-3/2.
\label{co8eq23}
\ea

As $\tau<1$, equations \eq{co8eq22} and \eq{co8eq23} follow from
\eq{co8eq21}. Thus, we need to choose $\vartheta\in(0,2)$ such
that \eq{co8eq21} holds for $i=1,\ldots,n$. This is possible as
$\ka$ was defined to satisfy the first equation of \eq{co6eq29}
and the second equation of \eq{co6eq30} with strict inequalities.
With these values of $\ka$ and $\vartheta$, the theorem follows
from Propositions \ref{co8prop2} and \ref{co8prop3}
with~$A_2=\max\bigl(3(n+1)C',2C''\bigr)$.
\end{proof}

The inequality $\vartheta<2$ in Theorem \ref{co8thm1} is
important for the following reason. If $\varpi\ne 0$ then
\eq{co8eq9} implies that $t^2\le D\md{s}$ for some $D>0$ and all
small $(s,t)\in\G$. Combining this with $\md{s}\le t^\vartheta$
gives $\md{s}\le D^{\vartheta/2}\md{s}^{\vartheta/2}$. If
$\vartheta>2$, or $\vartheta=2$ and $D<1$, then this fails for
all small nonzero~$s$.

That is, if we had $\vartheta\ge 2$ in Theorem \ref{co8thm1}
and $\varpi\ne 0$ then the condition $\md{s}\le t^\vartheta$
would have {\it excluded\/} all small $(s,t)\in\G$, so the
construction would produce {\it no} compact SL $m$-folds
$\smash{\ti N^{s,t}}$. But as $\vartheta<2$ the condition
$\md{s}\le t^\vartheta$ includes all small $(s,t)\in\G$ with
$\md{s}=O(t^2)$, which will in general be a nonempty subset
of~$\G$.

We may extend Definition \ref{co7def4} and Theorems
\ref{co7thm2}--\ref{co7thm4} to the situation of this
section in a straightforward way, by including the
modifications introduced in Definition \ref{co6def4}
and Theorems \ref{co6thm2}--\ref{co6thm4}, and there
are no significant new issues. Thus we prove:

\begin{thm} Theorems \ref{co7thm2}--\ref{co7thm4} hold in
the situation of Definitions \ref{co8def1}--\ref{co8def3}.
\label{co8thm2}
\end{thm}

\subsection{The main result for families when $Y(L_i)\ne 0$}
\label{co83}

Here is our main result, the analogue of Theorem \ref{co6thm5}
for families.

\begin{thm} Let\/ $(M,J,\om,\Om)$ be an almost Calabi--Yau
$m$-fold for $2\!<\!m\!<\nobreak\!6$, and\/ $X$ a compact SL\/
$m$-fold in $M$ with conical singularities at\/ $x_1,\ldots,x_n$
and cones $C_1,\ldots,C_n$. Define $\psi:M\ra(0,\iy)$ as in
\eq{co2eq3}. Let\/ $L_1,\ldots,L_n$ be Asymptotically Conical
SL\/ $m$-folds in $\C^m$ with cones $C_1,\ldots,C_n$ and rates
$\la_1,\ldots,\la_n$. Suppose $\la_i\le 0$ for $i=1,\ldots,n$.
Write $X'=X\sm\{x_1,\ldots,x_n\}$ and\/~$\Si_i=C_i\cap{\cal S}^{2m-1}$.

Set\/ $q=b^0(X')$, and let\/ $X_1',\ldots,X_q'$ be the connected
components of\/ $X'$. For $i=1,\ldots,n$ let\/ $l_i=b^0(\Si_i)$,
and let\/ $\Si_i^1,\ldots,\Si_i^{\smash{l_i}}$ be the connected
components of\/ $\Si_i$. Define $k(i,j)=1,\ldots,q$ by $\Up_i
\circ\vp_i\bigl(\Si_i^j\t(0,R')\bigr)\subset X'_{\smash{k(i,j)}}$
for $i=1,\ldots,n$ and\/ $j=1,\ldots,l_i$. Suppose that
\e
\sum_{\substack{1\le i\le n, \; 1\le j\le l_i: \\
k(i,j)=k}}\psi(x_i)^mZ(L_i)\cdot[\Si_i^j\,]=0
\quad\text{for all\/ $k=1,\ldots,q$.}
\label{co8eq24}
\e
Suppose also that the compact\/ $m$-manifold\/ $N$ obtained
by gluing $L_i$ into $X'$ at\/ $x_i$ for $i=1,\ldots,n$ is
connected. A sufficient condition for this to hold is that\/
$X$ and\/ $L_i$ for $i=1,\ldots,n$ are connected.

Suppose $\bigl\{(M,J^s,\om^s,\Om^s):s\in\F\bigr\}$ is a smooth
family of deformations of\/ $(M,J,\om,\Om)$, with base space
$\F\subset\R^d$, satisfying
\e
[\Im\Om^s]\cdot[\,\ov{\!X'_k\!}\,]=0
\quad\text{for all\/ $s\in\F$ and\/ $k=1,\ldots,q$.}
\label{co8eq25}
\e
Define $\varpi\in H^2_{\rm cs}(X',\R)$ to be the image
of\/ $\bigl(Y(L_1),\ldots,Y(L_n)\bigr)$ under the map
$\bigoplus_{i=1}^nH^1(\Si_i,\R)\ra H^2_{\rm cs}(X',\R)$
in \eq{co3eq3}. Define $\G\subseteq\F\t(0,1)$ to be
\e
\G=\bigl\{(s,t)\in\F\t(0,1):
[\om^s]\cdot\iota_*(\ga)=t^2\varpi\cdot\ga
\;\>\text{for all $\ga\in H_2(X,\R)$}\bigr\},
\label{co8eq26}
\e
where $\iota_*:H_2(X,\R)\ra H_2(M,\R)$ is the natural inclusion.

Then there exist\/ $\ep\in(0,1)$, $\ka>1$ and\/ $\vartheta\in(0,2)$
and a smooth family
\e
\bigl\{\smash{\ti N^{s,t}}:(s,t)\in\G,\quad t\in(0,\ep],
\quad \md{s}\le t^\vartheta\bigr\},
\label{co8eq27}
\e
such that\/ $\smash{\ti N^{s,t}}$ is a compact, nonsingular SL\/
$m$-fold in $(M,J^s,\om^s,\Om^s)$ diffeomorphic to $N$, which is
constructed by gluing $tL_i$ into $X$ at\/ $x_i$ for $i=1,\ldots,n$.
In the sense of currents in Geometric Measure Theory,
$\smash{\ti N^{s,t}}\ra X$ as~$s,t\ra 0$.
\label{co8thm3}
\end{thm}

The proof follows that of Theorem \ref{co7thm6}, but using
Theorems \ref{co8thm1} and \ref{co8thm2} to prove parts
(i)--(vii) of Theorem \ref{co5thm1} for $N^{s,t}$ and
$W^{s,t}$, for $(s,t)\in\G$. Here the restriction to
$(s,t)\in\G$ is because $N^{s,t}$ is only Lagrangian
in $(M,\om^s)$ if \eq{co8eq9} holds, from Definition
\ref{co8def2}, and Theorem \ref{co5thm1} applies only
when $N^{s,t}$ is Lagrangian.

When $\varpi=0$ the SL $m$-fold $\smash{\ti N^{0,t}}$ above
coincides with the SL $m$-fold $\smash{\ti N^t}$ of Theorem
\ref{co6thm5}. More generally, as in \S\ref{co75}, when
$\varpi=0$ much of Theorem \ref{co8thm3} follows from
Theorems \ref{co2thm2} and \ref{co6thm5}. However, when
$\varpi\ne 0$ there exist {\it no} desingularizations
$\smash{\ti N^{0,t}}$ of $X$ in $(M,J,\om,\Om)$, but
Theorem \ref{co8thm3} can still give desingularizations
$\smash{\ti N^{s,t}}$ in $(M,J^s,\om^s,\Om^s)$ for $s\ne 0$.
That is, by deforming $(M,J,\om,\Om)$ to $(M,J^s,\om^s,\Om^s)$
we can overcome the topological obstructions to desingularizing
$X$ in~$(M,J,\om,\Om)$.

Theorem \ref{co8thm3} is essentially a combination of
Theorems \ref{co6thm5} and \ref{co7thm5}. However,
Theorem \ref{co7thm5} assumes the $[\Im\Om^s]$ and
$Z(L_i)$ satisfy \eq{co7eq39}, whereas Theorem
\ref{co8thm3} assumes the stronger \eq{co8eq24} and
\eq{co8eq25}, which together imply \eq{co7eq39}. It
is natural to ask whether Theorem \ref{co8thm3} still
holds with the weaker assumption \eq{co7eq39} in place
of \eq{co8eq24} and~\eq{co8eq25}.

The answer to this ought really to be yes, as \eq{co8eq25}
was a technical condition introduced to ensure that $\varrho^s$
and $Z_i^s$ exist in Definition \ref{co8def1}, and with some
more work one should to be able to do without it. However,
such a revised theorem would suffer from the following problem.

Suppose $s,t$ satisfy both \eq{co7eq39} and equation
$[\om^s]\cdot\iota_*(\ga)=t^2\varpi\cdot\ga$ in \eq{co8eq26}.
If $\varpi\ne 0$ then \eq{co8eq26} implies that $t^2=O\bigl(
\md{s}\bigr)$, and we expect $\md{s}\approx t^2$. If the
right hand side of \eq{co7eq39} is nonzero then it gives
$t^m=O\bigl(\md{s}\bigr)$, and we expect $\md{s}\approx t^m$.
These conditions are not compatible, as $m>2$. The two kinds
of obstruction need to be resolved at different length scales.

Actually this is not a serious problem, provided the family
$\F$ has large enough dimension. There could still exist a
family of solutions $(s,t)$ to \eq{co7eq39} and \eq{co8eq26}
with $\md{s}\approx t^2$ and $[\Im\Om^s]\cdot[\,\ov{\!X'_k\!}\,]
=O(\md{s}^{m/2})=O(t^m)$ for~$k=1,\ldots,q$.

If $X'$ is connected then \eq{co8eq24} holds automatically
as in \S\ref{co65}, and \eq{co8eq25} simplifies to
$[\Im\Om^s]\cdot[X]=0$, giving an analogue of Theorem
\ref{co6thm6} for families:

\begin{thm} Let\/ $(M,J,\om,\Om)$ be an almost Calabi--Yau
$m$-fold for $2\!<\!m\!<\nobreak\!6$, and\/ $X$ a compact SL\/
$m$-fold in $M$ with conical singularities at\/ $x_1,\ldots,x_n$
and cones $C_1,\ldots,C_n$. Let\/ $L_1,\ldots,L_n$ be
Asymptotically Conical SL\/ $m$-folds in $\C^m$ with cones
$C_1,\ldots,C_n$ and rates $\la_1,\ldots,\la_n$. Suppose 
$\la_i\le 0$ for $i=1,\ldots,n$, and\/ $X'=X\sm\{x_1,\ldots,
x_n\}$ is connected.

Suppose $\bigl\{(M,J^s,\om^s,\Om^s):s\in\F\bigr\}$ is a smooth
family of deformations of\/ $(M,J,\om,\Om)$, with base space
$\F\subset\R^d$, satisfying
\e
[\Im\Om^s]\cdot[X]=0
\quad\text{for all\/ $s\in\F$, where $[X]\in H_m(M,\R)$.}
\label{co8eq28}
\e
Define $\varpi\in H^2_{\rm cs}(X',\R)$ to be the image
of\/ $\bigl(Y(L_1),\ldots,Y(L_n)\bigr)$ under the map
$\bigoplus_{i=1}^nH^1(\Si_i,\R)\ra H^2_{\rm cs}(X',\R)$
in \eq{co3eq3}. Define $\G\subseteq\F\t(0,1)$ to be
\e
\G=\bigl\{(s,t)\in\F\t(0,1):
[\om^s]\cdot\iota_*(\ga)=t^2\varpi\cdot\ga
\;\>\text{for all\/ $\ga\in H_2(X,\R)$}\bigr\},
\label{co8eq29}
\e
where $\iota_*:H_2(X,\R)\ra H_2(M,\R)$ is the natural inclusion.

Then there exist\/ $\ep\in(0,1)$, $\ka>1$ and\/ $\vartheta\in(0,2)$
and a smooth family
\e
\bigl\{\smash{\ti N^{s,t}}:(s,t)\in\G,\quad t\in(0,\ep],
\quad \md{s}\le t^\vartheta\bigr\},
\label{co8eq30}
\e
such that\/ $\smash{\ti N^{s,t}}$ is a compact, nonsingular
SL\/ $m$-fold in $(M,J^s,\om^s,\Om^s)$, which is constructed
by gluing $tL_i$ into $X$ at\/ $x_i$ for $i=1,\ldots,n$. In
the sense of currents in Geometric Measure Theory,
$\smash{\ti N^{s,t}}\ra X$ as~$s,t\ra 0$.
\label{co8thm4}
\end{thm}

When $d=0$ and $\F=\{0\}=\R^d$ equations \eq{co8eq25} and
\eq{co8eq28} hold automatically, and $\G$ in \eq{co8eq26}
and \eq{co8eq29} is nonempty if and only if $\varpi=0$.
But by exactness in \eq{co3eq3}, $\varpi=0$ is the
necessary and sufficient condition for the existence
of $\varrho$ in Theorems \ref{co6thm5} and \ref{co6thm6}.
Therefore, when $d=0$ Theorems \ref{co8thm3} and
\ref{co8thm4} reduce to Theorems \ref{co6thm5} and
\ref{co6thm6} respectively.

\end{document}